\documentclass[11pt]{amsart}
\usepackage{frenchineq}
\usepackage{hyperref}
\usepackage{cleveref}
\usepackage{fullpage}
\usepackage{fixmath}
\usepackage{mathtools}
\usepackage{comment}
\usepackage[english]{babel}
\usepackage[T1]{fontenc}
\usepackage{xargs}
\usepackage{amscd}
\usepackage{amsmath}
\usepackage{amsfonts}
\usepackage{amssymb}

\usepackage[colorinlistoftodos,bordercolor=orange,backgroundcolor=orange!20,linecolor=orange,textsize=scriptsize,disable]{todonotes}
\providecommand{\doi}[1]{Eprint \href{http://dx.doi.org/#1}{doi:#1}}
      
\usepackage{algorithm}
\usepackage{algpseudocode}

\usepackage{pgf} \usepackage{tikz}
\usetikzlibrary{calc} \usetikzlibrary{arrows}
\usetikzlibrary{shapes} \usepackage{barycentric}
\usetikzlibrary{shapes}
\usetikzlibrary{plotmarks}
\newcommand{\Minimize}{\operatornamewithlimits{Min}}

\newcommand{\supp}{\operatorname{supp}}
\newcommand{\proj}{\mathbb{P}}

\newcommand{\R}{\mathbb{R}}

\newcommand{\rmax}{\mathbb{R}_{\max}}
\newcommand{\zmax}{\mathbb{Z}_{\max}}
\newcommand{\rmin}{\mathbb{R}_{\min}}

\newcommand{\nrows}{m}
\newcommand{\ncols}{n}

\newtheorem{theorem}{Theorem}%
\newtheorem{corollary}[theorem]{Corollary}
\newtheorem{lemma}[theorem]{Lemma}
\newtheorem{proposition}[theorem]{Proposition}

\newtheorem{problem}{Problem}

\newcommand{\cW}{\mathcal{W}}
\newcommand{\cC}{\mathcal{C}}
\theoremstyle{remark}
\newtheorem{remark}{Remark}

\newtheorem{assumption}{Assumption}

\newcommand{\Col}{\operatorname{Col}}
\newcommand{\Row}{\operatorname{Row}}

\newcommand{\cw}{\operatorname{cw}}
\newcommand{\cV}{\mathcal{V}}
\newcommand{\cH}{\mathcal{H}}
\newcommand{\I}{\operatorname{I}}

\newcommand{\dist}{\operatorname{dist}}

\newcommand{\argmax}{\operatorname{argmax}}
\newcommand{\Sp}{\operatorname{Sp}}
\newcommand{\inr}{\operatorname{in}}
\newcommand{\inradius}{\operatorname{in-rad}}
\newcommand{\inrv}[1]{r^{\inr}_{#1}} %
\newcommand{\NEW}[1]{{\em #1}}

\let\ifX\iffalse
\newcommand{\zero}{\mathbf{0}}
\newcommand{\maxset}{S^{\max}} %
\newcommand{\minset}{S^{\min}} %
\newcommand{\minsetprime}{S^{'\max}} %
\newcommand{\Ttype}{T^{\mathrm{ty}}}
\newcommand{\Ttypei}{T^{\mathrm{ty},i}}
\def\coord#1#2#3{{-sqrt(3)/2*(#1-#2)} ,{ -(1/2)*#1 - (1/2)*#2 + #3}}

\usepackage{graphicx}
\usepackage{subfigure}
\usepackage[labelformat=parens,labelsep=quad,skip=3pt]{caption}

\title{Tropical linear regression and mean payoff games: or, how to measure the distance to equilibria}
\author{Marianne Akian}
\author{St\'ephane Gaubert}
\author{Yang Qi}
\author{Omar Saadi}
\address{INRIA and CMAP, \'Ecole polytechnique, IP Paris, CNRS}
\email{Marianne.Akian@inria.fr}
\email{Stephane.Gaubert@inria.fr}
\email{Yang.Qi@inria.fr}
\email{Omar.Saadi@polytechnique.edu}
\thanks{O. Saadi acknowledges the support of the Hassan II Academy of Science and Technology. We acknowledge the support of Fondation Math\'ematique Jacques Hadamard and EDF through the Gaspard Monge Program for Operations Research and their Interactions with Data Sciences (PGMO)}

\catcode`\@=11
\@namedef{subjclassname@2020}{%
  \textup{2020} Mathematics Subject Classification}
\catcode`\@=12
\begin{document}
\maketitle
\begin{abstract}
  We study a tropical linear regression problem
  consisting in finding the best approximation of a set of points by a tropical hyperplane.
  We establish a strong duality theorem, showing that the value of this
  problem coincides with the maximal radius of a Hilbert's ball
  included in a tropical polyhedron. We also show that this
  regression problem is polynomial-time equivalent to mean payoff games.
  We illustrate our results by solving an inverse problem
  from auction theory. 
  In this setting, a tropical hyperplane represents
  the set of equilibrium prices. Tropical
  linear regression allows us to quantify the distance
  of a market to the set of equilibria, and infer
  secret preferences of a decision maker.
\end{abstract}
\section{Introduction}
\subsection{The tropical linear regression problem}
A {\em tropical hyperplane} in the $n$-dimensional tropical vector space $(\mathbb{R}\cup\{-\infty\})^n$ is a set
of vectors of the form
\begin{align}\label{e-def-troph}
\cH_a = \{x\in (\mathbb{R}\cup\{-\infty\})^n,\qquad
\max_{1\leq i\leq n} a_i + x_i \text{ is achieved at least twice} \} \enspace .
\end{align}
Such a hyperplane is parametrized by the vector $a=(a_1,\dots,a_n)\in (\R\cup\{-\infty\})^n$, which is required to be non-identically $-\infty$. 

Tropical hyperplanes are among the most basic objects in tropical
geometry. They are images by the valuation of hyperplanes
over non-archimedean fields, and so, they are the simplest examples
of tropical linear spaces~\cite{SpeyerSturmfels04,rinconfink}
and tropical hypersurfaces~\cite{kapranov}. Tropical
hyperplanes arise in tropical convexity~\cite{cgq02,DS04}, since
closed tropical convex sets can be described as intersections
of tropical half-spaces. A further motivation
arises from the study of pricing problems: tropical hypersurfaces
have been used in~\cite{Baldwin2019} to represent
the influence of prices on the decision of
agents buying bundles of elementary products. The ``unit demand'' case
(bundles of cardinality one) is modelled by tropical hyperplanes.

In this paper, we address the following tropical analogue of the linear
regression problem. Given a finite set of points $\cV\subset (\R\cup\{-\infty\})^n$,
we look for the best approximation of these points by a tropical hyperplane. 
Of course, the notion of ``best approximation'' depends on the metric.
A canonical choice in tropical geometry is the (additive version of)
Hilbert's projective metric. Its restriction to $\R^n$ is induced
by the so called {\em Hilbert's seminorm} or {\em Hopf oscillation}
\[
\|x\|_H:= \max_{i\in [n]} x_i -\min_{i\in [n]} x_i \enspace .
\]
It is a projective metric, in the sense that the distance between two points is zero
if and only if these two points differ by an additive constant.
Hence, we formulate the {\em tropical linear regression problem}
as the following optimization problem:
\begin{align}
  \operatornamewithlimits{Min}_{a}
  \max_{v\in \cV} \min_{x\in \cH_a \cap \R^n} \|v-x\|_H
  \label{def-tropreg}
\end{align}
where the minimum is taken over the space of parameters of tropical hyperplanes.
For simplicity, we assume for the moment that the vectors
$v\in \cV$ have finite entries, this assumption will be relaxed
in the body of the article.

\Cref{def-tropreg} is a non-convex optimization problem,
which is of a disjunctive nature since
a tropical hyperplane is a union of convex cones.

The tropical linear regression problem~\eqref{def-tropreg} is not only of theoretical interest. We shall see that it allows one to quantify
the ``distance to equilibrium'' of a market model,
and to infer hidden preferences of a decision maker.

\subsection{Results}
We show that tropical linear regression is tractable,
theoretically, and to some extent, computationally.
Our main result is a strong duality theorem,
\Cref{th-strong-duality}, showing that the infimum of the distance of the set of points $\cV$
to a tropical hyperplane coincides with the supremum of the radii of Hilbert's balls
included in the tropical convex cone generated by the elements
of $\cV$. This provides optimality certificates
which can be interpreted geometrically as collections of $n$ ``witness'' points
among the elements of $\cV$. Our approach also entails that
tropical linear regression is polynomial-time equivalent to
solving mean payoff games. The latter games,
originally studied in~\cite{ehrenfeucht,gurvich}, are among the problems in the complexity class NP $\cap$ co-NP~\cite{zwick}
for which no polynomial time algorithm is known.
However, several effective methods are available~\cite{gurvich,zwick,bjorklund,DG-06}. In particular, policy iteration allows to solve large scale instances~\cite{chaloupka}, even if it is generally super-polynomial~\cite{friedman}.
Thus, the present results lead to a practical
solution of the tropical linear regression problem.

We subsequently study variants
of the tropical linear regression problem, involving in particular
the {\em signed} notion of tropical hyperplane, obtained
by requiring the maximum in~\eqref{e-def-troph} to be achieved
by two indices $i,j$ belonging to prescribed disjoint subsets
$I,J$ of $[n]$. We also establish a strong duality theorem in this
setting, and provide reductions to mean payoff games for these
variants.

We finally illustrate tropical linear regression by an
application to an auction model. We consider a market governed
by an invitation to tender procedure. 
We suppose that a decision maker selects repeatedly bids
made by firms, based on the bid prices, which are ultimately
made public (after the decision is taken),
and also on other criteria
(assessments of the technical quality of each firm or of environmental
impact) or influence factors (like bribes).
This is a variant of the classical ``first-price sealed-bid auction''~\cite{krishna}, with a bias induced by the secret preference.
Here, we define the market to be at equilibria if for each invitation,
there are at least two best offers. Hence, in the simplest model (unit demand),
the set of equilibria prices can be
represented by a tropical hyperplane.  We distinguish two versions
of this problem, one in which only the prices are public,
and the other, in which the identities of the winners of the successive
invitations are also known. In both cases, we show that solving a tropical
linear regression problem allows an observer to quantify the distance
of such a market to equilibrium, and also to infer secret preference factors.
This solves, in the special case of unit-demand,
an inverse problem, consisting in identifying
the agent preferences and utilities in auction
models, like the one of~\cite{Baldwin2019}.
This might be of interest to a regulation authority wishing
to quantify anomalies, or to a bidder, who, seeing the history of the market,
would wish to determine how much he should have bidded to win a
given invitation or to get the best price for an invitation that he won,
thus avoiding the ``winner's curse''.

\subsection{Related work, and discussion}
Several ``best approximation'' problems have been studied in tropical
geometry. The simplest one consists in finding the nearest point
in a (closed) tropical module, in the sense of Hilbert's metric.
The solution is given by the tropical projection~\cite{cgq02},
see also~\cite{AGNS10}. The best approximation in the space of ultrametrics,
which is a fundamental example of tropical module in view of its application to phylogenetics, has been thoroughly studied~\cite{chepoi,Lin2017,bernstein}.
Another important special case is the best approximation of a point
by a tropical linear space~\cite{ardila,joswigyuberndt}. In contrast
with the regression problem studied here, these problems concern the
approximation of a {\em single} point.

It is a general
principle that regression (best approximation)
is somehow dual to separation. Hence, tropical linear regression
should be compared with the tropical support vector machines (SVM) introduced
in~\cite{gartnerJaggi}, and further studied in~\cite{tropicalSVM}. Whereas
the input of the tropical SVM problem (a configuration
of points in dimension $n-1$ partitioned in $n$ color classes) is the same as the
one of the version ``with types'' of the tropical linear regression problem,
we explain in~\Cref{rk-SVM} why both problems differ in essential ways.

A different problem of tropical regression
consists in finding a vector $x$ minimizing the sup-norm
$\|y-Ax\|_\infty$ where $y$ is a vector of observations,
and $A$ is a known matrix acting tropically on $x$. This can be
solved in (strongly)
polynomial time, again by a tropical projection~\cite{butkovic}.
See also~\cite{chepoi} for a general version
of this result. Tropical linear regression problems of this nature
have been studied in the context of learning~\cite{maragossurvey}.
The sparse version is of practical interest; it arose in the
approximation of solutions of Hamilton-Jacobi PDE, where
it was shown to be equivalent to a non-metric infinite dimensional
facility location problem~\cite{zheng1}. The finite dimension
version which is NP-hard is studied in~\cite{maragos}.

A different tropical regression problem, with a $L_1$-type error term
(instead of sup-norm here), has been solved in~\cite[Theorem~4]{yoshidatropPCAphylo}, in the special case of a configuration of $n$ points in dimension $n-1$.
The value is given by a {\em tropical volume}~\cite{depersin},
instead of an inner radius.

Tropical geometry has been applied to economics in~\cite{Baldwin2019},
see also~\cite{Tran2019}, and~\cite{Danilov2001} for early results
in this direction. Our modelling of agent's responses to prices
is inspired by~\cite{Baldwin2019}. Auction
models taking into account bribery have been studied in particular in~\cite{compte,bribery,Rachmilevitch2013}.

We build on the results of~\cite{AGGut10},
showing the equivalence between tropical linear
programming and mean payoff games. Further reductions
and equivalences, concerning in particular
the problem of the emptiness of
tropical linear prevarieties, were given in~\cite{podolskii}. The relation
between the mean payoff of a game and the inner
radius of a Shapley operator was first observed
in~\cite{skomra:tel-01958741,mtns}, where it was applied to
define a condition number and derive complexity results
for games.  The inner radius of tropical polyhedra
defined by $n$ generators in dimension $n-1$ was 
initially characterized in~\cite{Sergeev2007},
as a tropical eigenvalue.

Open problems related to the present work are discussed in the concluding
section.

\subsection{Organization}
In~\Cref{sec-prelim}, we recall the needed
results concerning tropical algebra, mean payoff games,
and non-linear Perron-Frobenius theory. In~\Cref{sec-inner},
we show that computing the inner radius of a tropical
polyhedron given by generators is equivalent
to solving a mean payoff game. \Cref{sec-strong} contains
our main results, including~\Cref{th-strong-duality},
the strong duality theorem for tropical linear
regression. Several variants of the tropical linear regression problem
are dealt with in~\Cref{sec-signed}. In~\Cref{sec-algo},
we explain how to solve tropical linear regression
problems in practice, using mean payoff games
algorithms. In~\Cref{sec-example}, we give an application
to an auction problem. The appendix provides sufficient
conditions for the existence of finite eigenvectors
of a class of Shapley operators.
These conditions are helpful when dealing with regression
problems for configurations of points with $-\infty$ coordinates.

 \section{Preliminaries}\label{sec-prelim}
\subsection{Tropical cones}
The max-plus semifield $\rmax$ is the set of real numbers, 
completed by $-\infty$ and equipped with the
addition $(a,b)\mapsto \max(a,b)$ and the multiplication $(a,b)\mapsto a \odot b := a+b$. 
The name ``tropical'' will be used in the sequel as a synonym of ``max-plus''.
We shall occasionally use variants of this semifield. These
include {\em the min-plus semifield} $\rmin$,  which is the set $\R\cup\{+\infty\}$,  equipped with the addition $(a,b)\mapsto \min(a,b)$ and the multiplication $(a,b)\mapsto a \odot b:= a+b$.
This semifield is isomorphic to $\rmax$. These also include
the subsemifield $\zmax\subset \rmax$, with ground set $\mathbb{Z}\cup \{-\infty\}$.
We refer the reader to~\cite{BCOQ92,butkovic,maclagan_sturmfels}
for background on tropical algebra.

For any integer $n$, we set $[n]:=\{1,\ldots, n\}$.
For all $x,y\in(\rmax)^n$, $A\in(\rmax)^{n\times m}$, and $\lambda\in \rmax$, $\lambda+x\in (\rmax)^n$ denotes the vector with entries $\lambda+x_i$, for $i\in [n]$, 
$\lambda+A\in (\rmax)^{n\times m}$ denotes the matrix with entries $\lambda+A_{ij}$, for $i\in [n], j \in [m]$,
$x\vee y=\sup (x,y)$ denotes the vector  with entries $\max(x_i,y_i)$, for $i\in [n]$, and $x\wedge y=\inf (x,y)$ denotes the vector  with entries $\min(x_i,y_i)$, for $i\in [n]$.
The set $(\rmax)^n$ equipped with the addition 
$(x,y)\mapsto x\vee y$ and the action $(\lambda,x)\mapsto \lambda+x$
of $\rmax$ is a tropical module, i.e.\  a module over the semifield $\rmax$.

A subset $\cC$ of $(\rmax)^n$ is a \NEW{tropical (convex) cone} or equivalently a
\NEW{tropical submodule}
of $(\rmax)^n$ if it satisfies
$x,y\in \cC$ and $\lambda\in\rmax$ implies $\lambda+x\in \cC$  and
$x \vee y\in \cC$.
We endow $(\rmax)^n$ with the topology defined by the metric
$\delta(x,y)=\max_{i\in [n]} |e^{x_i}-e^{y_i}|$.
It induces the usual topology in $\R^n$.
For any given subset $\cV$ of $(\rmax)^n$, we denote by
$\Sp (\cV)$ the tropical submodule of $(\rmax)^n$ 
generated by $\cV$, that is
the minimal  tropical submodule of $(\rmax)^n$ containing $\cV$.
A tropical polyhedral cone $\cC$ is a tropical cone which is 
finitely generated, that is such that there exists a finite subset $\mathcal V$ 
such that $\cC=\Sp (\cV)$.
For any given matrix $V$, we also denote by $\Col(V)$
the column space of $V$, that is the tropical polyhedral
cone generated by the columns of $V$, and we denote by $\Row(V)$ the row space of $V$, that is  the tropical polyhedral
cone generated by the rows of $V$.

A tropical polyhedral cone can also be defined externally
by a system of tropical linear inequalities of the form           
\begin{align}\label{sys-trop}
  \max_{j\in [\ncols]}(A_{ij}+x_j)
 \leq \max_{j\in[\ncols]} (B_{ij}+x_j),\qquad i \in [\nrows] \enspace ,
\end{align}
where $A_{ij},B_{ij}$ belong
to $\rmax$, see~\cite{SGKatz}. Then, $A$ and $B$ will be thought
of as $\nrows\times \ncols$ matrices with entries in $\rmax$.

Let $A\in (\rmax)^{m\times n}$
and $x\in (\rmax)^\ncols$. We denote by $Ax$ the vector in $(\rmax)^\nrows$ with entries $(Ax)_i=\max_{j\in [\ncols]}(A_{ij}+x_j)$, for $i \in [\nrows]$.
To a matrix $A\in (\rmax)^{m\times n}$, we associate the operator $A^\sharp : (\rmin)^\nrows \to (\rmin)^\ncols$, given by:
\[
\forall y \in (\rmin)^\nrows, \forall j \in [\ncols], (A^\sharp y)_j = \min_{i\in [\nrows]}(- A_{ij}+y_i) \enspace,
\]
with the convention $-\infty+\infty=+\infty$.
The operator $A^\sharp$ is called the {\em adjoint} of $A$ and we can easily check that it satisfies the following property:
\[
\forall x\in (\rmax)^\ncols, \forall y\in (\rmin)^\nrows, A x \leq y \Leftrightarrow x\leq A^\sharp y \enspace .
\]
We define the identity matrix $\I \in (\rmax)^{n\times n}$ by $\forall i \in[n], \I_{ii}=0$, and $\forall i,j\in [n], i\neq j, \I_{ij} = -\infty$.

A scalar $\mu$ is a {\em tropical eigenvalue} of a matrix
$M\in (\rmax)^{n\times n}$ if there exists a vector $u\in (\rmax)^n$,
not identically $-\infty$, such that $Mu=\mu + u$ in the tropical sense.
The eigenvalue is known to be unique when the digraph of $M$ is strongly
connected, then it coincides with
the maximum weight-to-length
ratio of the circuits of the digraph of $M$.
We denote it by $\lambda(M)$. See~\cite{BCOQ92,butkovic}
for more information.

\subsection{Mean payoff games}\label{sec-def-meanpayoff}

We consider zero-sum
deterministic games, with perfect information, defined
as follows. There are two players, ``Max'' and ``Min''
(the maximizer and the minimizer), who will move
a token on a weighted digraph.
We assume this digraph is finite and bipartite: the node
set is the disjoint union of two non-empty sets $\maxset$ and $\minset$,
and the arc set $\mathcal{A}$ is included in $(\maxset \times \minset)\cup (\minset\times \maxset)$. The set of states of the game is the set of
nodes of the digraph. We associate a real weight
$w_{rs}$ to each arc $(r,s)$.

The two players
alternate their actions.
When the token is in node $i\in \minset$,
Player Min must choose an arc $(i,j)$ in the digraph, meaning
he moves the token to node $j$, and pays $w_{ij}$
to player Max.
When Player Min has no possible action, that is, when there are no arcs of
the form $(i,j)$ in the digraph,
the game terminates, and Player Max receives $+\infty$.
Similarly, when the token is in node $j\in \maxset$, 
Player Max must choose an arc $(j,i)$ in the digraph. Then he moves the token from node $j$ to node $i$,
and receives $w_{ji}$ from Player Min.
When Player Max has no possible action, that is there are no arcs of
the form $(j,i)$ in the digraph, the game terminates,
and Player Max receives $-\infty$.

We measure the time in turns, i.e.,
a time step consists of two half-turns (a move made by Player Min
followed by a move made by Player Max). We consider the following game in horizon $k$: starting from an initial state $\bar{\imath}\in \minset$. 
the two players make $k$ moves each, unless
the game terminates before. So, if the game does not terminate
before time $k$, the {\em history} of the game is described by a sequence
of nodes $\bar{\imath}=i_0,j_1,i_1,\dots,j_k,i_k$, belonging alternatively
to $\minset$ and $\maxset$, 
and the total payment received by Player Max is given
by
\[
R^k_{\bar{\imath}}
= w_{i_0j_1}+ w_{j_1i_1}+w_{i_1j_2} + \cdots + w_{j_k i_k}
\enspace .
\]
If the game terminates before time $k$, we set $R^k_{\bar{\imath}}=\pm\infty$
depending on the player who had no available action.
The following assumption requires
Player Min to have at least one available action in every state:
\begin{assumption}\label{assump1} 
For all $i\in \minset$, there exists $j\in \maxset$ such that
$(i,j)$ is an arc of the digraph of the zero-sum deterministic game.
\end{assumption}
In this way, we always have $R^k_{\bar{\imath}} \in \R\cup \{-\infty\}$.
We shall also consider the dual assumption.
\begin{assumption}\label{assump2} 
For all $j\in \maxset$, there exists $i\in \minset$ such that
$(j,i)$ is an arc of the digraph of the zero-sum deterministic game.
\end{assumption}
In most works on mean payoff games, both assumptions
are required to hold, which entails in particular
that $R^k_{\bar{\imath}}$ is finite.
Here, we shall occasionally relax~\Cref{assump2},
but always require~\Cref{assump1}, so that
$R^k_{\bar{\imath}}\in \R\cup\{-\infty\}$.
This leads to an unpleasant symmetry
breaking. However, we shall see that this generality will be sometimes
needed to handle the application to tropical linear regression.
Indeed, from a tropical perspective, $-\infty$ is the zero element, hence
a meaningful value.

A strategy of a player is a map which associates to the history of the game
an action of this player.
Assuming that Player Min plays according to strategy $\sigma$,
and that Player Max plays according to strategy $\tau$,
we shall write $R^k_{\bar{\imath}}=R^k_{\bar{\imath}}(\sigma,\tau)$ to indicate the dependence on these strategies.
It follows from standard dynamic programming arguments
that the game in horizon $k$ starting from node $\bar{\imath}$ has a value $v^k_{\bar{\imath}}$ and that Players Min and Max have optimal strategies
$\sigma^*$ and $\tau^*$, respectively, see e.g.~\cite[Th.~IV.3.2]{SMZ}.
This means that the payment
function has the following saddle point property:
\[
R^k_{\bar{\imath}}(\sigma,\tau^*)
\leq v^k_{\bar{\imath}} = R^k_{\bar{\imath}}(\sigma^*,\tau^*) \leq R^k_{\bar{\imath}}(\sigma^*,\tau)
\]
for all strategies $\sigma,\tau$.
Moreover, the {\em value vector} $v^k:=(v^k_i)_{i  \in \minset}$ is determined by the following
dynamic programming equation
\[
v^k=T(v^{k-1}), \qquad v^0 =0 
\]
where $T: (\R\cup\{-\infty\})^n \to (\R\cup\{-\infty\})^n $ is the
{\em Shapley
operator}, defined, for $i\in \minset$, by
\begin{align}
T_i(x) = \min_{j, \; (i,j)\in \mathcal{A} } (w_{ij}+ \max_{l,\; (j,l)\in \mathcal{A}}(w_{jl} + x_l)) \enspace .\label{e-def-shapley0}
\end{align}
Owing to~\Cref{assump1}, the above minimum is never taken over an empty set,
whereas the above maximum is never taken over an empty set when~\Cref{assump2}
is made. By convention, the maximum of an empty set is $-\infty$.
When both assumptions hold, $T$ sends $\R^n$ to $\R^n$.

We are interested in the limit
\[
\chi(T):= \lim_{k\to \infty } T^k(0)/k = \lim_{k\to\infty} v^k/k \enspace .
\]
Thus, $\chi_i(T)$ yields the limit of the mean payoff per time unit, for
the game starting from the initial state $i$, when the horizon tends
to infinity. It follows from~\cite{kohlberg}
that the limit does exist, and that $\chi(T)\in \R^n$,
when Assumption~\ref{assump1} and Assumption~\ref{assump2} hold
(and more generally, when $T$ is a piecewise affine
self-map of $\R^n$ that is non-expansive in some norm,
see~\cite{AGGut10} and~\Cref{subsec-pf} below for details).
Alternatively, under the same assumptions,
$\chi_i(T)$ can be characterized as the value
of an infinite mean-payoff game, in which player Max
wishes to maximize the liminf of the average payment received
per time unit, whereas player Min wishes to minimize the liminf
of the same quantity --- this is the approach 
originally described by Ehrenfeucht and Mycielski~\cite{ehrenfeucht}.
It follows from a general result of Mertens and Neyman,
on the existence of the so called {\em uniform value}~\cite{mertens_neyman},
that this approach leads to the same notion of mean payoff.
Hence, we shall refer to $\chi_i(T)$ as the (asymptotic) mean payoff
starting from node $i$. 

More generally, the limit $\chi(T) \in (\R\cup\{-\infty\})^n$
does exist as soon as Assumption~\ref{assump1} is satisfied.
To see this, observe that we can always construct in polynomial
time an equivalent game satisfying
also Assumption~\ref{assump2}. Indeed,
let us delete any node $i$ of $\minset$ in which
Player Min has at least one action $(i,j)\in\mathcal{A}$
such that Player Max has no available action in state $j$.
After at most $|\minset|$ of such deletions,
we arrive at a new game, played on a bipartite subdigraph
of the original graph, induced by
a subset of nodes belonging to Player Min $\minsetprime\subset\minset$.
Note that $\minsetprime$ may be empty.
It is immediate that this subdigraph satisfies
both Assumption~\ref{assump1} and~\ref{assump2}.
So, for  $i\in \minset$, the existence of $\lim_k v^k_i/k$
follows from the result already established,
whereas for $i\in\minset\setminus\minsetprime$, we
have $v^k_i= -\infty$ for $k$ large enough, implying $\lim_k v^k_i/k=-\infty$.

A (stationary) {\em policy} of Player Min is a map $\sigma: \minset \to \maxset$ such that
$(i,\sigma(i))\in \mathcal{A}$ for all $i\in \minset$. Such a policy
determines a one-player game, in which Player Min
always selects moves $i\to \sigma(i)$. This one-player
game corresponds to the Shapley operator $T^\sigma$, defined by
\[
T^\sigma_i(x) =
w_{i\sigma(i)}+ \max_{l,\; (\sigma(i),l)\in \mathcal{A}}(w_{\sigma(i)l} + x_l) \enspace .
\]
Similarly, a policy of Player 
Max is a map  $\tau: \maxset \to \minset$ such that
$(j,\tau(j))\in \mathcal{A}$ for all $j\in \maxset$.
It determines a one-player game, with Shapley operator $^{\tau}T$ defined by
\[
^{\tau}T_i(x) = 
\min_{j, \; (i,j)\in \mathcal{A} } (w_{ij}+ w_{j\tau(j)} + x_{\tau(j)} ) \enspace .
\]
A result of Liggett and Lippman~\cite{liggett_lippman}
entails that each player has
optimal strategies in a mean payoff game,
which are obtained by applying a stationary policy.
This entails in particular that
\[
\chi(T) = \min_\sigma \chi(T^\sigma) = \max_{\tau} \chi(^{\tau}T) \enspace,
\]
under Assumptions~\ref{assump1} and~\ref{assump2}.
The mean payoff $\chi_i(T)$ is known to coincide
with the weight-to-length ratio of a circuit
of the bipartite digraph of the game, the ``length''
being measured as the number of full turns, i.e.,
as the number of Min nodes of the circuit
(one half of the ordinary length).
In particular, if the payments $w_{rs}$ are integers, the mean payoff 
is a rational number $p/q$, where $p,q$ are integers,
and $q$ is a positive integer bounded by the maximal
length of a circuit of the bipartite digraph of the game (measuring
the length as the number of nodes of Min that are visited)
and $|p/q|$ is bounded by $2|\max_{rs}w_{rs}|$.

Now we formalize the following problem.

\begin{problem}[Mean payoff games]\label{pb-1}
  Input: A finite bipartite directed graph with integer weights, satisfying Assumptions~\ref{assump1} and~\ref{assump2}, together with an initial node $\imath$.
  Question: Is the mean payoff $\chi_\imath(T)$ starting from node $\imath$ nonnegative?
\end{problem}
As discussed in the introduction, \Cref{pb-1} is a fundamental problem in algorithmic game theory~\cite{gurvich}. It belongs to the class NP $\cap$ coNP~\cite{zwick}, no polynomial time algorithm is known.

It will be useful to keep in mind several equivalent versions of this problem.

As a first variant, one may ask whether
$\chi_\imath(T)$ is positive, instead of non-negative. This variant
is equivalent to the negated version of~\Cref{pb-1}: considering
$\tilde{T}(x)=-T(-x)$, i.e., the Shapley operator of the game
in which all payments are negated, we have that $\chi(\tilde{T})=-\chi(T)$,
and so, $\chi_\imath(T)>0$ iff $\chi_\imath(\tilde{T})\leq 0$.

As observed above, the variant of mean payoff games
in which Assumption~\ref{assump2} is relaxed reduces
to the variant in which this assumption holds
by a preprocessing, so there is no restriction
on requiring Assumption~\ref{assump2} in~\Cref{pb-1}.

Another variant consists in computing $\chi_\imath(T)$, instead of deciding
whether $\chi_\imath(T)$ is nonnegative. This problem
of computation polytime Turing-reduces to~\Cref{pb-1}
by binary search.
Indeed, given a rational number $\alpha=p/q$,
we can consider the modified game with integer weights $w^\alpha_{rs}= 2q(w_{rs}-\alpha/2)$,
which corresponds to replacing the Shapley operator
$T$ by $T^\alpha:= 2q(-\alpha + T)$. Thus $\chi_\imath(T^\alpha)\geq 0$
iff $\chi_\imath(T)\geq \alpha$. Then, since the mean payoff $\chi_\imath(T)$
is a rational number whose absolute value is bounded
by $2\max_{rs}|w_{rs}|$ and whose denominator
is bounded by $|\minset|$, we can compute $\chi_\imath(T)$
by a dichotomy argument, calling at each step an oracle
solving~\Cref{pb-1} for a modified game with weights $w^\alpha_{rs}$.

\subsection{Perron-Frobenius tools}\label{subsec-pf}
We now recall some tools from Perron-Frobenius theory, in relation
with mean payoff games. We refer the reader to~\cite{AGGut10}
for more information.

We denote by $\bot$ the vector of $(\rmax)^n$ identically equal to $-\infty$.
We consider the {\em Hilbert's projective metric},
defined for vectors $x,y\in (\rmax)^n$ where at least
one of them is not equal to $\bot$, by
\[ d(x,y)=\inf\{\lambda-\mu\mid \lambda,\mu\in \R,\;
\mu+y_i\leq x_i\leq \lambda +y_i\; \forall i\in [n]\}\in \R_{\geq 0}\cup\{+\infty\}
\enspace .\]
In addition, we set $d(\bot,\bot):=0$.

The {\em support} of a vector $x\in (\rmax)^n$ is defined
by $\supp x:= \{i\in [n]\mid x_i \neq -\infty\}$.
Each subset $I\subset [n]$ yields
a {\em part} $P_I$ of $(\rmax)^n$,
consisting of vectors with support $I$.

Observe that $d(x,y)$ is finite if and only if $x$ and $y$ belong
to the same part $P_I$. Moreover, if $I\neq \emptyset$,
\[ d(x,y)=\max_{i\in I}(x_i-y_i)-\min_{i\in I}(x_i-y_i)\enspace .\]

We denote by $\proj(\rmax)^n$
the {\em tropical projective space}, i.e.,
the quotient of the set of non-identically $-\infty$ vectors of $(\rmax)^n$
by the equivalence relation $\sim$ which identifies tropically proportional vectors. We shall abuse notation and denote by the same symbol a vector
and its equivalence class.
Similarly, we shall think of a part $P_I$ with $I\neq \emptyset$ as a subset
of the tropical projective space. %

Observe that $d(x,y)$ vanishes if and only $x$ and $y$ represent
the same point of the tropical projective space,
so that $d$ yields a well defined
metric on each part of the tropical projective space.
We denote by $B(a,r)$ the closed ball centered at $a\in \R^n$ with radius $r$
under Hilbert's projective metric.

It will be convenient to consider
an abstract version of the concrete Shapley
operators used so far. 
We call (abstract) {\em Shapley operator} a map $T: (\rmax)^n\to (\rmax)^n$
that is order preserving, continuous, and such that $T(\alpha + x)=\alpha + T(x)$ for all $\alpha \in \rmax$ and $x\in(\rmax)^n$. Observe that the operator $T$
defined by~\eqref{e-def-shapley0}, with $\minset=[n]$, is a special
case of abstract Shapley operator, as soon as \Cref{assump1} holds.
We shall often consider situations in which an abstract
Shapley operator restricts to a map $\R^n\to \R^n$,
we will still use the term {\em Shapley operator} for the restricted
map. 

We are interested in the non-linear spectral problem for $T$,
consisting in finding a vector $u \in (\rmax)^n$, non-identically
$-\infty$, and a scalar $\lambda \in \rmax$ such that
$T(u) = \lambda + u$. The {\em spectral radius} of $T$ is defined
as 
 \begin{align}
   \label{e-def-rho}
   \rho(T)= \sup\{\lambda \in \R\cup \{-\infty\}\mid \exists u \in (\R\cup\{-\infty\})^n, u\neq \bot, \; T(u) = \lambda +u\} \enspace.\end{align} %
Variants of this spectral radius are given by the Collatz-Wielandt number $\cw$, 
 \begin{align}
   \label{e-def-cw}
   \cw(T)= \inf\{\lambda \in \R\mid \exists u \in \R^n, \; T(u) \leq \lambda +u\} \enspace.\end{align}
 and by the dual Collatz-Wielandt number
 \begin{align}
   \label{e-def-cwp}
   \cw'(T)= \sup\{\lambda \in \R\cup \{-\infty\}\mid \exists u \in (\R\cup\{-\infty\})^n, u\neq \bot, \; T(u) \geq \lambda +u\} \enspace.\end{align}

 For all $x\in (\rmax)^n$, we define $\operatorname{top}x:=\max_{i\in[n]}x_i$.
 We shall also consider
 \[
 \overline{\chi}(T):= \lim_k \operatorname{top}(T^k(0))/k = \inf_{k\geq 1} \operatorname{top}(T^k(0))/k \enspace.
 \]
 The existence of the limit and the fact it coincides with the infimum follow
 from the subadditivity property $\operatorname{top}(T^{k+l}(0))
 \leq \operatorname{top}(T^{k}(0)) + \operatorname{top}(T^{l}(0))$. Of course,
 when the limit $\chi(T) = \lim_{k}T^k(0)/k$ exists, we have
 $\overline{\chi}(T) = \operatorname{top} \chi(T)= \max_{i\in [n]}\chi_i(T)$.
 Then, $\overline{\chi}(T)$ may be interpreted as the value
 of a modified mean payoff game, in which Player Max chooses first
 the initial state $i\in[n]$, and then, the games starts from this
 state as described in~\Cref{sec-def-meanpayoff}. Thus, in the sequel,
 we shall refer to $\overline{\chi}(T)$ as the {\em upper mean payoff}
 associated to the operator $T$.

 The following result, which follows from~\cite{AGGut10},
 provides several spectral characterizations of
 this upper mean payoff. We say that a map $F$
 from $\R^n $ to $(\R\cup\{-\infty\})^n$ is {\em piecewise affine}
 if we can cover $\R^n$ by finitely many polyhedra in such a way
 that each coordinate map $F_i$ is either affine, or identically $-\infty$,
 on each of these polyhedra.

\begin{theorem}\label{th-duality}
  Let $T: (\rmax)^n\to (\rmax)^n$ be a Shapley operator.
  Then,
  \begin{align}\label{e-th-cw}
\cw'(T)=  \rho(T)=\overline{\chi}(T) = \cw(T)\enspace ,
  \end{align}
  and the suprema in ~\eqref{e-def-rho} and ~\eqref{e-def-cwp} are always achieved.
  
  Moreover, if the restriction of $T$ to $\R^n$ is piecewise affine,
  and if $\rho(T)\neq -\infty$,
  then the infimum in~\eqref{e-def-cw} is also achieved.
\end{theorem}
Before giving the details of the derivation of~\Cref{th-duality}
from~\cite{AGGut10}, we need to recall a result
of Kohlberg. An {\em invariant half-line} of a Shapley operator $T:\R^n\to \R^n$
is a pair $(u,\eta)\in \R^n\times\R^n$ such that
 \[
 T(u+s\eta) = u+(s+1)\eta , \qquad \forall s\geq 0 \enspace .
 \]
 Recall that
 a self-map of $\R^n$ is {\em non-expansive} for a fixed norm $\|\cdot\|$
 if $\|T(x)-T(y)\|\leq \|x-y\|$. Observe that a Shapley
 operator that preserves $\R^n$ is automatically
 non-expansive in the sup-norm (see e.g.~\cite{sgjg04}).
 \begin{theorem}[\cite{kohlberg}]\label{th-kohlberg}
   A piecewise affine map $T:\R^n\to \R^n$ that is nonexpansive
   in some norm admits an invariant half-line.
   \end{theorem}
 If a Shapley operator $T:\R^n\to \R^n$ has an invariant half-line $(u,\eta)$,
 it is immediate, using the fact that $T$
 is nonexpansive in the sup-norm, that
 $\chi(T) = \lim_k T^k(0)/k = \lim_k T^k(u)/k= \lim_k (u+k\eta)/k = \eta$.
Thus, the invariant half-line determines the mean payoff vector.

 \begin{proof}[Proof of~\Cref{th-duality}]
  The equalities in~\eqref{e-th-cw} are established in~\cite[Lemma 2.8]{AGGut10}, where they are derived from a theorem of Nussbaum concerning
  continuous, order preserving and positively homogeneous self-maps of the orthant, see ~\cite[Theorem~3.1]{MR818894} and also~\cite[Prop.~1]{sgjg04}.
  If $T(u) =\rho(T) + u$ with $u\neq \bot$, we also
  have $T(v) = \rho(T) + v$, where $v:= u-\operatorname{top}u$
  is such that $\operatorname{top} v = 0$. Using the compactness
  of $\{v\in (\rmax)^n\mid \operatorname{top}v=0\}$ and the
  continuity of $T$ on this set, we deduce that the supremum
  in~\eqref{e-def-rho} is always achieved. A similar argument
  shows that the supremum in~\eqref{e-def-cwp} is also achieved.

  Consider now $F(x) = T(x) \vee (\cw(T)+ x)$, which
  sends $\R^n$ to $\R^n$, and which is piecewise affine
  because the action spaces are finite. It is immediate
  that $\cw(F)=\cw(T)$. Let us take an invariant half-line
  $(u,\eta)$ of $F$.
  Then, it follows from $F^k(u)=u+k\eta$, and from
  the nonexpansiveness of $F$ in the sup-norm that $\bar{\chi}(T)=
  \lim_{k\to\infty} \operatorname{top}F^k(0)/k=
  \lim_{k\to\infty} \operatorname{top}F^k(u)/k= \operatorname{top}\eta$.
  Moreover,  $F(u)= u+\eta\leq u+\operatorname{cw}(F)$,
  and so, $T(u) \leq u+ \operatorname{cw}(T)$.
 \end{proof}
 \begin{proposition}\label{prop-finiteeig}
   A piecewise affine Shapley operator $T:\R^n\to\R^n$
   admits a finite eigenvector if and only if the mean payoff
   $\chi_i(T)$ is independent of the choice of the initial state $i\in[n]$.
 \end{proposition}
 \begin{proof}
   By~\Cref{th-kohlberg}, $T$ has an invariant half-line
   $(u,\eta)$ and $\chi(T)=\eta$. So, if $\chi_i(T)=\lambda$
   for all $i$, we have $T(u)=\lambda +u$, showing
   that $u$ is a finite eigenvector of $T$. Conversely,
   if $T(u)=\lambda + u$ for some $u\in\R^n$,
   then, using the nonexpansiveness of $T$,
   $\chi(T)=\lim_k T^k(0)/k=\lim_k T^k(u)/k=\lim_k (u+k\lambda) /k = (\lambda,\dots,\lambda)$. 
   \end{proof}

\section{Inner radius of a tropical polyhedron defined by generators}
\label{sec-inner}
For any subset $\mathcal{W}$ of $(\rmax)^n$,
we define the \NEW{inner radius} of $\mathcal{W}$, 
denoted $\inradius(\mathcal{W})$, as the supremum of the radii
of Hilbert's balls centered at a point in $\R^n$ and included
in $\Sp (\cW)$.
More generally, for all non-empty subsets $I\subset [n]$,
we define the {\em relative inner radius} of $\mathcal{W}$,
denoted by  $\inradius_I(\cW)$, as the
supremum of the radii
of Hilbert's balls centered at a point in the part $P_I$ of $(\rmax)^n$
and included in $\Sp (\cW)$. Thus, in particular,
$\inradius_{[n]}(\cW)=\inradius(\cW)$.
Observe that the relative inner radius 
depends only on the image of $\mathcal{W}\cap P_I$
in the tropical projective space $\proj(\rmax)^n$.

In~\cite{mtns}, it is shown that computing
the inner radius of
a tropical polyhedral cone given by an external description
$P=\{x\in (\rmax)^n \mid Ax\leq Bx \}$
reduces to computing the Collatz-Wielandt number $\cw(T)$
of a Shapley operator.

In this paper, we consider the somehow dual situation
in which the tropical polyhedral cone is
given by an internal description,
\[ \Col (V)=\{V x\mid x\in (\rmax)^p\} \enspace ,
\]
where $V$ is a $n\times p$ matrix with entries
in the tropical semifield $\rmax$, rather by an external description.
Recall that the size of an external description
of a tropical polyhedral cone can be exponential
in the size of an internal description,
and vice versa~\cite{AlGK09}.
This leads us to consider the following problem.
\begin{problem}\label{pb-innerradius}
  Input: a matrix $V\in \zmax^{n\times p}$. Goal: Compute the inner radius of $\Col(V)$. 
\end{problem}

We shall make the following assumption.
\begin{assumption}
  \label{assum-finrow}
  The matrix $V$ has no identically $-\infty$ row and no identically $-\infty$ column.
\end{assumption}
This is not restrictive. Indeed, let $I\subset [n]$ (resp.\
$J\subset[p]$) denote
the set of indices of non-identically $-\infty$ rows (resp.\ columns)
of $V$ and $V'$ denote the $I\times J$ submatrix of $V$.
For $K\subset [n]$, if $K$ is not included in $I$,
we have $\inradius_K(V)=-\infty$, whereas if $K=I$, then
$\inradius_K(V)=\inradius(V')$. More generally, the relative
inner radii of $V$ for $K\subset I$ coincide with the ones of $V'$
(up to permutations of rows of $V$).

In~\cite{AGGut10}, the tropical linear independence of the columns
of the matrix $V$ was studied by means of a specific Shapley operator,
which will also play a key role in our approach.
We set $E=\{(i,k)\in [n]\times [p]\mid V_{ik}\neq -\infty\}$.
Consider the operator $T: (\rmax)^n\to (\rmax)^n$, defined by
\begin{align}
T_i(x) = \inf_{k\in [p], (i,k)\in E} \Big[-V_{ik}+ \max_{j\in [n],j\neq i} (V_{jk}+x_j) \Big]\enspace .\label{e-def-T}
\end{align}
Owing to \Cref{assum-finrow}, the latter infimum is never
taken over an empty family,
so the operator does send $(\rmax)^n$ to $(\rmax)^n$.
We shall sometimes write $T_V$ instead of $T$ to emphasize
the dependence on $V$. Observe that $T$ is exactly the Shapley
operator of a mean payoff game defined in~\Cref{sec-def-meanpayoff}: the set of nodes
belonging to Player Min is $\minset:=[n]$,
the set of nodes belonging to Player Max
is $\maxset:= E$, %
with the set of allowed
moves
\begin{equation}\label{setAforT}
 \mathcal{A} = \{(i,(i,k)) \mid i\in [n],k\in [p], (i,k)\in E \}
\cup \{((i,k), j) \mid (i,k)\in E, (j,k)\in E, \; i\neq j\}
\enspace .
\end{equation}
The payment associated with the arc $(i,(i,k))$ is $w_{i,(i,k)} = -V_{ik}$,
whereas the payment associated with $((i,k),j)$ is $w_{(i,k),j} = V_{jk}$.

\begin{remark}\label{rem-nonnegative}
  In this game, the mean payoff $\chi_\imath (T)$ starting from any
  state $\imath$ is always nonpositive. Indeed, Player Min can always play
  a ``tit for tat'' policy, moving to state $(j,k)$ from
  state $j$, and thus, paying $-V_{jk}$ to Max, if the last move
  of Max was $(i,k)\to j$, so that Min paid $V_{jk}$. In this
  way, Min can cancel the last payment he made,
which guarantees a nonpositive mean payoff.
\end{remark}

Given a vector $a\in (\rmax)^n$, $a\neq \bot$, we define the
{\em tropical hyperplane}:
\[
\cH_a:= \{x\in (\rmax)^n\mid \max_{i\in [n]}(a_i+x_i) \text{ achieved at least twice} \} \enspace .
\]
Observe that $\cH_a$ depends only on the point in the tropical projective
space represented by $a$.  Moreover, $\cH_a$ is stable under the additive action
of scalars, so that $\cH_a$ can be identified with the subset of the tropical
projective space consisting of the equivalence classes of non-identically
$-\infty$ vectors of $\cH_a$.

For a finite vector $a \in \R^n$, the tropical hyperplane $\cH_a$ divides $(\rmax)^n$ into $n$ sectors $(S_i(a))_{i\in [n]}$, defined by
\begin{equation}\label{e-sectors}
S_i(a) :=\{ x \in (\rmax)^n~|~\forall j \in [n], \; x_i+a_i \geq x_j + a_j \} \enspace .
\end{equation}
The vector $-a$, which is unique up to an additive constant,
is called the {\em apex}
of $\cH_a$. Indeed, the set $\cH_a\cap \R^n$ modulo the scalar additions is the support of a polyhedral 
complex and $-a\in \R^n$ is the unique vertex (cell of dimension $0$)
of this complex. Then, we shall say that $\cH_a$ has a {\em finite apex}.
See~\Cref{fig-hyperplane} for an illustration.

\begin{figure}[htbp]
\begin{tikzpicture}
\coordinate (v1) at (\coord{-3.5}{0}{-1});
\coordinate (v2) at (\coord{0}{-3.5}{-1});
\coordinate (v3) at (\coord{0}{0}{-3});

\fill[red!30,opacity=0.7] (v1) -- (\coord{0}{0}{-1}) -- (v2)--cycle;
\fill[blue!30,opacity=0.7] (v1) -- (\coord{0}{0}{-1}) -- (v3)--cycle;
\fill[green!30,opacity=0.7] (v2) -- (\coord{0}{0}{-1}) -- (v3)--cycle;

\draw[thick] (\coord{0}{0}{-1}) -- (v1);
\draw[thick] (\coord{0}{0}{-1}) -- (v2);
\draw[thick] (\coord{0}{0}{-1}) -- (v3);
\filldraw[black] (\coord{0}{0}{-1}) circle (1.5pt) node[below,right] {$-a$};
\draw node[above] at (\coord{-1/2}{0}{-1}) {$\cH_a$};

\draw[dashed,->] (\coord{0}{0}{0}) -- (\coord{3.5}{0}{0}) node[above] {$x_1$};
\draw[dashed,->] (\coord{0}{0}{0}) -- (\coord{0}{3.5}{0}) node[above] {$x_2$};
\draw[dashed,->] (\coord{0}{0}{0}) -- (\coord{0}{0}{1.5}) node[above, right] {$x_3$};
\filldraw[black] (\coord{0}{0}{0}) circle (1.5pt) node[below,right] {0};

\end{tikzpicture}

\caption{The hyperplane $\cH_a$ with finite apex $a=(0,0,1)^\top$ and the sectors that $\cH_a$ defines in the projective space $\proj(\rmax)^3$.}
\label{fig-hyperplane}
\end{figure}

The following result shows that verifying whether there is a tropical
hyperplane containing a given collection of vectors reduces to solving
a mean payoff game. 
\begin{proposition}\cite[Corollary 4.8]{AGGut10}
  For $a\in (\rmax)^n$ such that $a\neq \bot$, suppose that
  $V\in (\rmax)^{n\times p}$
    satisfies \Cref{assum-finrow}, 
    and let $T$ be defined as above.
    Then, the following assertions are equivalent:
\begin{enumerate}
  \item $a \leq T(a) $;
  \item The column space $\Col(V)$ is included in $\cH_a$.
   \end{enumerate}
\end{proposition}
\begin{corollary}\label{cor-iff}
The columns of $V$ are contained in a tropical hyperplane
iff $\rho(T)$ is nonnegative.
\end{corollary}
\begin{proof}
  This follows from the equality
  $\rho(T)=\cw'(T)$ in \Cref{th-duality}
  and the fact the supremum is achieved
  in~\eqref{e-def-cwp}.
  \end{proof}

\begin{theorem}\label{th-inradius}
  Let $T=T_V$ be the Shapley operator associated to the matrix $V\in (\rmax)^{n\times p}$ defined in~\eqref{e-def-T}. Then, $\rho(T)\leq 0$. Moreover,
  \[ -\rho(T) = \inradius(\Col(V)).
  \]
If $\rho(T)$ is finite, a maximal Hilbert's ball
  included in $\Col(V)\cap \R^n$ is given by $B(-a,-\rho(T))$ where
  $a$ is any vector in $\R^n$ such that $T(a) \leq \rho(T) + a$.
\end{theorem}
We will deduce \Cref{th-inradius} from the following lemma:
\begin{lemma}\label{lem-ball}
  For all $\lambda\in [-\infty,0]$ and $a\in \R^n$,
\[ B(-a,-\lambda) \subset \operatorname{Col} V\iff T(a) \leq \lambda + a
\]
\end{lemma}
\begin{proof}
  Suppose first that $\lambda$ is finite.
  Then, considering~\eqref{e-def-T}, we see that $T(a) \leq \lambda + a$ is equivalent to 
\begin{equation}\label{SupEigenvalue}
\forall i\in [n], \enspace \exists k \in [p], \enspace \forall j\in [n],\enspace j\neq i, \quad -\lambda -a_i +a_j\leq V_{ik}-V_{jk} \enspace.
\end{equation}
Let $x\in \R^n$, we have $x\in B(-a,-\lambda)$ if and only if
\begin{equation}\label{Ball}
\forall i\in [n], \enspace \forall j\in [n] , \quad x_i-x_j \leq -\lambda -a_i +a_j\enspace.
\end{equation}
Moreover, the basic properties of residuation entail that $VV^\sharp \leq \I$,
where $V^\sharp x$ is the maximal element $y$ such that $Vy\leq x$.
It follows that $x\in \Col(V)$ if and only if $x=VV^\sharp x$,
or equivalently, $x\leq VV^{\sharp} x$.
The latter property can be rewritten as $x_i\leq \max_{k\in [p]} \{V_{ik}+\min_{j\in[n]}(-V_{jk}+x_j)\}$, for all $i\in [n]$,
which is equivalent to 
\begin{equation}\label{ColV}
\forall i\in [n], \enspace \exists k \in [p], \enspace \forall j\in [n] , \quad x_i-x_j \leq V_{ik}-V_{jk} \enspace.
\end{equation}

We can see that if \cref{SupEigenvalue} and \cref{Ball} are true then \cref{ColV} follows, which shows the ``$\Leftarrow$'' direction of the lemma.

Now, we suppose that $B(-a,-\lambda) \subset \operatorname{Col} V$. For a given $i\in [n]$, we consider the vector $x^{(i)}\in \R^n$ given by $x^{(i)}_i=-\lambda -a_i$ and $x^{(i)}_j=-a_j$ for all $j\neq i$. Since $\lambda \leq 0$, we have $x^{(i)}\in B(-a,-\lambda)$, then $x^{(i)}\in \Col(V)$. Therefore by \cref{ColV}, there exists $k\in [p]$ such that $\forall j\in [n] , \enspace x^{(i)}_i-x^{(i)}_j \leq V_{ik}-V_{jk}$. Moreover, we have $\forall j\in [n], j\neq i$, $x^{(i)}_i-x^{(i)}_j=-\lambda -a_i+a_j$. Finally this yields \cref{SupEigenvalue}, which proves that $T(a) \leq \lambda + a  $.

We finally show that the conclusion of the lemma is still
true when $\lambda = -\infty$. This follows from $B(-a,+\infty)
= \cup_{\mu\in (-\infty,0)} B(-a,-\mu)$ and $-\infty + a = \inf_{\mu\in (-\infty,0)} \mu+a$.
\end{proof}
\begin{proof}[Proof of \Cref{th-inradius}]
  If $B(-a,-\lambda)\subset \Col(V)$
  for some finite $a$, with $\lambda\leq 0$, by \Cref{lem-ball},
  we see that $T(a)\leq \lambda +a$, and we deduce
  from the Collatz-Wielandt property (\Cref{th-duality})
  that $\rho(T)\leq \lambda$, and so, the radius
  of the ball, $-\lambda$, is bounded above by $-\rho(T)$.

  Moreover, it follows from \Cref{assum-finrow} that $\Col(V)$
  has a finite vector $a$; indeed, we can take for $a$
  the supremum of the columns of $V$. Then, $B(-a,0)\subset \Col(V)$,
  and by the previous observation, $0\leq -\rho(T)$.

If $\rho(T)=-\infty$, then using the expression of the Collatz-Wielandt number of $T$, we get that for all finite $\lambda\leq 0$, there exists a finite vector $a\in \R^n$
  such that $T(a)\leq \lambda +a$. By~\Cref{lem-ball},
  this implies that $B(-a,-\lambda )\subset \Col(V)$, and so
$\inradius(\Col(V))\geq -\lambda$. Since this holds for all $\lambda\leq 0$,
we deduce that $\inradius(\Col(V))=+\infty= -\rho(T)$ is the supremum of the
radius of a Hilbert's ball included in $\Col(V) \cap \R^n$.

  Finally, if $\rho(T)$ is finite, since the infimum is attained in the expression of the Collatz-Wielandt number of $T$ (see \Cref{th-duality}), there exists a finite vector $a\in \R^n$
  such that $T(a)\leq \rho(T)+a$. By~\Cref{lem-ball},
  this entails that $B(-a,-\rho(T))\subset \Col(V)$.

  This shows that $-\rho(T)$ is the maximum radius of a Hilbert's
  ball included in $\Col(V) \cap \R^n$.
\end{proof}
\begin{remark}
  One can give an alternative, less direct proof, of~\Cref{th-inradius}
  by deriving it from Theorem~16 of~\cite{mtns}.
  The latter result shows that if $T$ is a Shapley operator
  which satisfies the technical condition 
  ($T$ must be ``diagonal free''),
  then, the supremum of the radii of Hilbert's balls
  included
  in $S(T):= \{x\in \R^n\mid x\leq T(x)\}$ coincides with
  $\sup \{\mu \in \R\mid \exists v\in \R^n,
  \mu +v \leq T(v)\}$. The initial part of the proof
  of~\Cref{lem-ball}, up to~\eqref{ColV},
  entails that $\Col(V)$ is precisely the set
  of vectors $x$ such that $x\leq -T(-x)$.
  \end{remark}
The following is an immediate corollary of \Cref{th-inradius}
\begin{corollary}\label{cor-eig}
  The set $\operatorname{Col}(V)\cap \R^n$ is of empty interior
  if and only $\rho(T)=0$.\hfill\qed
\end{corollary}
By combining \Cref{cor-eig} and~\Cref{cor-iff}, we
recover the following known result,
established in~\cite{DSS} (when the entries
of the matrix $V$ are finite).
\begin{corollary}[Compare with Th.~4.2 of~\cite{DSS}]\label{cor-nempty}
  The set $\Col(V)\cap \R^n$ is of empty interior if and only if $\Col(V)$
  is included in a tropical hyperplane. \hfill\qed
\end{corollary}

The following additional corollary implies that we can check in polynomial
time whether the inner radius of $\Col(V)$ is finite.
\begin{corollary}
  The following assertions are equivalent:
  \begin{enumerate}
  \item\label{i-1} The inner radius of $\Col(V)$ is infinite;
  \item\label{i-2} There is no part of $(\rmax)^n$ that is left invariant by the operator $T$;
  \item\label{i-3} $T^n(0)$ is the vector identically equal to $-\infty$;
    \item\label{i-4} $\rho(T)=-\infty$.
    \end{enumerate}
\end{corollary}
\begin{proof}
  \eqref{i-3}$\Rightarrow$\eqref{i-1}: Suppose that $T^n(0)$ is equal to $\bot$, the identically $-\infty$ vector. Let us take $u\in (\rmax)^n$, not identically $-\infty$, such that
  $T(u) = \rho(T) +u$. Then, there is a constant $\alpha\in \mathbb{R}$
  such that $u_i\leq \alpha$, for all $i\in [n]$, and
  so $n\rho(T)+u = T^n(u) \leq T^n(0)+\alpha$ is the identically $-\infty$ vector. It follows that $\rho(T)=-\infty$. Then, by~\Cref{th-inradius}, the inner
  radius of $\Col(V)$ is infinite.

  \eqref{i-1}$\Rightarrow$\eqref{i-2}: Let $I$ be a non-empty subset
  of $[n]$, and suppose that the part $P_I$ consisting of vectors
  of $(\rmax)^n$ of support $I$ is left invariant by $T$.
  Let $u$ be the vector in this part such that $u_i=0$ for all $i\in I$.
  Since $T(u)\in P_I$, there exists a real number $\alpha$ such that
  $T(u)\geq \alpha + u$. Hence, $\rho(T)=\cw'(T)\geq \alpha>-\infty$,
  and, by~\Cref{th-inradius}, $\inradius(\Col(V))=-\rho(T)<+\infty$.

  \eqref{i-2}$\Rightarrow$\eqref{i-3}:
  Consider the map $\pi: (\rmax)^n \to \mathcal{P}([n])$,
  which sends a vector $u$ to its support,
  $\pi(u)=\{i\in [n]\mid u_i \neq -\infty\}$, and consider the equivalence relation $\ker \pi$ on $(\rmax)^n$,
  such that $(x,y)\in \ker \pi$ iff $\pi(x)=\pi(y)$.
  The quotient set $(\rmax)^n/\ker \pi$  can be
  identified to $\mathcal{P}([n])$, and the order
  on $(\rmax)^n$ induces an order on $(\rmax)^n/\ker \pi$,
  corresponding to the inclusion order on $\mathcal{P}([n])$.
  The elements
  if $(\rmax)^n/\ker \pi$ are precisely the parts of $(\rmax)^n$,
  together with the singleton consisting of the identically
  $-\infty$ vector. Let $\bot=\pi^{-1}(\emptyset)$ denote this singleton, and let $\top=\pi^{-1}([n])$.
  Observe that $\top$ is the
  maximal element of $(\rmax)^n/\ker \pi$,
  and that $\bot$ is its minimal element.
  
  Since the operator $T$ is order preserving and commutes
  with the addition of a constant, it induces
  a map $T_\pi$ from $(\rmax)^n/\ker \pi$ to itself,
  which is still order preserving. 
  Moreover, the fixed points of $T_\pi$ distinct from
  $\bot$ are precisely the parts of $(\rmax)^n$ that are
  invariant by $T$.
We have, $T_\pi(\top)\leq \top$, from which we deduce
 that $((T_\pi)^k(\top))_{k\geq 0}$ is a nonincreasing sequence.
 If $(T_\pi)^k(\top)=(T_\pi)^{k+1}(\top)\neq \bot$, for some
 $k$, then $(T_\pi)^k(\top)$ would be an invariant part of $T$,
 contradicting the assumption.
 It follows that the sequence $(T^k_\pi(\top))_{k\geq 0}$ strictly decreases
 until it reaches $\bot$. Since the maximal cardinality of a chain
 in the lattice $\mathcal{P}(n)$ is $n+1$, it follows
 that $(T_\pi)^n(\top)=\bot$. Hence,
 $T^n(0)$ is the identically $-\infty$ vector.

 Finally, the equivalence between \eqref{i-1} and~\eqref{i-4}
follows from~\Cref{th-inradius}.
  \end{proof}

Recall that a vector $u$ in a tropical cone $\cV\subset (\rmax)^n$ is {\em extreme}~\cite{GK06,BSS} if $u=v\vee w$ with $v,w\in  \cV$
implies that $u=v$ or $u=w$. An {\em extreme direction} of $\cV$
is of the form $\rmax + u$, for some extreme vector of $\cV$, i.e.,
it consists of the tropical scalar multiples of $u$.
We say that a tropical cone in $(\rmax)^n$ is {\em simplicial}
if it has precisely $n$ extreme directions. 
\begin{proposition}\label{coro-simplicial}
If a Hilbert's ball of positive radius is included in $\Col(V)$,
then it is also included in a simplicial tropical cone
generated by some $n$ columns of $V$.
\end{proposition}
\begin{proof}
   For all maps $\sigma: [n]\to [p]$, such that $(i,\sigma(i))\in E$,
   we consider the Shapley
   operator of the one-player game obtained when player MIN
   selects the action $k=\sigma(i)$ in state $i$, that is,
   \[
  T^{\sigma}: (\rmax)^n\to (\rmax)^n, \qquad T_i^\sigma(x)=
-V_{i\sigma(i)}+ \max_{j\in [n],j\neq i} (V_{j\sigma(i)}+x_j) \enspace .
\]
If $B(-a,-\lambda)\subset \Col(V)$, then, by \Cref{lem-ball},
$T(a) \leq \lambda +a$. So, by choosing $k = \sigma(i)$ that achieves
the minimum in the expression of $T(a)$ in~\eqref{e-def-T},
we get $T^\sigma(a) \leq \lambda +a$. 
Let $J:= \sigma([n])$,
so that $|J|\leq n$. Since $(i,\sigma(i))\in E$ holds for all $i\in [n]$,
the submatrix $V[J]$ of $V$, obtained
by keeping the columns in $J$, cannot have a $-\infty$ row. Hence~\Cref{lem-ball}
can be applied to $V[J]$. We deduce that $B(-a,-\lambda)\subset
\Col(V[J])$. 
Up to eliminating elements of $J$, we may assume that
the set $J$ is minimal to generate $\Col(V[J])$.

Let $u^j$ denote the $j$th column of $V$. 
Then, every $u^j$ must be extreme
in $\Col(V[J])$. Indeed, suppose that $u^j=v\vee w$ with $v,w\in \cV$ with $u^j\neq v$ and $u^j\neq w$.
Then, we can write $v=\vee_{k\in J} (\lambda_k + u^k)$
and $w=\vee_{k\in J} (\mu_k + u^k)$, for some $\lambda_k,\mu_k \in \rmax$.
Moreover, we must have $\lambda_j <0$, otherwise, $v\geq u^j$,
and since $v\leq v\vee w=u^j$, $v=u^j$, a contradiction.
A similar result holds for $w$. 
Since $u^j=v\vee w=\vee_{k\in J} ((\lambda_k\vee\mu_k) + u^k)$,
and $\lambda_j\vee \mu_j<0$, we deduce that 
$u^j=\vee_{k\in J\setminus\{j\}} ((\lambda_k\vee\mu_k) + u^k)$
is generated by the columns $\{u^k\mid k\in J\setminus\{j\}\}$,
contradicting the minimality of $J$. It follows that every column
of $V[J]$ is extreme in $\Col(V[J])$.

To show that $\Col(V[J])$
is simplicial, it remains to check that $|J|\geq n$.
It is known that a collection of at most $n-1$ vectors
in $(\rmax)^n$ is included in a tropical hyperplane -- this
follows for instance from a tropical analogue of the Radon theorem,
see e.g.~\cite{But} or~\cite[Coro.~6.13]{AGG08};
or this can be deduced from the characterization
of the tropical rank~\cite{DSS,IRowen,AGGut10}.  
So if $|J|<n$, then $\Col(V[J])$ is of empty interior,
contradicting $B(-a,-\lambda)\subset \Col(V[J])$.
\end{proof}

We get as a corollary the following result.

\begin{corollary}\label{cor-simplicial}
  We have
  \begin{align}
    \inradius(\Col(V)) = \max_J \inradius(\Col(V[J]))
    \label{e-morphism-inradius}
  \end{align}
  where the maximum is taken over all subsets $J\subset [p]$ of cardinality
  $n$. Moreover, if the inner radius is positive, the maximum
  is achieved by $J$ such that $\Col(V[J])$ is simplicial. 
\end{corollary}
By convention, if $p<n$, the maximum in~\eqref{e-morphism-inradius}
is zero.
\begin{proof}
  The inequality $\geq$ in~\eqref{e-morphism-inradius}
  is trivial. If $\inradius(\Col(V))=0$, the equality trivially
  holds in~\eqref{e-morphism-inradius}. If $\inradius(\Col(V))>0$,
then for all  $0\leq \lambda < \inradius(\Col(V))$, there exists
a Hilbert's ball of radius $\lambda$ included in $\Col(V)$.
By \Cref{coro-simplicial}, this ball is also
included in a simplicial tropical cone generated by columns of $V$,
which means that there exists $J\subset [p]$ of cardinality $n$
such that $\lambda\leq \inradius(\Col(V[J]))\leq \max_J \inradius(\Col(V[J]))$.
Since this holds for all $0\leq \lambda < \inradius(\Col(V))$,
we deduce  the inequality 
$\inradius(\Col(V)) \leq  \max_J \inradius(\Col(V[J]))$ and so the 
equality.
  \end{proof}

\begin{corollary}\label{coro-equiv}
  Computing the inner-radius of a tropical polyhedron
(\Cref{pb-innerradius})
  is polynomial-time Turing equivalent to mean payoff games
  (\Cref{pb-1}).
\end{corollary}
\begin{proof}
  We observed immediately after stating~\Cref{pb-1} that the problem
  of computing $\chi_i(T)$, where $T$ is the Shapley operator
  of a deterministic mean payoff game, satisfying Assumption~\ref{assump1},
  polynomially Turing-reduces to mean payoff games.
  By~\Cref{th-inradius}, the opposite of the inner-radius is equal to $\rho(T)$. Since, $\rho(T)=\max_{i\in[n]} \chi_i(T)$, computing the inner-radius polynomially Turing-reduces to mean payoff games.
  
  Conversely, Corollary~3.11 of~\cite{podolskii}
  shows in particular that mean payoff games (\Cref{pb-1})
  polynomially Turing-reduces to checking
  whether a collection of vectors $v^1,\dots,v^p$ of $(\zmax)^n$ are included
  in a tropical hyperplane. By~\Cref{cor-nempty} and~\Cref{cor-eig},
  the latter problem is equivalent to checking whether
  the inner-radius of a tropical polyhedral cone $\Col(V)$ vanishes.
\end{proof}
\begin{corollary}\label{cor-center}
  Computing the center of a Hilbert's ball of maximal radius
  included in $\Col(V)$, where $V\in \zmax^{n\times p}$,
  polynomially Turing-reduces to mean payoff games.
\end{corollary}
\begin{proof}
  We first compute the maximal radius, $-\rho(T)$, which has been
  noted above, polynomially Turing-reduces to mean payoff games.
  We can also obtain by the same type of reduction
  an optimal policy $\sigma$ of Player Min,
  which satisfies $\rho(T)= \rho(T^\sigma)$. Indeed,
  for each move of player Min $i\to j$, we can consider
  a modified Shapley operator $T^{(i,j)}$ corresponding
  to the game in which player Min makes the move $i\to j$
  when in node $i$ (i.e., this player has no choice in node $i$),
  and all the other allowed moves are unchanged. By checking
  whether $\rho(T^{(i,j)})=\rho(T)$, we can verify if
  the move $i\to j$ belongs to an optimal policy of Player Min.
  By repeatedly restricting the freedom of moves of Player Min,
  we arrive, after a polynomial number of evaluation
  of $\rho(\cdot)$, at such an optimal policy $\sigma$.
  We showed that the center of an optimal Hilbert's
  ball is of the form $-u$ where $u\in \R^n$ and $T(u)\leq \rho(T)+u$. Since $T\leq T^\sigma$, and $\rho(T)=\rho(T^\sigma)$,
  it suffices
  to construct a vector in $\R^n$ such that $T^\sigma(u)\leq \rho(T^\sigma)+u$.
  Considering the tropically linear map
  $B:= -\rho(T^\sigma)+T^\sigma$, we see this is equivalent
  to $Bu\leq u$. A standard result of tropical spectral
  theory
  shows that one can compute such a vector $u$ by solving
  a shortest path problem.
  Actually, a tropical generating family of the set of such vectors $u$ is the set of columns of the so called ``metric closure'' or ``Kleene star'' $B^*$ of the matrix $B$, defined as the tropical sum of the tropical powers of $B$, 
  see e.g.~\cite[Th.~3.101]{BCOQ92} and~\cite[\S~4.4]{butkovic}.
  Moreover, the tropical sum $u$ of the columns of $B^*$
  is a finite vector. In this way, we constructed
  $u$ such that $Bu\leq u$, and so $T(u) \leq \rho(T)+u$. 
    \end{proof}
\begin{remark}
  A subset $J\subset [p]$ satisfying
  $\inradius(\Col(V)) = \max_J \inradius(\Col(V[J]))$
  can be computed by using any mean payoff game algorithm
  that returns, together with the mean payoff $\overline{\chi}(T)$,
  a vector $u\in \R^n$ such that $T(u) \leq \overline{\chi}(T)+ u$.  Indeed,
  we saw in the proof of~\Cref{coro-simplicial} that, taking any policy
  $\sigma$ such that $T(u)=T^{\sigma}(u)$, and setting $J:=\sigma([n])$,
  we have $B(-u,-\overline{\chi}(T))\subset \Col(V[J])$.
\end{remark}

\begin{figure}[htbp]
\begin{center}

\def\coord#1#2#3{{-sqrt(3)/2*(#1-#2)} ,{ -(1/2)*#1 - (1/2)*#2 + #3}}

\begin{tikzpicture}

\coordinate (v1) at (\coord{-3}{0}{-1});
\coordinate (v2) at (\coord{0}{-3}{-1});
\coordinate (v3) at (\coord{0}{0}{-4});
\coordinate (v4) at (\coord{1}{0}{-2});
\coordinate (v5) at (\coord{1}{-1}{-1});
\coordinate (v6) at (\coord{-1}{1}{-1});
\coordinate (v7) at (\coord{0}{1}{-2});
\coordinate (v8) at (\coord{0}{-1}{0});
\coordinate (v9) at (\coord{-1}{0}{0});

\coordinate (f1) at (\coord{-2}{0}{-1});
\coordinate (f2) at (\coord{0}{-2}{-1});
\coordinate (f3) at (\coord{0}{0}{-3});

\coordinate (b1) at (\coord{0}{0}{0});
\coordinate (b2) at (\coord{0}{1}{0});
\coordinate (b3) at (\coord{0}{1}{-1});
\coordinate (b4) at (\coord{0}{0}{-2});
\coordinate (b5) at (\coord{1}{0}{-1});
\coordinate (b6) at (\coord{1}{0}{0});

\filldraw[black,draw=black,opacity=0.8,very thick] (v1) -- (f1) -- (v9) -- (b1) -- (v8) -- (f2) -- (v2) -- (f2) -- (v5) -- (b5) -- (v4) -- (f3) -- (v3) -- (f3) -- (v7) -- (b3) -- (v6) -- (f1) -- (v1) -- cycle;
\filldraw[gray] (f1) -- (v9) -- (b1) -- (v8) -- (f2) -- (v5) -- (b5) -- (v4) -- (f3) -- (v7) -- (b3) -- (v6) -- (f1) -- cycle;

\filldraw[blue!80,draw=blue!80,opacity=0.7, thick] (b1) -- (b2) -- (b3) -- (b4) -- (b5) -- (b6) -- (b1) -- cycle;

\filldraw[black] (v1) circle (1.5pt) node[above] {$V_{\cdot 1}$};
\filldraw[black] (v2) circle (1.5pt) node[above] {$V_{\cdot 2}$};
\filldraw[black] (v3) circle (1.5pt) node[below,right]{$V_{\cdot 3}$};
\filldraw[black] (v4) circle (1.5pt) node[below,left] {$V_{\cdot 4}$};
\filldraw[black] (v5) circle (1.5pt) node[below,left] {$V_{\cdot 5}$};
\filldraw[black] (v6) circle (1.5pt) node[below,right] {$V_{\cdot 6}$};
\filldraw[black] (v7) circle (1.5pt) node[below,right] {$V_{\cdot 7}$};
\filldraw[black] (v8) circle (1.5pt) node[above] {$V_{\cdot 8}$};
\filldraw[black] (v9) circle (1.5pt) node[above] {$V_{\cdot 9}$};

\draw[dashed,->] (\coord{0}{0}{0}) -- (\coord{5}{0}{0}) node[above] {$x_1$};
\draw[dashed,->] (\coord{0}{0}{0}) -- (\coord{0}{5}{0}) node[above] {$x_2$};
\draw[dashed,->] (\coord{0}{0}{0}) -- (\coord{0}{0}{2}) node[above, right] {$x_3$};
\filldraw[black] (\coord{0}{0}{0}) circle (1.5pt) node[below,right] {0};

\end{tikzpicture}
\end{center}

\caption{Example of an inner ball of the column space $\Col (V)$ in the projective space $\proj(\rmax)^3$, where $V=\left(\begin{array}{ccccccccc}
  -3 & 0  &  0 &  1  & 1 & -1 & 0 & 0 & -1\\
   0 & -3 &  0  & 0  & -1 & 1 & 1 & -1 & 0\\
   -1 & -1  & -4 & -2  &  -1  &  -1  & -2  & 0  & 0
\end{array}
\right)$.}
\label{fig-inner-ball}
\end{figure}

\if{
\begin{figure}[htbp]
\renewcommand{\baryx}{x_1}
\renewcommand{\baryy}{x_2}
\renewcommand{\baryz}{x_3}
\colorlet{mygreen}{gray!50!black}

\begin{center}
\begin{tikzpicture}%
[scale=0.65,>=triangle 45
,vtx/.style={mygreen},
ray/.style={myred}]
\equilateral{7}{100};

\barycenter{g1}{\expo{-6}}{\expo{0}}{\expo{0}};
\barycenter{g2}{\expo{0}}{\expo{-6}}{\expo{0}};
\barycenter{g3}{\expo{0}}{\expo{0}}{\expo{-6}};
\barycenter{h1}{\expo{2}}{\expo{0}}{\expo{0}};
\barycenter{h1a}{\expo{1}}{\expo{0}}{\expo{-1}};
\barycenter{h1b}{\expo{1}}{\expo{-1}}{\expo{0}};
\barycenter{h1c}{\expo{1}}{\expo{0}}{\expo{0}};
\barycenter{h2}{\expo{0}}{\expo{2}}{\expo{0}};
\barycenter{h2a}{\expo{-1}}{\expo{1}}{\expo{0}};
\barycenter{h2b}{\expo{0}}{\expo{1}}{\expo{-1}};
\barycenter{h2c}{\expo{0}}{\expo{1}}{\expo{0}};
\barycenter{h3}{\expo{0}}{\expo{0}}{\expo{2}};
\barycenter{h3a}{\expo{0}}{\expo{-1}}{\expo{1}};
\barycenter{h3b}{\expo{-1}}{\expo{0}}{\expo{1}};
\barycenter{h3c}{\expo{0}}{\expo{0}}{\expo{1}};
\barycenter{f1}{\expo{-2}}{\expo{0}}{\expo{0}};
\barycenter{f2}{\expo{0}}{\expo{-2}}{\expo{0}};
\barycenter{f3}{\expo{0}}{\expo{0}}{\expo{-2}};

\barycenter{e1}{\expo{0}}{0}{0};
\barycenter{e2}{0}{\expo{0}}{0};
\barycenter{e3}{0}{0}{\expo{0}};

\barycenter{k1}{\expo{-1}}{\expo{0}}{\expo{0}};
\barycenter{k2}{\expo{0}}{\expo{-1}}{\expo{0}};
\barycenter{k3}{\expo{0}}{\expo{0}}{\expo{-1}};

\filldraw[gray,draw=black,opacity=0.9,very thick] (g1) -- (f1) -- (h3b) -- (h3c) -- (h3a) -- (f2) -- (g2) -- (f2) -- (h1b) -- (h1c) -- (h1a) -- (f3) -- (g3) -- (f3) -- (h2b) -- (h2c) -- (h2a) -- (f1) -- cycle;

\filldraw[red,draw=red,opacity=0.9, thick] (k1) -- (h3c) -- (k2) -- (h1c) -- (k3) -- (h2c) -- (k1) -- cycle;

\filldraw[vtx] (g1) circle (0.75ex) node[below] {$G_{1\cdot }$};
\filldraw[vtx] (g2) circle (0.75ex) node[below] {$G_{2\cdot }$};
\filldraw[vtx] (g3) circle (0.75ex) node[below] {$G_{3\cdot }$};

\filldraw[vtx] (h1a) circle (0.75ex) node[below,left] {$G_{4a\cdot }$};
\filldraw[vtx] (h1b) circle (0.75ex) node[below,left] {$G_{4b\cdot }$};
\filldraw[vtx] (h2a) circle (0.75ex) node[below,right] {$G_{5a\cdot }$};
\filldraw[vtx] (h2b) circle (0.75ex) node[below,right] {$G_{5b\cdot }$};
\filldraw[vtx] (h3a) circle (0.75ex) node[above] {$G_{6a\cdot }$};
\filldraw[vtx] (h3b) circle (0.75ex) node[above] {$G_{6b\cdot }$};

\end{tikzpicture}
\end{center}
\caption{A template inspired by https://arxiv.org/pdf/1004.2778.pdf}
\label{fig-inner-ball}
\end{figure} }\fi

We can verify easily that $\lambda=-1$ and
$a=(0, 0, 1)^\top$
satisfy $T(a)=\lambda + a$.
Moreover, a policy $\sigma$ such that $T^\sigma(a)=T(a)$
is given by $\sigma(1)=4$, $\sigma(2)=6$ and $\sigma(3)=8$.
Therefore, by~\Cref{th-inradius} the maximal radius of a Hilbert's
ball included in $\operatorname{Col} V$ is $-\lambda=1$.
Moreover, a maximal Hilbert's ball is given by $B(a,1)$,
and $B(a,1)$ is included in the simplicial cone $\Col(V[J])$
where $J=\{4,6,8\}=\sigma([3])$. 
This Hilbert's ball, together with the simplicial cone $\Col(V[J])$, are 
shown in~\Cref{fig-inner-ball}. Observe
that the set $J$ such that $\inradius(\Col(V[J]))=\inradius(\Col(V))$
is not unique, indeed, every $J'=\{i,j,k\}$ with $i\in \{4,5\}$,
$j\in \{6,7\}$ and $k\in \{8,9\}$ is a candidate.

\section{The strong duality theorem for tropical linear regression}
\label{sec-strong}

In this section we will study the best approximation of a set of points in the tropical projective space by a tropical hyperplane.
We will show that the best error of approximation is equal to the inner radius of the tropical module generated by this set of points. 

Let $\cV=\{v^{(1)},\cdots,v^{(p)}\} \subset \proj(\rmax)^n$ be
a finite subset of the tropical projective space. Since we mainly focus on $\Sp (\cV)$, by abusing notions, we denote by $V\in (\rmax)^{n\times p}$ the matrix whose columns are given by some representatives of $v^{(1)},\cdots,v^{(p)}$.
Note that $\Sp (\cV)= \Col (V)$, which does not depend on the choice of the representatives of $v^{(1)},\cdots,v^{(p)}$.
In the following, we use the notation $\inrv{\cV}=\inradius (\Sp (\cV))$.

We introduce a one-sided Hausdorff distance from a set $A \subset \proj(\rmax)^n$ to a set $B\subset \proj(\rmax)^n$ with respect to the Hilbert's projective metric, which we shall call the \NEW{Hilbert's distance} from $A$ to $B$:
\begin{align}
\dist_H(A,B) := \sup_{a\in A} \dist_H(a,B)\enspace,\quad \text{with}\; \dist_H(a,B):= \inf_{b\in B} d(a,b) \enspace.\label{def-hausdorff}
\end{align}

Note that $\dist_H(A,B)=0$ if and only if for each part $P_I$
of the tropical projective space, $A\cap P_I$ is included in
the closure of $B\cap P_I$ with respect to the relative topology of $P_I$.

We are interested in the following
{\em tropical linear regression problem}, consisting
of finding a best hyperplane approximation
of the set $\cV$ in Hilbert's distance:
\begin{equation}\label{e-regression}
\inf_{a\in \proj(\rmax)^n} \dist_H (\cV,\cH_a) \enspace .
\end{equation}
Observe that if there is an index $i\in [n]$ such that
$v^{(1)}_i=\dots=v^{(p)}_i=-\infty$, then the tropical
linear regression problem is trivially solved by
setting $a_i =0$ and $a_j=-\infty$ for $j\neq i$.
Hence, in the sequel, we shall assume that the matrix $V$ satisfies~\Cref{assum-finrow}.
In particular, considering the operator $T$ defined in~\Cref{e-def-T},
we know from~\Cref{th-inradius} that the inner radius of $\Sp (\cV)$ is $-\rho(T)$.

The following lemma gives a simple formula for the Hilbert's distance from a point to a hyperplane.
\begin{lemma}\label{lem-dist-hyperplane}
  For $x,a\in \proj(\rmax)^n$,
  let $i^*\in \argmax_{i\in [n]} (x_i+a_i)$. Then the Hilbert's distance from the point $x$ to the hyperplane $\cH_a$ is
\begin{equation}\label{dist-hyperplane}
\dist_H (x,\cH_a) = x_{i^*}+a_{i^*} - \max_{i\in [n], i\neq i^*} (x_i+a_i) \enspace ,
\end{equation}
where we use the convention $(-\infty)-(-\infty)=0$.
\end{lemma}

\begin{proof}
If $\max_{i\in [n]} (x_i+a_i) = -\infty$, then $x\in \cH_a$ and ~\Cref{dist-hyperplane} holds with the convention $(-\infty)-(-\infty)=0$.
If $\max_{i\in [n]} (x_i+a_i)\neq -\infty$ and the maximum in the expression is attained twice, then $x\in \cH_a$ and~\Cref{dist-hyperplane} holds. 

Now we focus on the case $\max_{j\in [n], j\neq i^*} (x_j+a_j)  < x_{i^*}+a_{i^*}\in \R$, which implies $x_{i^*}\in \R$ and $a_{i^*}\in \R$. We split the argument into the following two cases.

Case 1: $\max_{j\in [n], j\neq i^*} (x_j+a_j) \in \R$. 
Then $\delta:= x_{i^*}+a_{i^*} - \max_{j\in [n], j\neq i^*} (x_i+a_i) > 0$. Consider the point $\tilde{x}$ given by
\[
\left\{
\begin{array}{ll}
\tilde{x}_{i^*}=x_{i^*} - \delta &,\\
\tilde{x}_j=x_j &, \text{ for } j\in [n], j\neq i^*
\end{array}
\right.
\]
Then $\tilde{x}\in \cH_a$ and $d(x,\tilde{x})=\delta$, implying $\dist_H (x,\cH_a) \leq \delta$.
Now, let $x'\in \cH_a$, then the maximum in $\max_{j\in [n]} (x'_j+a_j)$ is achieved at least twice. So there exists $i\neq i^*$, such that $\max_{j\in [n]} (x'_j+a_j)=x'_i+a_i$. 
Since $i \neq i^*$, we have $\delta \leq x_{i^*}+a_{i^*}-(x_i+a_i)$, then $x'_i-x_i \geq x'_i+a_i -(x_{i^*}+a_{i^*}) +\delta$. 
Since $x'_{i^*}+a_{i^*} \leq \max_{k\in [n]} (x'_k+a_k) = x'_i+a_i$, then $x'_{i^*}-x_{i^*}\leq x'_i+a_i -(x_{i^*}+a_{i^*})$.
Therefore $d(x,x') \geq (x'_i-x_i)-(x'_{i^*}-x_{i^*}) \geq \delta$, which proves 
$\dist_H(x,\cH_a)\geq \delta$.

Case 2: $\max_{j\in [n], j\neq i^*} (x_j+a_j) = - \infty$. For $x'\in \cH_a$, there exists $i\neq i^*$, such that $\max_{j\in [n]} (x'_j+a_j)=x'_i+a_i$. If $x'_i+a_i= - \infty$, then $x'_{i^*}-a_{i^*}=-\infty$. Since $a_{i^*} \in \R$, then $x'_{i^*}=-\infty$. Thus the fact that $x_{i^*}\in \R$ forces $d(x,x')=+\infty$, i.e., $\dist_H(x,\cH_a)=+\infty$ and~\Cref{dist-hyperplane} holds. Now if $x'_i+a_i\in \R$, then $x'_{i}\in \R$ and $a_{i}\in \R$. Since the assumption $\max_{j\in [n], j\neq i^*} (x_j+a_j) = - \infty$ and $i\neq i^*$ gives us $x_i+a_i=-\infty$, we have $x_i=-\infty$, which leads to $d(x,x')=+\infty$. Therefore $\dist_H(x,\cH_a)=+\infty$ and~\Cref{dist-hyperplane} holds.
\end{proof}

The next lemma shows that the distance from a Hilbert's ball to any tropical hyperplane is bounded below by the radius of this ball. 
\begin{lemma}\label{lem-ball-hyperplane}
  For $a,b\in \proj(\rmax)^n$, suppose that the supports
  of $a$ and $b$ are not disjoint. Then,
  for all $r\geq 0$, we have
\begin{equation}
\dist_H (B(a,r),\cH_b) \geq r \enspace .
\end{equation}
\end{lemma} 

\begin{proof}
  Let $i^*\in\argmax_{i\in[n]}(a_i + b_i)$.
  Since the supports of $a,b$ are not disjoint,
  we have $a_{i^*}+b_{i^*}>-\infty$.
Define $x\in (\rmax)^n$ by $x_{i^*}=r+a_{i^*}$ and $x_i=a_i$ for all $i\neq i^*$. Then $x\in B(a,r)$, and for all $i\neq i^*$, $x_{i^*}+b_{i^*}= r+a_{i^*}+b_{i^*}\geq r+ a_i + b_i= r + x_i + b_i$. So by~\Cref{lem-dist-hyperplane} we deduce that $\dist_H(x,\cH_b) \geq r$, which implies that $\dist_H (B(a,r),\cH_b) \geq r$.
\end{proof}

\begin{lemma}\label{lem-dist-span}
  Suppose that $\mathcal{W}$ is a tropical cone in $(\rmax)^n$.
  Then, 
\begin{equation}
  \dist_H (\Sp (\cV),\mathcal{W}) =
  \dist_H (\cV,\mathcal{W}) \enspace .
\end{equation}
\end{lemma} 

\begin{proof}
  Consider an element $x\in \Sp (\cV)$, so that
  there exists a finite subset of points $(v^{(j)})_{j\in J}$ of $\cV$ and $(\alpha_j)_{j\in J}\in \R^{J}$, satisfying $x=\vee_{j\in J} (\alpha_j+v^{(j)})$.
  Take $\lambda>\dist_H (\cV,\mathcal{W})$.
  Then, for any $j \in J$, there exists $w^{(j)}\in \mathcal{W}$ such that $d(v^{(j)},w^{(j)})\leq \lambda$, and so, there are real numbers $\gamma^j,\beta^j$ such that
  $\gamma^j+ w^{(j)} \leq v^{(j)} \leq \beta^j+ w^{(j)}$,
  and $\beta^j-\gamma^j\leq \lambda$. After replacing
  $w^{(j)}$ by $\gamma^j+ w^{(j)}\in \mathcal{W}$, we may assume
  that $\gamma^j=0$. Then,
  $\vee_{j\in J} \alpha_j + w^{(j)} \leq x \leq \lambda +
  \vee_{j\in J} \alpha_j + w^{(j)}$, which entails
  that $\dist_H(x,\mathcal{W}) \leq \lambda$. Since this
  holds for all $\lambda>\dist_H(\cV,\mathcal{W})$,
  we deduce that $\dist_H(x,\mathcal{W})
  \leq \dist_H(\cV,\mathcal{W})$,
  and so,
$\dist_H (\Sp (\cV),\cW) \leq \dist_H (\cV,\cW)$.

The other inequality follows from $\cV \subset \Sp (\cV)$.
\end{proof}

The next lemma shows that the distance from the set $\cV$ to any tropical hyperplane is always greater than or equal to the radius of any Hilbert's ball included in the module $\Sp(\cV)$.

\begin{lemma}[Weak duality]\label{lem-weak-duality}
We have the following inequality
\begin{equation}\label{e-weak-duality}
\inrv{\cV} = \sup \{r \geq 0~|~\exists a \in \R^n, B(a,r) \subset  \Sp(\cV) \} \leq 
\inf_{b\in \proj(\rmax)^n} \dist_H (\cV,\cH_b) \enspace .
\end{equation}
\end{lemma} 

\begin{proof}
Let $a \in \R^n$ and $r\geq 0$ such that $B(a,r) \subset \Sp(\cV)$, and let $b \in \proj(\rmax)^n$. Since the supports of $a$ and $b$ are not disjoint, by~\Cref{lem-ball-hyperplane}, we have $r\leq \dist_H (B(a,r),\cH_b)$. Since $B(a,r) \subset \Sp(\cV)$, then $\dist_H (B(a,r),\cH_b)\leq \dist_H (\Sp(\cV),\cH_b)$. Therefore, by 
using~\Cref{lem-dist-span}, we conclude that $r \leq \dist_H (\cV,\cH_b)$.
\end{proof}

\begin{lemma}\label{lem-hyperplane}
For all $\lambda \in [-\infty,0]$ and $b \in \proj(\rmax)^n$, we have
\begin{equation}
T(b)\geq \lambda +b \Leftrightarrow \dist_H (\cV,\cH_b) \leq -\lambda \enspace .
\end{equation}
\end{lemma} 

\begin{proof}
  The equivalence is trivial if $\lambda=-\infty$,
  so, we suppose that $\lambda\in (-\infty,0]$.
    Suppose in addition that $T(b)\geq \lambda +b$, i.e., for any $i\in [n]$,
\[
\min_{k\in [p], (i,k)\in E}[-V_{ik}+\max_{j\in[n], j\neq i} (V_{jk}+b_j)] \geq \lambda +b_i.
\]
Then for any $i \in [n]$ and any $k \in [p]$,
\[
V_{ik}+b_i \leq \max_{j\in[n], j\neq i} (V_{jk}+b_j) - \lambda.
\] 
For each $k\in [p]$, by taking $i \in  \argmax_{j\in [n]}(V_{jk}+b_j)$ and using~\Cref{lem-dist-hyperplane},
we deduce that the distance from the column $V_{\cdot k}= v^{(k)}$ to the hyperplane $\cH_b$ is $\leq -\lambda$, which implies $\dist_H (\cV,\cH_b) \leq -\lambda$.

Now we suppose that $\dist_H (\cV,\cH_b) \leq -\lambda$. For $k\in [p]$ and $i\in [n]$,
if $i\not\in \argmax_{j\in [n]}(V_{jk}+b_j)$, then
\[
V_{ik}+b_i \leq \max_{j\in[n], j\neq i} (V_{jk}+b_j)\leq \max_{j\in[n], j\neq i} (V_{jk}+b_j) -\lambda.
\]
Otherwise, if $i\in \argmax_{j\in [n]}(V_{jk}+b_j)$, then knowing that $\dist_H (v^{(k)},\cH_b) \leq -\lambda$ and using~\Cref{lem-dist-hyperplane}, we get $V_{ik}+b_i \leq \max_{j\in[n], j\neq i} (V_{jk}+b_j) -\lambda$.
Therefore, we deduce that for all $i \in [n]$ and $k \in [p]$,
\[
V_{ik}+b_i \leq \max_{j\in[n], j\neq i} (V_{jk}+b_j) - \lambda.
\]
Thus for all $i\in [n]$,
\[
\min_{k\in [p], (i,k)\in E}[-V_{ik}+\max_{j\in[n], j\neq i} (V_{jk}+b_j)] \geq \lambda +b_i,
\]
namely $T(b)\geq \lambda +b$.
\end{proof}

The following theorem presents a strong duality result between a best tropical hyperplane approximation of a set $\cV$ of points and the largest inner balls that its module $\Sp (\cV)$ contains.

\begin{theorem}[Strong duality]\label{th-strong-duality}
  We have 
\begin{equation}\label{e-duality}
\min_{b\in \proj(\rmax)^n} \dist_H (\cV,\cH_b) =  \inrv{\cV} = \sup \{r \geq 0\mid \exists a \in \R^n, B(a,r) \subset  \Sp(\cV) \} .
\end{equation}
The minimum is achieved by any vector $b \in \proj(\rmax)^n$
such that $T(b)\geq \rho(T) + b$. Moreover, if $\inrv{\cV}$ is finite, the supremum is achieved by a ball $B(-c,\inrv{\cV})$ 
where $c \in \R^n$ is any vector such that $T(c)\leq \rho(T) +c$. 
\end{theorem}

\begin{proof}%
  \Cref{th-inradius} entails that $\inrv{\cV}=-\rho(T)$
  and that the last assertion of the theorem holds.
Moreover, the existence of a vector $b\in \proj(\rmax)^n$ such
$T(b)\geq \rho(T)+ b$ follows from \Cref{th-duality}.
Then, by~\Cref{lem-hyperplane}, we have $\dist_H (\cV,\cH_b) \leq \inrv{\cV}$,
which combined with the weak duality property~\eqref{e-weak-duality} implies
that the equality holds in~\eqref{e-weak-duality}, and that
$b$ such that $T(b)\geq \rho(T)+b$ achieves the minimum in~\eqref{e-duality}.
\end{proof}

The following lemma allows us to bound from below the value of the tropical
linear regression problem by looking at points in the sectors
of a hyperplane $\cH_a$.

\begin{lemma}\label{lem-hyperplane-ball}
If $a\in \R^n$ and $r\in [0,+\infty]$ are such that 
\[
\forall i \in [n], \exists k \in [p], \quad  v^{(k)} \in S_i(a) \; \text{\and} \; \dist_H(v^{(k)}, \cH_a) \geq r  \enspace,
\]
then $B(-a,r) \subset \Col (V)$ and $\min_{b\in \proj(\rmax)^n} \dist_H (\cV,\cH_b) \geq r$ .
\end{lemma}
\begin{proof}
If $r=+\infty$, then for any $i\in [n]$, there is some $\sigma_i \in [p]$ such that $v^{(\sigma_i)} \in S_i(a)$ and $\dist_H(v^{(\sigma_i)}, \cH_a) = +\infty$. Since $a$ is finite, for any $i \in[n]$, we have $v^{(\sigma_i)}_i\in \R$ and $v^{(\sigma(i))}_j=-\infty$ for any $j \neq i$. We deduce that $\Sp(v^{(\sigma_1)}, \cdots, v^{(\sigma_n)})=\rmax$. Then $\Col (V)=\rmax$, and so  $B(-a,+\infty)= \rmax \subset \Col (V)$.

Now we consider $r\in [0,+\infty)$. For $i\in [n]$, by the assumption of this lemma, there exists $k\in [p]$ such that $v^{(k)} \in S_i(a)$ and $\dist_H(v^{(k)}, \cH_a) \geq r $. 
Hence, by using~\eqref{e-sectors} and~\Cref{lem-dist-hyperplane}, we deduce that the column $V_{\cdot k}=v^{(k)}$ satisfies $V_{ik}+a_i \geq r + \max_{j\neq i} (V_{jk}+a_j)$. 
Therefore, we have $-V_{ik} +\max_{j\neq i} (V_{jk}+a_j) \leq -r +a_i$, which implies for any $i \in [n]$, 
\[
T_i(a)=\inf_{l\in [p], (i,l)\in E} \Big[-V_{il}+ \max_{j\in [n],j\neq i} (V_{jl}+a_j) \Big] \leq -r + a_i,
\]
i.e., $T(a)\leq -r +a$. Therefore, by~\Cref{lem-ball}, we deduce that $B(-a,r) \subset \Col (V)$.

Finally by~\Cref{th-strong-duality}, we have
\[
\min_{b\in \proj(\rmax)^n} \dist_H (\cV,\cH_b) = \inrv{\cV} = \sup \{s \geq 0\mid \exists x \in \R^n, B(x,s) \subseteq \Sp(\cV) \} \geq r.
\]
\end{proof}

Given a hyperplane $\cH_b$, we call {\em witness point} of $\cH_b$ any point $p$ in $\cV$ such that the distance from $p$ to the hyperplane $\cH_b$ equals the distance from the set $\cV$ to this hyperplane.

\begin{theorem}[Optimality certificates]\label{th-witness}
Let $a\in \R^n$, then the following assertions are equivalent:
\begin{enumerate}
  \item $T(a)=\rho(T) +a$;
  \item The hyperplane $\cH_a$ admits a witness point in each sector,
    meaning that
    $\forall i \in [n], \exists k \in [p], \;  v^{(k)} \in S_i(a)  \text{ \and }  \dist_H(v^{(k)}, \cH_a) = 
\dist_H(\cV, \cH_a)$.
\end{enumerate}
Moreover, if these assertions hold,
then, $\rho(T)=-\dist_H(\cV,\cH_a)$,
$\cH_a$ is an optimal
  solution of the tropical linear regression problem,
  and $B(-a,\dist_H(\cV,\cH_a))$ is a Hilbert's ball
  of maximal radius included in $\Sp(\cV)$.
\end{theorem}

\begin{proof}
If $a\in \R^n$ satisfies $T(a)=\rho(T) +a=-\inrv{\cV} +a$, then by~\Cref{th-strong-duality}, $\cH_a$ is optimal in~\Cref{e-duality}, i.e., $\dist_H(\cV, \cH_a)=\inrv{\cV}$, and for all $i\in [n]$ we have
\[
\min_{k\in [p], (i,k)\in E}[-V_{ik}+\max_{j\in[n], j\neq i} (V_{jk}+a_j)]=-\inrv{\cV}+a_i.
\]
Then for all $i\in [n]$, there exists $k\in [p]$ such that
\[
-V_{ik}+\max_{j\in[n], j\neq i} (V_{jk}+a_j)=-\inrv{\cV}+a_i,
\]
i.e., $V_{ik}+a_i = \inrv{\cV}+ \max_{j\in[n], j\neq i} (V_{jk}+a_j)$. This implies that $v^{(k)} = V_{\cdot k} \in S_i(a)$, and also by~\Cref{lem-dist-hyperplane}, that $\dist_H(v^{(k)}, \cH_a)= \inrv{\cV}=\dist_H(\cV, \cH_a)$.

Now, we suppose that we have assertion $(2)$. By~\Cref{lem-hyperplane-ball}, we have
\[
\min_{b\in \proj(\rmax)^n} \dist_H (\cV,\cH_b) \geq \dist_H(\cV, \cH_a),
\]
which means that $\cH_a$ achieves the minimum in~\eqref{e-duality}, so that $\dist_H(\cV, \cH_a)=\inrv{\cV}$.
Hence, $\forall k \in [p], \dist_H(v^{(k)}, \cH_a)\leq \inrv{\cV}$, so that by~\Cref{lem-dist-hyperplane} we have $\forall k \in [p], \forall i \in [n]$, 
\[
V_{ik}+a_i \leq \inrv{\cV}+ \max_{j\in[n], j\neq i} (V_{jk}+a_j).
\]
Therefore, we obtain
\begin{equation}\label{e-a}
\forall i \in [n], \forall k \in [p],(i,k)\in E; -\inrv{\cV}+a_i \leq -V_{ik} + \max_{j\in[n], j\neq i} (V_{jk}+a_j) \enspace.
\end{equation}
Assertion $(2)$ also implies $\forall i \in [n], \exists k \in [p], V_{ik}+a_i \geq \max_{j\in[n], j\neq i} (V_{jk}+a_j)$ and  $\dist_H(v^{(k)}, \cH_a)= \inrv{\cV}$, with $V_{ik}\neq -\infty$ because $a\in \R^n$ and $V_{\cdot k}\neq \bot$.  This means, by~\Cref{lem-dist-hyperplane}, that $\forall i \in [n], \exists k \in [p],(i,k)\in E, V_{ik}+a_i = \inrv{\cV}+ \max_{j\in[n], j\neq i} (V_{jk}+a_j)$. Then
\begin{equation}\label{e-b}
\forall i \in [n], \exists k \in [p],(i,k)\in E; -\inrv{\cV}+a_i = -V_{ik} + \max_{j\in[n], j\neq i} (V_{jk}+a_j) \enspace .
\end{equation}
From~\eqref{e-a} and~\eqref{e-b}, we conclude that
\[
\forall i \in [n], -\inrv{\cV}+a_i= \min_{k\in[p],(i,k)\in E} [-V_{ik} + \max_{j\in[n], j\neq i} (V_{jk}+a_j)] = T_i(a).
\]
Therefore $T(a)=\rho(T)  + a$.

The final part of the theorem follows from~\Cref{th-strong-duality}.
\end{proof}
\begin{remark}\label{rk-rk}
  When $T(a)=\rho(T)+a$ and $a\in\R^n$, \Cref{th-witness} and~\Cref{th-strong-duality} entail the following remarkable property: there is an optimal Hilbert's
  ball whose center coincides with the apex of an optimal
  regression hyperplane. This property is illustrated in~\Cref{fig-Multi-Hyperplanes} below.
\end{remark}
\begin{remark}
  The situation in which $T(a)=\bot$ holds for some finite vector
  $a$ (or equivalently, for all finite vectors $a$) is degenerate.
  Indeed, we observe from the proof of~\Cref{th-witness}
  that $T(a)=\bot$ for some finite vector $a$
  if and only if, for all $i\in [n]$,
  there is a vector $v^{(k)}$ such that $v^{(k)}_i$
  is finite and all $v^{(k)}_j$ with $j\neq i$ are $-\infty$.
  Then, $V$ contains a $n\times n$ diagonal submatrix,
  and so, $\Col(V)=\Sp(\cV)=(\rmax)^n$. 
  \end{remark}
We next exhibit a situation in which the existence
of a finite eigenvector, required
to apply~\Cref{th-witness}, is guaranteed.
\begin{proposition}\label{prop-finite}
  Suppose that all the vectors $v\in \cV$ have finite entries.
  Then, the operator $T$ has a finite eigenvector $a$.
\end{proposition}
\begin{proof}
 Theorems~9 and~13 of~\cite{sgjg04}
 imply that  an order preserving and additively
homogeneous map $T:\R^n\to \R^n$ has a finite eigenvector if
 the recession function $\hat{T}(x):=\lim_{s\to \infty} s^{-1}T(sx)$
 has only fixed points on the diagonal. When
 the matrix $V$ is finite, considering $T:=T_V$,
 we have $\hat{T}_i(x)=\max_{j\in[n], j\neq i}x_j$,
for all $i\in[n]$,
so the latter condition is trivially satisfied.
This entails that there exists a vector $a\in\R^n$
such that $T(a)=\rho(T)+a$.
\end{proof}
A more general condition, involving the notion of dominions,
is given in \Cref{sec-appendix}.
\if{
\begin{corollary}\label{coro-witness}
Let $a \in \R^n$ such that $T(a)=-\inrv{\cV}  + a$. The tropical hyperplane $\cH_a$ is optimal in~\Cref{e-duality},%
and it has at least a witness point in each sector, i.e.
\begin{equation}
\forall i \in [n], \exists k \in [p] ; \; v^{(k)} \in S_i(a)  \; ,  \; \dist_H(v^{(k)}, \cH_a) =  \inrv{\cV} \enspace .
\end{equation}
\end{corollary}
\begin{proof}
We have $T(a)=-\inrv{\cV}  + a$, then by~\Cref{th-strong-duality}, $\cH_a$ is optimal in~\Cref{e-duality}, and $\forall i\in [n]$, $\min_{k\in [p]}[-V_{ik}+\max_{j\in[n], j\neq i} (V_{jk}+a_j)]=-\inrv{\cV}+a_i$. 
Then for all $i\in [n]$, there exists $k\in [p]$ such that $-V_{ik}+\max_{j\in[n], j\neq i} (V_{jk}+a_j)=-\inrv{\cV}+a_i$, i.e. $V_{ik}+a_i = \inrv{\cV}+ \max_{j\in[n], j\neq i} (V_{jk}+a_j)$. This implies that $v^{(k)}= V_{\cdot k} \in S_i(a)$ and also, by~\Cref{lem-dist-hyperplane}, that $\dist_H(v^{(k)}, \cH_a)= \inrv{\cV}$.
\end{proof}

The following theorem characterizes the hyperplanes $\cH_a$ associated to a finite vector $a$ that is and eigenvector of the operator $T$, i.e. $T(a)=-\inrv{\cV}  + a$.

\begin{theorem}\label{th-witnesses}
A vector $a\in \R^n$ is optimal in~\Cref{e-duality} and the hyperplane $\cH_a$ has at least a witness point in each sector $(S_i(a))_{i\in [n]}$ if and only if $T(a)=-\inrv{\cV}  + a$ .
\end{theorem}

\begin{proof}
The ``If" implication is given by~\Cref{coro-witness}.

Now, if $a\in \R^n$ is optimal in~\Cref{e-duality} then by~\Cref{th-strong-duality}, $\dist_H(\cV, \cH_a)=\inrv{\cV}$.
Then $\forall k \in [p], \dist_H(v^{(k)}, \cH_a)\leq \inrv{\cV}$, and by~\Cref{lem-dist-hyperplane} we deduce that $\forall k \in [p], \forall i \in [n], V_{ik}+a_i \leq \inrv{\cV}+ \max_{j\in[n], j\neq i} (V_{jk}+a_j)$. Therefore, $\forall i \in [n], \forall k \in [p], -\inrv{\cV}+a_i \leq -V_{ik} + \max_{j\in[n], j\neq i} (V_{jk}+a_j)$.

The existence of a witness point in each sector of $\cH_a$, implies that $\forall i \in [n], \exists k \in [p], \dist_H(v^{(k)}, \cH_a)= \inrv{\cV}$ and $V_{ik}+a_i \geq \max_{j\in[n], j\neq i} (V_{jk}+a_j)$. This means, by~\Cref{lem-dist-hyperplane}, that $\forall i \in [n], \exists k \in [p], V_{ik}+a_i = \inrv{\cV}+ \max_{j\in[n], j\neq i} (V_{jk}+a_j)$. Then $\forall i \in [n], \exists k \in [p], -\inrv{\cV}+a_i = -V_{ik} + \max_{j\in[n], j\neq i} (V_{jk}+a_j)$.

We conclude that $\forall i \in [n], -\inrv{\cV}+a_i= \min_{k\in[p]} [-V_{ik} + \max_{j\in[n], j\neq i} (V_{jk}+a_j)] = T_i(a)$. Therefore $T(a)=-\inrv{\cV}  + a$.
\end{proof}
}\fi

The following proposition shows that we can determine witness points from a policy $\sigma:[n] \mapsto [p]$, that satisfies $T(a)=T^\sigma(a)$ where $a$ is a finite eigenvector of the operator $T$. For an illustration of this lemma see~\Cref{fig-inner-ball}.
\begin{proposition}
Let $a \in \R^n$ such that $T(a)=-\inrv{\cV}  + a$, and $\sigma:[n] \mapsto [p]$ a map, such that $\forall i \in [n], (i,\sigma(i)) \in E$. We have $T(a)=T^\sigma(a)$ if and only if for all $i\in [n]$, $V_{\cdot \sigma(i)}$ is a witness point of $\cH_a$ that belongs to the sector $S_i(a)$.
\end{proposition}

\begin{proof}
If $T(a)=T^\sigma(a)$, then $T^\sigma (a)=-\inrv{\cV}  + a$. Therefore, we have for all $i\in [n]$,  $-V_{i \sigma(i)} + \max_{j\neq i} (V_{j \sigma(i)} +a_j) = -\inrv{\cV}  + a_i$, i.e. $V_{i \sigma(i)}+a_i = \inrv{\cV} + \max_{j\neq i} (V_{j \sigma(i)} +a_j)$, which means that $V_{\cdot \sigma(i)} \in S_i(a)$ and, by~\Cref{lem-dist-hyperplane}, that $\dist_H(V_{\cdot \sigma(i)}, \cH_a) = \inrv{\cV}$, i.e for all $i\in [n]$, $V_{\cdot \sigma(i)}$ is a witness point in the sector $S_i(a)$.

Conversely, if for all $i\in [n]$, $V_{\cdot \sigma(i)}$ is a witness point in the sector $S_i(a)$. Let $i\in [n]$, we have then $V_{i \sigma(i)}+a_i = \inrv{\cV} + \max_{j\neq i} (V_{j \sigma(i)} +a_j)$, i.e. $-V_{i \sigma(i)} + \max_{j\neq i} (V_{j \sigma(i)} +a_j) = -\inrv{\cV}  + a_i$.
We know that for all $k\in [p]$, $\dist_H(V_{\cdot k}, \cH_a) \leq \inrv{\cV}$, then by~\Cref{lem-dist-hyperplane}, $V_{i k}+a_i \leq \inrv{\cV} + \max_{j\neq i} (V_{j k} +a_j)$, i.e. $-V_{i k} + \max_{j\neq i} (V_{j k} +a_j) \geq -\inrv{\cV}  + a_i$.
Therefore, $T_i(a) = \inf_{k\in [p], (i,k)\in E} \Big[-V_{ik}+ \max_{j\in [n],j\neq i} (V_{jk}+a_j) \Big]  = -V_{i \sigma(i)} + \max_{j\neq i} (V_{j \sigma(i)} +a_j) = T^\sigma_i (a)$.
\end{proof}

\begin{figure}[htbp]
\begin{center}


\def\coord#1#2#3{{-sqrt(3)/2*(#1-#2)} ,{ -(1/2)*#1 - (1/2)*#2 + #3}}
\begin{tikzpicture}[scale=0.95]

\coordinate (v1) at (\coord{-3}{0}{-1});
\coordinate (v2) at (\coord{0}{-3}{-1});
\coordinate (v3) at (\coord{0}{0}{-4});
\coordinate (v4) at (\coord{1}{0}{-2});
\coordinate (v5) at (\coord{1}{-1}{-1});
\coordinate (v6) at (\coord{-1}{1}{-1});
\coordinate (v7) at (\coord{0}{1}{-2});
\coordinate (v8) at (\coord{0}{-1}{0});
\coordinate (v9) at (\coord{-1}{0}{0});

\coordinate (v1p) at (\coord{-5}{0}{-1});
\coordinate (v2p) at (\coord{0}{-5}{-1});
\coordinate (v3p) at (\coord{0}{0}{-5});

\coordinate (f1) at (\coord{-2}{0}{-1});
\coordinate (f2) at (\coord{0}{-2}{-1});
\coordinate (f3) at (\coord{0}{0}{-3});

\coordinate (b1) at (\coord{0}{0}{0});
\coordinate (b2) at (\coord{0}{1}{0});
\coordinate (b3) at (\coord{0}{1}{-1});
\coordinate (b4) at (\coord{0}{0}{-2});
\coordinate (b5) at (\coord{1}{0}{-1});
\coordinate (b6) at (\coord{1}{0}{0});

\fill[red!30,opacity=0.7] (v1p) -- (\coord{0}{0}{-1}) -- (v2p)--cycle;
\fill[blue!30,opacity=0.7] (v1p) -- (\coord{0}{0}{-1}) -- (v3p)--cycle;
\fill[green!30,opacity=0.7] (v2p) -- (\coord{0}{0}{-1}) -- (v3p)--cycle;

\filldraw[black,draw=black,opacity=0.8,very thick] (v1) -- (f1) -- (v9) -- (b1) -- (v8) -- (f2) -- (v2) -- (f2) -- (v5) -- (b5) -- (v4) -- (f3) -- (v3) -- (f3) -- (v7) -- (b3) -- (v6) -- (f1) -- (v1) -- cycle;
\filldraw[gray] (f1) -- (v9) -- (b1) -- (v8) -- (f2) -- (v5) -- (b5) -- (v4) -- (f3) -- (v7) -- (b3) -- (v6) -- (f1) -- cycle;

\filldraw[blue!80,draw=blue!80,opacity=0.7, thick] (b1) -- (b2) -- (b3) -- (b4) -- (b5) -- (b6) -- (b1) -- cycle;

\filldraw[black] (v1) circle (1.5pt) node[above] {$V_{\cdot 1}$};
\filldraw[black] (v2) circle (1.5pt) node[above] {$V_{\cdot 2}$};
\filldraw[black] (v3) circle (1.5pt) node[below,right]{$V_{\cdot 3}$};
\filldraw[red] (v4) circle (1.5pt) node[below,left,black] {$V_{\cdot 4}$};
\filldraw[red] (v5) circle (1.5pt) node[below,left,black] {$V_{\cdot 5}$};
\filldraw[red] (v6) circle (1.5pt) node[below,right,black] {$V_{\cdot 6}$};
\filldraw[red] (v7) circle (1.5pt) node[below,right,black] {$V_{\cdot 7}$};
\filldraw[red] (v8) circle (1.5pt) node[above,black] {$V_{\cdot 8}$};
\filldraw[red] (v9) circle (1.5pt) node[above,black] {$V_{\cdot 9}$};

\draw[dashed,->] (\coord{0}{0}{0}) -- (\coord{5}{0}{0}) node[above] {$x_1$};
\draw[dashed,->] (\coord{0}{0}{0}) -- (\coord{0}{5}{0}) node[above] {$x_2$};
\draw[dashed,->] (\coord{0}{0}{0}) -- (\coord{0}{0}{2.5}) node[above, right] {$x_3$};
\filldraw[black] (\coord{0}{0}{0}) circle (1.5pt) node[below,right] {0};

\draw[very thick,red] (\coord{0}{0}{-1}) -- (v1p);
\draw[very thick,red] (\coord{0}{0}{-1}) -- (v2p);
\draw[very thick,red] (\coord{0}{0}{-1}) -- (v3p);
\filldraw[black] (\coord{0}{0}{-1}) circle (1.5pt) node[below,right] {$-a$};
\draw node[above] at (\coord{-1/2}{0}{-1}) {$\cH_a$};

\end{tikzpicture}
\end{center}
\caption{The inner ball of a column space $\Col (V)$ and the linear regression of the columns of $V$.}
\label{fig-witnesses}
\end{figure}
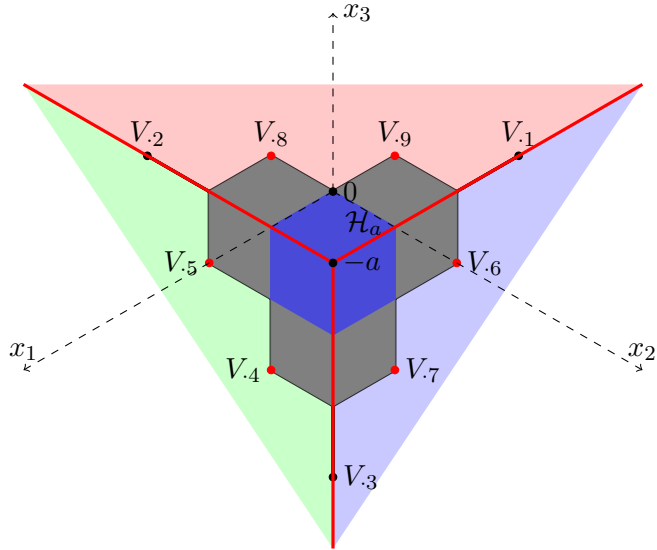

We now formalize the tropical linear regression problem:
\begin{problem}[Tropical linear regression]\label{pb-regression}
  Input: a finite set of vectors $\cV\subset \zmax^{n}$.
Goal: compute
the infimum of the one-sided Hausdorff distance of $\cV$
to a tropical hyperplane, i.e.,
  the value of the optimization problem~\eqref{e-regression}.
\end{problem}
\begin{corollary}\label{cor-polytime}
  The tropical linear regression problem (\Cref{pb-regression})
  is polynomial time Turing-equivalent to mean payoff games
  (\Cref{pb-1}).
\end{corollary}
\begin{proof}
  This follows from the strong duality theorem (\Cref{th-strong-duality}) and~\Cref{coro-equiv}. 
\end{proof}
\begin{corollary}\label{cor-dualcenter}
  Computing an optimal regression hyperplane $\cH_a$ in~\eqref{e-regression}, given a finite set of vectors $\cV\subset \zmax^{n}$, polynomially Turing-reduces to mean payoff games. 
\end{corollary}
\begin{proof}
  By~\Cref{th-strong-duality}, we need to find a vector $a$
  such that $T(a)\geq \rho(T)+a$. Arguing as in the proof
  of~\Cref{cor-center},
  but exchanging the roles of Player Max and Min,
  we end up with an optimal policy $\tau$ of Player Max.
  Then, it suffices to find a vector $a\in (\R\cup\{-\infty\})^n$,
  $a\neq \bot$, such that
  $^\tau T(a) \geq \rho(^\tau T) + a$. Still arguing
  as in the proof of~\Cref{cor-center}, we are
  reduced to a problem of tropical (min-plus instead max-plus)
  spectral theory, which again reduces to a shortest path problem.
  \end{proof}
In~\Cref{fig-witnesses}, we consider the same matrix $V$ as in~\Cref{fig-inner-ball}. The~\Cref{fig-witnesses} shows the witness points in each of the sectors defined by the hyperplane $\cH_a$ where $a=(0,0,1)^\top$ satisfies $T(a)=\lambda  + a$ with $\lambda=-1$.
In this example, we have two witness points in each sector: $V_{\cdot 4}$ and $V_{\cdot 5}$ are the witness points in the sector $S_1(a)$ (in green), $V_{\cdot 6}$ and $V_{\cdot 7}$ are the witness points in the sector $S_2(a)$ (in blue) and $V_{\cdot 8}$ and $V_{\cdot 9}$ are the witness points in the sector $S_3(a)$ (in red).

In~\Cref{fig-Multi-Hyperplanes}, we consider the following matrix $U\in \R^{3 \times 4}$:
\[
U=\left(\begin{array}{cccc}
  -1 & 0  &  1 &  0  \\
   0 & -1 &  0  & 1  \\
   0 & 0  & -2 & -2  
\end{array}
\right)
\]

The operator associated to $U$ is the following map $T:(\rmax)^n \to (\rmax)^n$ :
\[
\resizebox{.9\hsize}{!}{$T \begin{pmatrix}
x_1 \\ x_2 \\ x_3
\end{pmatrix}
 = \begin{pmatrix}
  \min[1+\max(x_2,x_3), \max(-1+x_2,x_3), -1+\max(x_2,-2 + x_3), \max(1+x_2,-2+x_3) ]  \\
  \min[\max(-1+x_1,x_3), 1+\max(x_1,x_3), \max(1+x_1,-2 + x_3),-1+ \max(x_1,-2+x_3) ]  \\
  \min[\max(-1+x_1,x_2), 1+\max(x_1,-1+x_2), 2+\max(1+x_1, x_2), 2+ \max(x_1,1+x_2) ]    
\end{pmatrix} \enspace .$}
\]

We verify easily that $\lambda=-1$ and $a=(0,0,1)^\top$ satisfy $T(a)=\lambda +a$, so that the inner radius of $\Col (U)$ is $\inrv{\mathcal{U}} =1$. In this example, other hyperplanes like $\cH_b$ and $\cH_c$, with $b=(0,0,-1)^\top$ and $c=(0,0,-\infty)^\top$, are also optimal solutions of the tropical linear regression
problem,
but $\cH_a$ is the only hyperplane such that $a$ is a finite eigenvector of the operator $T$ and, hence, that satisfies also $B(-a,1) \subset \Col (U)$. 

\begin{figure}[htbp]
\begin{center}


\def\coord#1#2#3{{-sqrt(3)/2*(#1-#2)} ,{ -(1/2)*#1 - (1/2)*#2 + #3}}
\begin{tikzpicture}[scale=0.95]

\coordinate (v1) at (\coord{-1}{0}{0});
\coordinate (v2) at (\coord{0}{-1}{0});
\coordinate (v3) at (\coord{1}{0}{-2});
\coordinate (v4) at (\coord{0}{1}{-2});

\coordinate (Ha1) at (\coord{-5}{0}{-1});
\coordinate (Ha2) at (\coord{0}{-5}{-1});
\coordinate (Ha3) at (\coord{0}{0}{-4.5});

\coordinate (Hb1) at (\coord{-4.5}{0}{0.5});
\coordinate (Hb2) at (\coord{0}{-4.5}{0.5});
\coordinate (Hb3) at (\coord{0}{0}{-4});

\coordinate (f) at (\coord{0}{0}{-3});

\coordinate (b1) at (\coord{0}{0}{0});
\coordinate (b2) at (\coord{0}{1}{0});
\coordinate (b3) at (\coord{0}{1}{-1});
\coordinate (b4) at (\coord{0}{0}{-2});
\coordinate (b5) at (\coord{1}{0}{-1});
\coordinate (b6) at (\coord{1}{0}{0});

\fill[red!30,opacity=0.7] (Ha1) -- (\coord{0}{0}{-1}) -- (Ha2)--cycle;
\fill[blue!30,opacity=0.7] (Ha1) -- (\coord{0}{0}{-1}) -- (Ha3)--cycle;
\fill[green!30,opacity=0.7] (Ha2) -- (\coord{0}{0}{-1}) -- (Ha3)--cycle;

\filldraw[black,draw=black,opacity=0.8,very thick] (v1) -- (b1) -- (v2) -- (v3) -- (f) -- (v4) -- (v1) -- cycle;
\filldraw[gray!70] (v1) -- (b1) -- (v2) -- (v3) -- (f) -- (v4) -- (v1) -- cycle;

\filldraw[black!60,draw=black!60,opacity=0.7, thick] (b1) -- (b2) -- (b3) -- (b4) -- (b5) -- (b6) -- (b1) -- cycle;

\filldraw[red] (v1) circle (1.5pt) node[below right,black] {$U_{\cdot 1}$};
\filldraw[red] (v2) circle (1.5pt) node[below left,black] {$U_{\cdot 2}$};
\filldraw[red] (v3) circle (1.5pt) node[below,left,black] {$U_{\cdot 3}$};
\filldraw[red] (v4) circle (1.5pt) node[below,right,black] {$U_{\cdot 4}$};

\draw[dashed,->] (\coord{0}{0}{0}) -- (\coord{5}{0}{0}) node[above] {$x_1$};
\draw[dashed,->] (\coord{0}{0}{0}) -- (\coord{0}{5}{0}) node[above] {$x_2$};
\draw[dashed,->] (\coord{0}{0}{0}) -- (\coord{0}{0}{3.5}) node[above, right] {$x_3$};
\filldraw[black] (\coord{0}{0}{0}) circle (1.5pt) node[below,right] {0};

\draw[dashed] (b1) -- (\coord{-3}{0}{0});
\draw[dashed] (b1) -- (\coord{0}{-3}{0});
\draw[dashed] (b3) -- (\coord{0}{1}{-3.5});
\draw[dashed] (b3) -- (\coord{-4}{0}{-2});
\draw[dashed] (b5) -- (\coord{1}{0}{-3.5});
\draw[dashed] (b5) -- (\coord{0}{-4}{-2});

\draw[black!30!green,line width=1.2mm] (\coord{0}{0}{-4}) -- (\coord{0}{0}{3.2});
\draw node[right] at (\coord{0}{0}{2.7}) {$\cH_c$};

\draw[very thick,red] (\coord{0}{0}{-1}) -- (Ha1);
\draw[very thick,red] (\coord{0}{0}{-1}) -- (Ha2);
\draw[very thick,red] (\coord{0}{0}{-1}) -- (Ha3);
\filldraw[black] (\coord{0}{0}{-1}) circle (1.5pt) node[below,right] {$-a$};
\draw node[above] at (\coord{-1/2}{0}{-1}) {$\cH_a$};

\draw[black,line width=0.3mm] (\coord{0}{0}{0.5}) -- (Hb1);
\draw[black,line width=0.3mm] (\coord{0}{0}{0.5}) -- (Hb2);
\draw[black,line width=0.3mm] (\coord{0}{0}{0.5}) -- (Hb3);
\filldraw[black] (\coord{0}{0}{0.5}) circle (1.5pt) node[above] {$-b$};
\draw node[above] at (\coord{-1/2}{0}{0.5}) {$\cH_b$};

\end{tikzpicture}
\end{center}
\caption{A column space $\Col (U)$ (light and dark gray regions) with multiple hyperplanes that are optimal solutions of the tropical linear regression problem, and multiple inner balls of maximal radius, but a unique optimal hyperplane with witness points in each sector, corresponding to the finite eigenvector $a=(0,0,1)^\top$ of $T$ and to the inner ball in dark gray.}
\label{fig-Multi-Hyperplanes}
\end{figure}
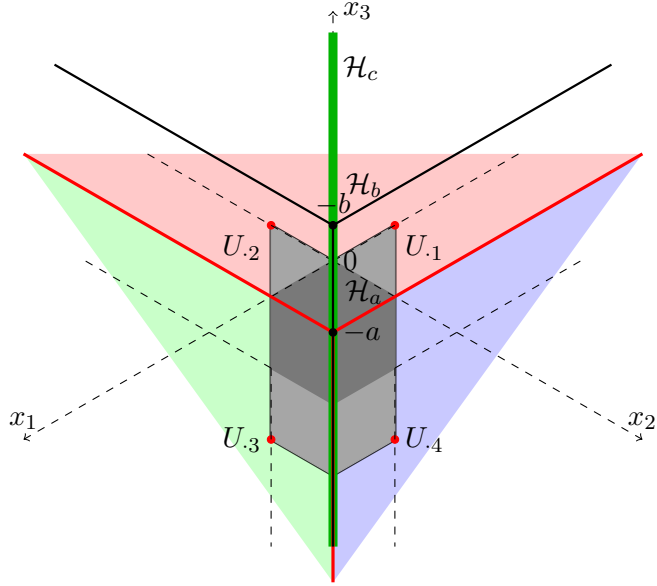

\section{Tropical linear regression with sign or type patterns}\label{sec-signed}

Here, we study several variants of the tropical linear regression
problem, which can also be solved by the present technique of reduction
to a mean payoff game. The second of these variants (with ``types'')
will arise in the economic application of~\Cref{sec-example}.

\subsection{Tropical linear regression with signs}
Given $I,J \subset [n]$ such that $I,J\neq \emptyset$, $I\cup J=[n]$ and $I\cap J = \emptyset$ and $a\in \proj(\rmax)^n$, we define the {\em signed tropical hyperplane} of type $(I,J)$:
\begin{equation}\label{signed-hyperplane}
\cH^{IJ}_a : =\{ x\in (\rmax)^n \mid \max_{i\in I} (a_i+x_i) = \max_{j\in J} (a_j+x_j)\}
\end{equation}

Given a set $\cV\subset (\rmax)^n$, of cardinality $|\cV| =p $,
the {\em signed tropical linear regression problem} of type $(I,J)$
consists in finding the best approximation of $\cV$ by a signed hyperplane of type $(I,J)$: %

\begin{equation}
\min_{a\in \proj(\rmax)^n} \dist_H(\cV , \cH^{IJ}_a) \enspace .
\end{equation}

Let $M$ be a closed tropical cone of $(\rmax)^n$ and $x\in (\rmax)^n$. The projection $P_M(x)$ of the point $x$ onto $M$~\cite{cgq02} is defined by:
\begin{equation}\label{e-proj}
P_M(x):= \max \{z\in M \mid z\leq x \} \enspace . 
\end{equation}

The Hilbert's distance from $x$ to $M$ is achievd
by the projection $P_M(x)$.

\begin{theorem}[{\cite{cgq02}}]\label{th-proj}
Given a closed tropical semimodule $M\subset (\rmax)^n$ and $x\in (\rmax)^n$, we have:
\[
\dist_H(x,M) = d(x,P_M(x)) \enspace .
\]
\end{theorem}

The following lemma identifies the projection of a point $x\in \proj(\rmax)^n$ onto a signed tropical hyperplane $\cH^{IJ}_a$. 

\begin{lemma}\label{lem-proj-sh}
Let $x,a \in \proj(\rmax)^n$ and $K=\supp a$. The projection $P_{\cH^{IJ}_a}(x)$ of $x$ onto $\cH^{IJ}_a$ is given by:
\begin{equation}\label{e-proj-sh}
[P_{\cH^{IJ}_a}(x)]_l=
\left\{
\begin{array}{lll}
x_l & , \text{ for }   l\in K^c \\
\min\{ x_l   , -a_l + \max_{j\in J} (a_j+x_j)  \} &, \text{ for }   l\in I \cap K \\
\min \{ x_l   , -a_l + \max_{i\in I} (a_i+x_i)  \} &,  \text{ for }   l\in J \cap K
\end{array}
\right.
\end{equation}
where $K^c$ denotes the complementary of $K$ in $[n]$.
\end{lemma}

\begin{proof}
Denote the right hand side vector of~\eqref{e-proj-sh} by $\tilde{x}$.
From~\eqref{e-proj}, we have $P_{\cH^{IJ}_a}(x)=\max \{z\in \cH^{IJ}_a \mid z\leq x \}$. Let $z\in \cH^{IJ}_a $ such that $z\leq x$. We will prove that $z\leq \tilde{x}$.
Let $l\in I$, if $l\in I\cap K^c$, we have right away that $z_l\leq x_l = \tilde{x}_l$. Now if $l\in I \cap K$, knowing that $z\in \cH^{IJ}_a$ and using~\eqref{signed-hyperplane}, we have $a_l+z_l \leq \max_{i\in I} (a_i+z_i) = \max_{j\in J} (a_j+z_j)\leq \max_{j\in J} (a_j+x_j)$. Then, $z_l\leq -a_l + \max_{j\in J} (a_j+x_j)$.
We know also that $z_l\leq x_l$, then $z_l \leq \tilde{x}_l$. 
Similarly the inequality $z_l \leq \tilde{x}_l$ can also be proved for all $l\in J$. 
Therefore, for all $z\in \cH^{IJ}_a $, if $z\leq x$ then $z\leq \tilde{x}$. 
Using~\eqref{e-proj}, it suffices now to prove that $\tilde{x} \in \cH^{IJ}_a$. 
Indeed, %
$\max_{i\in I} (a_i+\tilde{x}_i) = \max_{i\in I\cap K} (a_i+\tilde{x}_i) =
\max_{i\in I\cap K}\{ \min (a_i + x_i, \max_{j\in J} (a_j + x_j)) \} = \min \{ \max_{i\in I\cap K} (a_i + x_i), \max_{j\in J} (a_j + x_j)\} = \min \{ \max_{i\in I} (a_i + x_i), \max_{j\in J} (a_j + x_j)\}$, and by symmetry we deduce that $\max_{j\in J} (a_j+\tilde{x}_j) $ is also equal to the same quantity,
and so $\tilde{x} \in \cH^{IJ}_a$.
\end{proof}
\begin{remark}
  The formula of~\Cref{lem-proj-sh} may be compared 
  with formula for the projection of a point onto a tropical half-space $\{ x\in (\rmax)^n \mid \max_{i\in I} (a_i+x_i) \leq  \max_{j\in J} (a_j+x_j)\}$,
  see~\cite[Th.~5.1]{AGNS10}.
  \end{remark}

\begin{proposition}\label{prop-dist-signed}
Let $x,a \in \proj(\rmax)^n$. The Hilbert's distance of the point $x$ to the signed hyperplane $\cH^{IJ}_a$ is:
\begin{equation}
\dist_H(x,\cH^{IJ}_a) = |\max_{i\in I} (x_i+a_i) -\max_{j\in J} (x_j+a_j) | \enspace ,
\end{equation}
if at least one of these maxima is finite, and $\dist_H(x,\cH^{IJ}_a) =0$
otherwise.
\end{proposition}
\begin{proof}
From~\Cref{th-proj}, we have $\dist_H(x,\cH^{IJ}_a) = d(x,\tilde{x})$ with $\tilde{x}=P_{\cH^{IJ}_a}(x)$.
Let $K=\supp a$ and $O=\supp x$. If $I\cap K \cap O = J\cap K \cap O = \emptyset$, then $K \cap O = \emptyset$, so $x+a\equiv -\infty$, and this means that $x\in \cH^{IJ}_a$ and so $\dist_H(x,\cH^{IJ}_a)=0$.

If  $I\cap K \cap O =  \emptyset$ and  $J\cap K \cap O \neq \emptyset$, then $\max_{i\in I}(x_i+a_i) = -\infty$ and $\max_{j\in J}(x_j+a_j) \neq -\infty$. Let $j\in J\cap K \cap O$, we have $\tilde{x}_j =\min \{ x_j , -a_j + \max_{i\in I} (a_i+x_i)  \} = -\infty$ and we have $x_j \neq -\infty$, then $d(x,\tilde{x})=+\infty=|\max_{i\in I} (x_i+a_i) -\max_{j\in J} (x_j+a_j) |$. 
By symmetry we treat the case when $I\cap K \cap O \neq \emptyset$ and  $J\cap K \cap O = \emptyset$.

Now, we suppose that $I\cap K \cap O \neq \emptyset$ and  $J\cap K \cap O \neq \emptyset$.
Let $i\in I \cap K \cap O$, we have $x_i - \tilde{x}_i = x_i + \max \{-x_i , a_i - \max_{j\in J} (a_j +x_j)\} = 
\max\{0, x_i+ a_i - \max_{j\in J} (a_j +x_j)\}$.
Then, we have $\max_{i\in I}(x_i - \tilde{x}_i) = \max_{i\in I \cap K \cap O}(x_i - \tilde{x}_i) = \max\{0, \max_{i\in I \cap K \cap O}(x_i+ a_i) - \max_{j\in J} (a_j +x_j)\} = \max\{0, \max_{i\in I}(x_i+ a_i) - \max_{j\in J} (a_j +x_j)\}$, and symmetrically, we have $\max_{j\in J}(x_j - \tilde{x}_j) = \max\{0, \max_{j\in J}(x_j+ a_j) - \max_{i\in I} (a_i +x_i)\}$. Therefore, we deduce that $\max_{l\in [n]} (x_l - \tilde{x}_l) = |\max_{i\in I} (x_i+a_i) -\max_{j\in J} (x_j+a_j) |$. 

To finish the proof we need now to show that $\min_{l\in [n]} (x_l - \tilde{x}_l) = 0$. This is a general property of the projection $\tilde{x}=P_M(x)$
of a vector on a closed tropical cone: since $\tilde{x}\leq x$, the minimum
is nonnegative, and if the minimum is positive, adding a small constant $\epsilon$ to every entry of $\tilde{x}$, we get a vector $\tilde{x}^\epsilon$ which
still belongs to $M$ and satisfies $\tilde{x}^\epsilon\leq x$, contradicting
$P_M(x) = \max\{z\in M\mid z\leq x\}$.
\end{proof}

In the sequel, we suppose that the following~\Cref{assum-signed} holds.
\begin{assumption}\label{assum-signed}
We suppose that for each $l\in [n]$, there exists $v\in \cV$, such that $v_l\neq -\infty$.
\end{assumption}

We now introduce the operator $T^{IJ}: (\rmax)^n \mapsto (\rmax)^n$, defined by:

\begin{equation}\label{e-T-IJ}
T^{IJ}_l(x) :=  \left\{
\begin{array}{ll}
\inf_{v\in \cV, v_l\neq -\infty} \{ -v_l +\max_{j\in J} (v_j+x_j) \}, \text{ if } \;  l\in I \enspace ,\\
\inf_{v\in \cV, v_l\neq -\infty} \{ -v_l +\max_{i\in I} (v_i+x_i) \},  \text{ if } \;  l\in J\enspace . 
\end{array}
\right.
\end{equation}

The following result, analogous to~\Cref{lem-hyperplane}, gives a metric
interpretation of the sub-eigenspace of the operator $T^{IJ}$.

\begin{lemma}\label{lem-signed-hyperplane}
Let $\lambda \in [-\infty, 0]$ and $a\in \proj(\rmax)^n$, we have
\[
T^{IJ} (a) \geq \lambda +a  \Leftrightarrow \dist_H(\cV, \cH^{IJ}_a)  \leq - \lambda \enspace.
\]
\end{lemma}

\begin{proof}
The equivalence is trivial if $\lambda=-\infty$,
so, we suppose that $\lambda\in (-\infty,0]$.
We have
\begin{align}
  T^{IJ} (a)  \geq \lambda +a
   & \Leftrightarrow  \left\{
 \begin{array}{ll}
\forall l\in I, \forall v \in \cV, v_l\neq -\infty;  -v_l + \max_{j\in J} (v_j + a_j) \geq \lambda + a_l\\
\forall l\in J, \forall v \in \cV, v_l\neq -\infty;  -v_l + \max_{i\in I} (v_i + a_i) \geq \lambda + a_l
\end{array}
\right.   \\
&  \Leftrightarrow  \left\{
\begin{array}{ll}
\forall v \in \cV, \max_{l\in I}( v_l+ a_l ) \leq  \max_{j\in J} (v_j + a_j) -\lambda \\
\forall v \in \cV, \max_{l\in J}( v_l+ a_l )\leq  \max_{i\in I} (v_i + a_i)  -\lambda 
\end{array}
\right.\label{lastcond}
\end{align}
Let $\cV'$ denote the set of vectors $v\in \cV$ for which at least
one of the latter maxima are finite, and observe that the vectors
of $\cV\setminus \cV'$ trivially belong to $\cH_a^{I,J}$.
Then, using~\Cref{prop-dist-signed}, we see that the last condition in~\eqref{lastcond}
is equivalent to
\[
\forall v \in \cV', d(v,\cH^{IJ}_a) = |\max_{i\in I} (v_i+a_i) -\max_{j\in J} (v_j+a_j) | \leq -\lambda
\]
i.e., 
$\dist_H(\cV, \cH^{IJ}_a)  \leq - \lambda$.
\end{proof}

Let $w\in \R^n$ and $r\geq 0$, we define the {\em vertical interval} of type $I,J$ centered at point $w$ and of radius $r$, 
\[
B_{IJ}(w,r) =
\{ \lambda + w + \mu e^J \mid \mu \in [-r, r],\lambda \in \R\} \enspace,
\]
where
$e^J$ is the vector of $\R^n$ such that $e^J_l =
0$ for $l\in I$ and $e^J_l=1$ for $l\in J$.
Using the identity $-\mu+\mu e^J=-\mu e^I$, we see
\[
B_{IJ}(w,r) =
\{ \lambda + w + \mu e^I \mid \mu \in [-r, r],\lambda \in \R\} \enspace.
\]
\begin{lemma}\label{lem-sball1}
Let $\lambda \in [-\infty, 0]$ and $a\in \R^n$, we have
\[
B_{IJ}(-a,-\lambda) \subset \Sp(\cV) \Rightarrow T^{IJ}(a) \leq \lambda +a \enspace.
\]
\end{lemma}

\begin{proof}
Suppose first that $\lambda$ is finite. 
If $B_{IJ}(-a,-\lambda) \subset \Sp(\cV)$, then 
\[
\forall \mu \in [-\lambda,\lambda], \exists (\alpha_v)_{v\in \cV} \in \R^p, -a+ \mu e^J = \max_{v\in \cV} (\alpha_v + v).
\]
Let $\mu \in [-\lambda,\lambda]$, we have 
\begin{equation}\label{e-sb1}
\forall i \in I, (\forall v \in \cV, -a_i \geq \alpha_v + v_i \text{ and }  \exists v^{(i)}\in \cV, -a_i = \alpha_{ v^{(i)}} +  v^{(i)}_i )\enspace,
\end{equation}
and also
\begin{equation}\label{e-sb2}
\forall j \in J, (\forall v \in \cV, -a_j + \mu \geq \alpha_v + v_j \text{ and }  \exists v^{(j)}\in \cV, -a_j + \mu = \alpha_{ v^{(j)}} +  v^{(j)}_j)\enspace .
\end{equation}
From~\eqref{e-sb2}, we have $\forall v \in \cV, \sup_{j\in J} (v_j+a_j)\leq  - \alpha_v + \mu$, and from~\eqref{e-sb1}, we have $\forall i \in I, v^{(i)}_i \neq -\infty$ because $\alpha_{ v^{(i)}} +  v^{(i)}_i =-a_i \in \R$.
Then, for all $i\in I$, we have $T^{IJ}_i(a)= \inf_{v\in \cV, v_i\neq -\infty} \{ -v_i +\sup_{j\in J} (v_j+a_j) \} \leq \inf_{v\in \cV, v_i\neq -\infty} \{ -v_i - \alpha_v + \mu \} \leq -v^{(i)}_i - \alpha_{v^{(i)}} + \mu = \mu +a_i$. 
This being true for all $\mu \in [-\lambda, \lambda]$, we take here $\mu = \lambda$ and we get that $\forall i \in I, T^{IJ}_i(a) \leq \lambda +a_i$.

Similarly, we have $\forall j\in J, T^{IJ}_j(a)= \inf_{v\in \cV, v_j\neq -\infty} \{ -v_j +\sup_{i\in I} (v_i+a_i) \} \leq \inf_{v\in \cV, v_j\neq -\infty} \{ -v_j - \alpha_v  \} \leq -v^{(j)}_j - \alpha_{v^{(j)}} = -\mu +a_j$. By taking here $\mu = - \lambda$, we get that $\forall j \in J, T^{IJ}_j(a) \leq \lambda +a_j$. 
Therefore, we get that $T^{IJ}(a) \leq \lambda +a$.

The conclusion of the lemma is still
true when $\lambda = -\infty$. This follows from $B_{IJ}(-a,+\infty)
= \cup_{\mu\in (-\infty,0)} B_{IJ}(-a,-\mu)$ and $-\infty + a = \inf_{\mu\in (-\infty,0)} \mu+a$.
\end{proof}
\begin{lemma}\label{lem-sball2}
Let $\lambda\in [-\infty,0]$, %
we have
\[
\exists u \in  \R^n; T^{IJ}(u) \leq \lambda +u  \Rightarrow
\exists w \in \R^n; B_{IJ}(w,-\lambda) \subset \Sp(\cV) \enspace.
\]
\end{lemma}

\begin{proof}
  Suppose first that $\lambda$ is finite.
  For simplicity of notation, we shall
  assume that $u=0$. The general case reduces
  to this one by replacing 
  every vector $v\in \cV$ by the vector $v+u$.
We denote by $V$ the matrix whose columns are the elements of $\cV$.
Since $T^{IJ}(u) \leq \lambda +u$, denoting
by $\sigma$ a map $[n]\to [p]$ such that for all $l\in [n]$, $v=V_{\cdot \sigma(l)}$
achieves the minimum
in~\eqref{e-T-IJ},  we get:
\begin{align}
  \forall l\in I, \;\forall j\in J,\qquad
  V_{j\sigma(l)} &\le V_{l \sigma(l)} + \lambda\enspace,\label{e-v1}\\
  \forall l \in J, \;\forall i\in I,
  \qquad V_{i\sigma(l)}& \le V_{l \sigma(l)} + \lambda
  \enspace .\label{e-v2}
\end{align}
Consider the vectors
\[
w^I:=\vee_{i \in I} -V_{i \sigma(i)} + V_{\cdot \sigma(i)},
\qquad
w^J:=\vee_{j \in J} -V_{j \sigma(j)} + V_{\cdot \sigma(j)},
\]
so that $w^I,w^J\in \Sp(\cV)$. By considering the values $i=l$
or $j=l$ in the suprema above,  we get
\begin{align}
w^I_l\geq 0,\forall l\in I,\qquad 
w^J_l\geq 0 ,\forall l\in J \enspace .
\label{e-lb}
\end{align}
Moreover, using~\eqref{e-v1}, we get
\begin{align}
w^I_j = \vee_{i\in I}
-V_{i \sigma(i)} + V_{j \sigma(i)} 
\leq \lambda ,
\qquad \text{ for all }j\in J\enspace,
\label{e-v7}
\end{align}
and similarly, using~\eqref{e-v2},
\begin{align}
w^J_i \leq \lambda, \text{ for all }i\in I \enspace .
\label{e-v8}
\end{align}
Define the vector $w$ by
\[
w_l = \left\{
\begin{array}{ll}
w^I_l, \text{ if } l\in I \enspace,\\
w^J_l, \text{ if } l\in J\enspace .
\end{array}\right.
\]
Using~\eqref{e-lb}--\eqref{e-v8}, we
deduce that for all $\mu\in [\lambda,-\lambda]$,
\[
w+\mu e^I = (w^I+\mu )\vee w^J \in \Sp(\cV) \enspace ,
\]
and so $B_{IJ}(w,-\lambda)\subset \Sp(\cV)$.

We finally show that the conclusion of the lemma is still
true when $\lambda = -\infty$. This follows from the fact that the above center $w$ depends only on the vectors of $\cV$ and does not depend on $\lambda$, and also from the facts that $B(w,+\infty) = \cup_{\mu\in (-\infty,0)} B(w,-\mu)$ and $-\infty + w = \inf_{\mu\in (-\infty,0)} \mu+w$.
\end{proof}

The next result is immediate from~\Cref{lem-sball1,lem-sball2}. It is analogous to~\Cref{lem-ball}. It shows that the existence of a super-eigenvector of $T^{I,J}$ is equivalent to the existence of a {\em vertical interval} included in the module $\Sp(\cV)$. 
\begin{proposition}\label{coro-signed-ball}
Let $\lambda\in [-\infty,0]$, and $a\in \R^n$, we have
\[
\exists u \in  \R^n; T^{IJ}(u) \leq \lambda +u  \Leftrightarrow \exists w \in \R^n; B_{IJ}(w,-\lambda) \subset \Sp(\cV) \enspace.
\pushQED{\qed}
\qedhere
\popQED
\]
\end{proposition}

We now derive a strong duality theorem for signed tropical regression.

\begin{theorem}\label{th-strong-duality-signed}
We have
\begin{equation}\label{e-strong-duality-IJ}
\min_{a\in \proj(\rmax)^n} \dist_H(\cV, \cH^{IJ}_a) = - \rho(T^{IJ}) = \sup \{r \geq 0\mid \exists w \in \R^n, B_{IJ}(w,r) \subset  \Sp(\cV) \} \enspace.
\end{equation}
The minimum is achieved by any vector $b \in \proj(\rmax)^n$
such that $T^{IJ}(b)\geq \rho(T^{IJ}) + b$. Moreover, if $\rho(T^{IJ})$ is finite, the supremum is achieved by a ball $B(c,\rho(T^{IJ}))$ 
where $c \in \R^n$ can be deduced from any vector $u$ such that $T^{IJ}(u)\leq \rho(T^{IJ}) +u$. 
\end{theorem}

\begin{proof}
From~\Cref{coro-signed-ball}, \Cref{lem-signed-hyperplane} and the Collatz-Wielandt property (\Cref{th-duality}), we deduce the strong duality property~\eqref{e-strong-duality-IJ}.
Moreover, the existence of a vector $b\in \proj(\rmax)^n$ such
$T^{IJ}(b)\geq \rho(T^{IJ}) + b$ follows from \Cref{th-duality}.
Then, by~\Cref{lem-signed-hyperplane}, we have $\dist_H (\cV,\cH_b) \leq -\rho(T^{IJ})$,
which implies that $b$ such that $T^{IJ}(b)\geq \rho(T^{IJ}) + b$ achieves the minimum in~\eqref{e-strong-duality-IJ}.

Finally, if $\rho(T^{IJ})$ is finite, since the infimum is attained in the expression of the Collatz-Wielandt number of $T^{IJ}$ (see \Cref{th-duality}), there exists a finite vector $u\in \R^n$ such that $T^{IJ}(u)\leq \rho(T^{IJ})+u$. By the proof of~\Cref{lem-sball2}, we can then construct a vector $c$ such that $B_{IJ}(c,-\rho(T^{IJ})) \subset \Sp(\cV)$.
\end{proof}

\begin{remark}\label{rk-ordinaryregression}
When the set $I=\{i\}$ is of cardinality one, the
regression problem for the signed hyperplane~\eqref{signed-hyperplane}
has the following special form:
\begin{align}
\Minimize_{a\in \R^n} \max_{v\in \cV} |v_i - (\max_{j\neq i} a_j -a_i + v_j)|
\enspace . \label{e-onesided}
\end{align}
This can be solved in a direct way~\cite{maragossurvey},
avoiding the recourse to mean payoff games.%
Indeed, \eqref{e-onesided} reduces to the following ``one-sided'' tropical linear
regression problem. Given sample points $(x^{(k)},y^{(k)})$ in $\R^n\times \R^m$, for $k\in [p]$, compute
  \begin{align}
    \Minimize_A \max_{k\in [p]} \|y^{(k)}-Ax^{(k)}\|_\infty\enspace ,\label{e-regression-onesided} 
  \end{align}
where the minimum is taken over tropical matrices $A$ of size $m\times n$,
and the product $Ax^{(k)}$ is understood tropically. 
Up to a straightforward duality, this problem was solved
in~\cite[Theorem 3.5.2]{butkovic}, the result being
attributed there to Cuninghame-Green~\cite{CG}.
Alternatively, this solution may be recovered by combining~\cite[Coro.~1]{chepoi} with the explicit formula of the tropical projection~\cite[Th.~5]{cgq02}.
More precisely, define the matrix $\bar{A}\in \R^{m\times n}$ by
$\bar{A}_{ij}:= \min_{k\in [p]} y^{(k)}_i - x^{(k)}_{j} $,
so that $\bar A$ is the maximal matrix such that $Ax^{(k)}\leq y^{(k)}$
for all $k\in [p]$. Let $\delta := \max_{k\in p} \|y^{(k)}-\bar{A}x^{(k)}\|_\infty$,
and $A^{\textrm{opt}}_{ij}= \bar{A}_{ij} + \delta/2$. Then, $A^{\textrm{opt}}$
is the greatest optimal solution. It can be computed in $O(mnp)$ arithmetic
operations. By specializing this formula, one can solve~\eqref{e-onesided}
in $O(np)$ arithmetic operations. We refer the reader to~\cite{maragossurvey}
for more information, and for the solution of further problems
of this category. 
\end{remark}
  \begin{remark}
  In contrast, when $I,J$ are part of the input,
  the signed linear tropical regression problem
  is polytime Turing equivalent to mean payoff games.
  This can be seen as follows.
  The reduction to mean payoff games is a consequence
  of~\Cref{th-strong-duality-signed}.
  Conversely, observe that finding a signed
  tropical hyperplane $\cH^{I,J}_a$ containing a set
  $\cV=\{v^{(1)},\dots,v^{(p)}\}$ in $\R^n$ is equivalent to solving a
  tropical linear system of the form
  $Bx = Cy$, where $x\in (\rmax)^I$,
  $y\in (\rmax)^J$, $B\in (\rmax)^{p\times I}$,
  $C\in (\rmax)^{p\times J}$, $B_{ki}=v^{(k)}_i$ for $i\in I$
  and $C_{kj}=v^{(k)}_j$ for $j\in I$. Indeed, the vector $a$
  defining $\cH^{I,J}_a$ is given by $a_i= x_i$ for $i\in I$
  and $a_j=y_j$ for $j\in J$. We know from~\cite{AGGut10}
  that deciding whether a mean payoff game has an initial
  winning position is equivalent to the existence
  of a non-identically $-\infty$ solution $z\in (\rmax)^s$
  of a system of tropical linear inequalities $Fz\leq Gz$,
  where $F,G\in (\zmax)^{r\times s}$ are given.
  Such a system $Fz\leq Gz$ can be rewritten as $Bx=Cy$
  by introducing lift variables $u,v\in (\rmax)^r$,
  so that $v=Fz$ and $u$ is a slack variable.
  Setting $y:=(u,v)$, identified to a column vector,
  $B:=\left(\begin{smallmatrix}\zero& \I \\ \I & \I \end{smallmatrix}
  \right)$ and
  $C:=\left(\begin{smallmatrix}F \\ G \end{smallmatrix}
  \right)$, where $\zero$ is a zero tropical matrix,
  and $\I$ the identity matrix, we see
  that $Fz\leq Gz$ has a non-identically $-\infty$
  solution iff $By=Cz$ has a non-identically $-\infty$ solution.
  It follows that mean payoff
  games reduce to checking whether
  there is a solution of a signed tropical linear regression problem
  with zero error.
  \end{remark}

\subsection{Tropical linear regression with type information}
\label{subsec-typed}
The following variant will be relevant to the application to economy considered below, to measure the ``distance to equilibria'' of a market. We suppose the
set of points $\cV$ is the disjoint union $\cV = \cup_{i\in [n]} \cV_i$, where each $\cV_i$ is non-empty. We shall say that the points of $\cV_i$ are of {\em type
$i\in[n]$}. Note that the set of types is the same as the set
of indices of vectors. For each type $i\in [n]$, we consider the
signed hyperplane:
\begin{equation*}
\cH_a^i :=\cH_a^{\{i\}\{i\}^c} = \{ x\in (\rmax)^n \mid a_i+x_i = \max_{j\neq i} (a_j+x_j)\}\enspace. 
\end{equation*}

The {\em typed tropical linear regression} problem associated
to the partition $\cV_1,\dots,\cV_n$ of $\cV$, is defined
as:
\begin{equation}\label{eq-reg-types}
\Minimize_{a\in \proj(\rmax)^n} \max_{i \in [n]} \dist_H(\cV_i , \cH_a^i) \enspace .
\end{equation}

The value of this problem is small if and only if
for each $i\in [n]$, the points of $\cV_i$ are close
to the signed tropical hyperplane $\cH_a^i$.

From~\Cref{prop-dist-signed}, we know that $\dist_H(v,\cH_a^i) = | v_i+a_i - \max_{j\neq i} (v_j+a_j) |$.

We suppose in the sequel that~\Cref{assum-signed} holds. For each type $i\in [n]$, we consider the Shapley operator $\Ttypei: (\rmax)^n \mapsto (\rmax)^n$, given by~\eqref{e-T-IJ} where the type considered is $(I,J)=(\{i\},\{i\}^c)$ and to the set of points is $\cV_i$: 

\begin{equation}\label{e-T-typesi}
\Ttypei_l(x) :=  \left\{
\begin{array}{ll}
\inf_{v\in \cV_i, v_i\neq -\infty} \{ -v_i +\max_{j\neq i} (v_j+x_j) \}, \text{ if } \;  l= i \enspace ,\\
\inf_{v\in \cV_i, v_l\neq -\infty} \{ -v_l + v_i \}+x_i ,  \text{ if } \;  l\neq i \enspace .
\end{array}
\right.
\end{equation}

We consider now the Shapley operator $\Ttype: (\rmax)^n \mapsto (\rmax)^n$ given by the infimum of the operators $\Ttypei, i\in [n]$. It is given by:

\begin{equation}\label{e-T-types}
\Ttype_l(x) := \min_{i\in [n]}  \Ttypei_l(x) \enspace.
\end{equation}

The following lemma, analogous to~\Cref{lem-hyperplane}, gives a metric
interpretation of the sub-eigenspace of the operator $T$.

\begin{lemma}\label{lem-hyperplane-types}
Let $\lambda \in [-\infty, 0]$ and $a\in \proj(\rmax)^n$, we have
\[
\Ttype (a) \geq \lambda +a  \Leftrightarrow \max_{i\in [n]}\dist_H(\cV_i, \cH^i_a)  \leq - \lambda \enspace.
\]
\end{lemma}

\begin{proof}
Let $\lambda \in [-\infty, 0]$ and $a\in \proj(\rmax)^n$. From~\eqref{e-T-types} and~\Cref{lem-signed-hyperplane}, we deduce the equivalence:\phantom{\qedhere}
\begin{align*}
\Ttype (a)  \geq \lambda +a
& \Leftrightarrow \forall i\in [n], \Ttypei (x) \geq \lambda + a \\
& \Leftrightarrow \forall i\in [n], \dist_H(\cV_i, \cH^i_a)  \leq - \lambda \enspace, \\
& \Leftrightarrow \max_{i\in [n]}\dist_H(\cV_i, \cH^i_a)  \leq - \lambda \enspace.
\pushQED{\qed}
\qedhere
\popQED
\end{align*}
\end{proof}

From~\Cref{lem-hyperplane-types,th-duality}, we deduce the following result, showing that the tropical linear regression problem with types,
associated to the sets
$\cV_1,\dots,\cV_n$, also reduces to a mean payoff game.

\begin{theorem}\label{th-typed}
We have,
\[
\min_{a\in \proj(\rmax)^n} \max_{i\in [n]}\dist_H(\cV_i, \cH^i_a) = -\rho(\Ttype) \enspace.
\]
Moreover, the minimum is achieved by any vector $a\in \proj(\rmax)^n$
such that $\Ttype (a) \geq \rho(\Ttype)+a$.
\end{theorem}
\begin{remark}\label{rk-SVM}
  Typed tropical linear regression should be compared
  with the tropical SVM problem introduced in~\cite{gartnerJaggi}.
  In the tropical SVM setting, we have a partition of the set
  of points in $n$ {\em color classes}, $\cV_{c_1},\dots,\cV_{c_n}$,
  and we are looking
  for a tropical hyperplane $\cH_a$, and for a permutation
  $\sigma$ of $\{1,\dots,n\}$ such that for
  all $i\in [n]$, all the points of color $c_i$ are in the same
  {\em sector} $S_{\sigma(i)}(a)$. In other words, we want the tropical
  hyperplane to separate the $n$ color classes. This is not
  possible in general, so one needs to consider metric versions,
    modelling the minimization of
    classification errors~\cite{tropicalSVM}.
  A possible metric formulation, in the spirit of the present approach,
  would be to consider
  \begin{align}
  \min_{\sigma \in \mathfrak{S}_n} \min_{a\in \R^n} \max_{i\in [n]}
  \dist_H(\cV_i,S_{\sigma(i)}(a))
  \qquad \textrm{(Metric Tropical SVM)}
  \end{align}
  where $\mathfrak{S}_n$ denotes the symmetric group on $n$ letters.
  By comparison with~\eqref{eq-reg-types}, we see that we have in addition
  a minimization over the symmetric group, but the subproblem with a fixed
  permutation $\sigma$ arising in the SVM problem
  is simpler than the analogous problem of typed tropical linear regression,
  since the
  sector $S_{\sigma(i)}$ is convex, whereas the set $\cH^i_a$ arising
  in~\eqref{eq-reg-types} is not a convex one. In the application described
  below, it is the set $\cH^i_a$ that is relevant to measure
  the ``distance to equilibrium''.
\end{remark}

In~\Cref{figTypedRegression1}, we consider the following matrix $V\in \R^{3 \times 11}$:
\begin{equation}\label{e-example-typed}
V=\left(\begin{array}{ccccccccccc}
    1 &  1 &  2 &  0  & 0 &  0 & -3 & -1 &  0 &  0 & -2 \\
    0 & -2 &  0 &  1 &  1 &  2 &  1 &  0 &  0 & -3 &  0 \\
    0 &  0 & -2 & -2 & -1 & -2 &  0 &  2 &  3 &  1 &  1
\end{array}
\right) ,
\end{equation}

and the types are given by the subsets of $\cV=[11]$ as follows $\cV_1=\{1,2,3,4\}$, $\cV_2=\{5,6,7,8\}$ and $\cV_3=\{9,10,11\}$.

The operators $\Ttypei: (\rmax)^n \mapsto (\rmax)^n$ given by~\eqref{e-T-typesi} and associated to the above matrix $V$ and partition $(\cV_i)_{i\in[3]}$ are given by:

\[
T^{\mathrm{ty},1} \begin{pmatrix}
x_1 \\ x_2 \\ x_3
\end{pmatrix}
 = \begin{pmatrix}
  \min[-1+\max(-2+x_2,x_3), -2+\max(x_2,-2 + x_3), \max(1+x_2,-2+x_3) ]  \\ %
  -1 +x_1  \\
  1 +x_1  
\end{pmatrix} \enspace ,
\]

\[
T^{\mathrm{ty},2} \begin{pmatrix}
x_1 \\ x_2 \\ x_3
\end{pmatrix}
 = \begin{pmatrix}
  1 +x_2  \\
  \min[-2+\max(x_1,-2+x_3), -1+\max(-3+x_1, x_3), \max(-1+x_1, 2+x_3) ]  \\ %
 -2 + x_2   
\end{pmatrix} \enspace ,
\]

\[
T^{\mathrm{ty},3} \begin{pmatrix}
x_1 \\ x_2 \\ x_3
\end{pmatrix}
 = \begin{pmatrix}
  1+x_3  \\
  1+x_3  \\
  \min[-3+\max(x_1,x_2), -1+\max(x_1,-3+x_2), -1+\max(-2+x_1, x_2)]    
\end{pmatrix} \enspace .
\]

Then the operator $\Ttype: (\rmax)^n \mapsto (\rmax)^n$ given by~\eqref{e-T-types} is in this example:

\[
\Ttype \begin{pmatrix}
x_1 \\ x_2 \\ x_3
\end{pmatrix}
 = \begin{pmatrix}
  \min[-1+\max(-2+x_2,x_3), -2+\max(x_2,-2 + x_3), 1 +x_2, 1+x_3, ]  \\ %
  \min[-2+\max(x_1,-2+x_3), -1+\max(-3+x_1, x_3), -1 +x_1, 1+x_3]  \\  %
  \min[-3+\max(x_1,x_2), -1+\max(x_1,-3+x_2), 1 +x_1, -2 + x_2]     %
\end{pmatrix} \enspace .
\]

We verify easily that $\lambda=-2$ and $a=(0,0,-1)^\top$ satisfy $\Ttype(a)=\lambda +a$, so that by~\Cref{th-typed} the apex $a$ is optimal for the typed tropical linear regression problem~\eqref{eq-reg-types}.

We notice that in this case, the tropical hyperplane $\cH_a$ has at least one witness point in each sector, which means, by~\Cref{th-witness}, that $\cH_a$ is also an optimal hyperplane in the sense of the usual tropical linear regression studied in~\Cref{sec-strong}.

Now, if we consider the same matrix $V$ in~\eqref{e-example-typed}, but we exchange the types of the points $V_{\cdot 8}$ and $V_{\cdot 10}$, i.e. we consider the partition $\widetilde{\cV}_1=\{1,2,3,4\}$, $\widetilde{\cV}_2=\{5,6,7,10\}$ and $\widetilde{\cV}_3=\{8,9,11\}$, then the new typed Shapley operator $\widetilde{\Ttype}$ is given by:
\[
\widetilde{\Ttype} \begin{pmatrix}
x_1 \\ x_2 \\ x_3
\end{pmatrix}
 = \begin{pmatrix}
  \min[-3 +x_2, 3+x_3 ]  \\ 
  \min[-2+\max(x_1,-2+x_3), -1+\max(-3+x_1, x_3), -1 +x_1, 1+x_3]  \\  
  \min[1 +x_1, -4 + x_2]   
\end{pmatrix} \enspace .
\]

We verify easily that $\mu=-5/2$ and $b=(0,1/2,-1)^\top$ satisfy $\widetilde{\Ttype}(b)=\mu +b$. This example is presented in~\Cref{figTypedRegression2}. Here, we notice that the hyperplane $\cH_b$ that is optimal in the typed tropical linear regression sense (\Cref{subsec-typed}) does not have witness points in each sector, which means that it is not optimal in the usual tropical linear regression framework (\Cref{sec-strong}).

\begin{figure}[htbp]
\centering
\subfigure{


\def\coord#1#2#3{{-sqrt(3)/2*(#1-#2)} ,{ -(1/2)*#1 - (1/2)*#2 + #3}}
\begin{tikzpicture}[scale=0.75]

\coordinate (v1) at (\coord{1}{0}{0});
\coordinate (v2) at (\coord{1}{-2}{0});
\coordinate (v3) at (\coord{2}{0}{-2});
\coordinate (v4) at (\coord{0}{1}{-2});

\coordinate (v5) at (\coord{0}{1}{-1});
\coordinate (v6) at (\coord{0}{2}{-2});
\coordinate (v7) at (\coord{-3}{1}{0});
\coordinate (v8) at (\coord{-1}{0}{2});

\coordinate (v9) at (\coord{0}{0}{3});
\coordinate (v10) at (\coord{0}{-3}{1});
\coordinate (v11) at (\coord{-2}{0}{1});

\coordinate (Ha1) at (\coord{-5}{0}{1});
\coordinate (Ha2) at (\coord{0}{-5}{1});
\coordinate (Ha3) at (\coord{0}{0}{-4});

\coordinate (f) at (\coord{0}{0}{-3});

\coordinate (b1) at (\coord{0}{0}{3});
\coordinate (b2) at (\coord{0}{2}{3});
\coordinate (b3) at (\coord{0}{2}{1});
\coordinate (b4) at (\coord{0}{0}{-1});
\coordinate (b5) at (\coord{2}{0}{1});
\coordinate (b6) at (\coord{2}{0}{3});

\fill[red!30,opacity=0.7] (Ha1) -- (\coord{0}{0}{1}) -- (Ha2)-- (\coord{0}{0}{4.5})--cycle;
\fill[blue!30,opacity=0.7] (Ha1) -- (\coord{0}{0}{1}) -- (Ha3)-- (\coord{0}{5}{1})--cycle;
\fill[green!30,opacity=0.7] (Ha2) -- (\coord{0}{0}{1}) -- (Ha3)-- (\coord{5}{0}{1})--cycle;



\draw[very thick,red] (\coord{0}{0}{1}) -- (Ha1);
\draw[very thick,red] (\coord{0}{0}{1}) -- (Ha2);
\draw[very thick,red] (\coord{0}{0}{1}) -- (Ha3);
\filldraw[black] (\coord{0}{0}{1}) circle (1.5pt) node[below,right] {$-a$};
\draw node[above] at (\coord{-1/2}{0}{1}) {$\cH_a$};
\filldraw[black] (\coord{0}{0}{0}) circle (1.5pt) node[below,right] {0};

\filldraw[black] (v1) circle (3pt) node[above,black] {$V_{\cdot 1}$};
\filldraw[black] (v2) circle (3pt) node[above,black] {$V_{\cdot 2}$};
\filldraw[black] (v3) circle (3pt) node[below,left,black] {$V_{\cdot 3}$};
\filldraw[black] (v4) circle (3pt) node[right,black] {$V_{\cdot 4}$};

\filldraw[black] ([xshift=-2.5pt,yshift=-2.5pt]v5) rectangle ++(5pt,5pt) node[right,black] {$V_{\cdot 5}$};
\filldraw[black] ([xshift=-2.5pt,yshift=-2.5pt]v6) rectangle ++(5pt,5pt) node[below right,black] {$V_{\cdot 6}$};
\filldraw[black] ([xshift=-2.5pt,yshift=-2.5pt]v7) rectangle ++(5pt,5pt) node[below right,black] {$V_{\cdot 7}$};
\filldraw[red] ([xshift=-2.5pt,yshift=-2.5pt]v8) rectangle ++(5pt,5pt) node[above,black] {\textcolor{red}{$V_{\cdot 8}$}};

\node[mark size=2.5pt,color=black] at (v9) {\pgfuseplotmark{triangle*}};
\draw node[above] at (v9) {$V_{\cdot 9}$};
\node[mark size=2.5pt,color=red] at (v10) {\pgfuseplotmark{triangle*}};
\draw node[above] at (v10) {\textcolor{red}{$V_{\cdot 10}$}};
\node[mark size=2.5pt,color=black] at (v11) {\pgfuseplotmark{triangle*}};
\draw node[right] at (v11) {$V_{\cdot 11}$};

\draw[dashed,->] (\coord{0}{0}{0}) -- (\coord{5}{0}{0}) node[above] {$x_1$};
\draw[dashed,->] (\coord{0}{0}{0}) -- (\coord{0}{5}{0}) node[above] {$x_2$};
\draw[dashed,->] (\coord{0}{0}{0}) -- (\coord{0}{0}{4.5}) node[above, right] {$x_3$};

\draw[dashed] (b1) -- (\coord{-3}{0}{3});
\draw[dashed] (b1) -- (\coord{0}{-3}{3});
\draw[dashed] (b3) -- (\coord{0}{2}{-3});
\draw[dashed] (b3) -- (\coord{-3.5}{2}{1});
\draw[dashed] (b5) -- (\coord{2}{0}{-3});
\draw[dashed] (b5) -- (\coord{2}{-3.5}{1});

\end{tikzpicture}
\label{figTypedRegression1}
}
\subfigure[]{


\def\coord#1#2#3{{-sqrt(3)/2*(#1-#2)} ,{ -(1/2)*#1 - (1/2)*#2 + #3}}
\begin{tikzpicture}[scale=0.75]

\coordinate (v1) at (\coord{1}{0}{0});
\coordinate (v2) at (\coord{1}{-2}{0});
\coordinate (v3) at (\coord{2}{0}{-2});
\coordinate (v4) at (\coord{0}{1}{-2});

\coordinate (v5) at (\coord{0}{1}{-1});
\coordinate (v6) at (\coord{0}{2}{-2});
\coordinate (v7) at (\coord{-3}{1}{0});
\coordinate (v8) at (\coord{-1}{0}{2});

\coordinate (v9) at (\coord{0}{0}{3});
\coordinate (v10) at (\coord{0}{-3}{1});
\coordinate (v11) at (\coord{-2}{0}{1});

\coordinate (Ha1) at (\coord{-5}{-0.5}{1});
\coordinate (Ha2) at (\coord{0}{-5.5}{1});
\coordinate (Ha3) at (\coord{0}{-0.5}{-4});

\coordinate (f) at (\coord{0}{0}{-3});

\coordinate (b1) at (\coord{0}{-0.5}{3.5}); 
\coordinate (b3) at (\coord{0}{2}{1}); 
\coordinate (b5) at (\coord{2.5}{-0.5}{1}); 

\fill[red!30,opacity=0.7] (Ha1) -- (\coord{0}{-0.5}{1}) -- (Ha2)-- (\coord{0}{-0.5}{4.5})--cycle; 
\fill[blue!30,opacity=0.7] (Ha1) -- (\coord{0}{-0.5}{1}) -- (Ha3)-- (\coord{0}{4.5}{1})--cycle;
\fill[green!30,opacity=0.7] (Ha2) -- (\coord{0}{-0.5}{1}) -- (Ha3)-- (\coord{5}{-0.5}{1})--cycle;



\draw[very thick,red] (\coord{0}{-0.5}{1}) -- (Ha1);
\draw[very thick,red] (\coord{0}{-0.5}{1}) -- (Ha2);
\draw[very thick,red] (\coord{0}{-0.5}{1}) -- (Ha3);
\filldraw[black] (\coord{0}{-0.5}{1}) circle (1.5pt) node[below,right] {$-b$};
\draw node[above] at (\coord{-1/2}{-0.5}{1}) {$\cH_b$};
\filldraw[black] (\coord{0}{0}{0}) circle (1.5pt) node[below,right] {0};

\filldraw[black] (v1) circle (3pt) node[above,black] {$V_{\cdot 1}$};
\filldraw[black] (v2) circle (3pt) node[above,black] {$V_{\cdot 2}$};
\filldraw[black] (v3) circle (3pt) node[below,left,black] {$V_{\cdot 3}$};
\filldraw[black] (v4) circle (3pt) node[right,black] {$V_{\cdot 4}$};

\filldraw[black] ([xshift=-2.5pt,yshift=-2.5pt]v5) rectangle ++(5pt,5pt) node[right,black] {$V_{\cdot 5}$};
\filldraw[black] ([xshift=-2.5pt,yshift=-2.5pt]v6) rectangle ++(5pt,5pt) node[below right,black] {$V_{\cdot 6}$};
\filldraw[black] ([xshift=-2.5pt,yshift=-2.5pt]v7) rectangle ++(5pt,5pt) node[below right,black] {$V_{\cdot 7}$};
\node[mark size=2.5pt,color=red] at (v8) {\pgfuseplotmark{triangle*}};
\draw node[above] at (v8) {\textcolor{red}{$V_{\cdot 8}$}};

\node[mark size=2.5pt,color=black] at (v9) {\pgfuseplotmark{triangle*}};
\draw node[above] at (v9) {$V_{\cdot 9}$};
\filldraw[red] ([xshift=-2.5pt,yshift=-2.5pt]v10) rectangle ++(5pt,5pt) node[above,red] {\textcolor{red}{$V_{\cdot 10}$}};
\node[mark size=2.5pt,color=black] at (v11) {\pgfuseplotmark{triangle*}};
\draw node[right] at (v11) {$V_{\cdot 11}$};

\draw[dashed,->] (\coord{0}{0}{0}) -- (\coord{5}{0}{0}) node[above] {$x_1$};
\draw[dashed,->] (\coord{0}{0}{0}) -- (\coord{0}{5}{0}) node[above] {$x_2$};
\draw[dashed,->] (\coord{0}{0}{0}) -- (\coord{0}{0}{5}) node[above, right] {$x_3$};

\draw[dashed] (b1) -- (\coord{-2.5}{-0.5}{3.5});
\draw[dashed] (b1) -- (\coord{0}{-3}{3.5});
\draw[dashed] (b3) -- (\coord{0}{2}{-2.5});
\draw[dashed] (b3) -- (\coord{-3}{2}{1});
\draw[dashed] (b5) -- (\coord{2.5}{-0.5}{-2.5});
\draw[dashed] (b5) -- (\coord{2.5}{-3.5}{1});

\end{tikzpicture}
\label{figTypedRegression2}
}
\caption{\Cref{figTypedRegression1}: A set of typed points $\cV$ with three types in $\proj(\rmax)^3$ with an optimal tropical hyperplane $\cH_a$ in the sense of the typed tropical regression, where $a=(0,0,-1)^\top$ satisfies $\Ttype(a)=-2+a$.
\Cref{figTypedRegression2}: The same set of typed points $\cV$ as~\Cref{figTypedRegression1} but with the types of the two points $V_{\cdot 8}$ and $V_{\cdot 10}$ being exchanged, and an optimal tropical hyperplane $\cH_b$ in the sense of the typed tropical regression, where $b=(0,1/2,-1)^\top$ satisfies $\widetilde{\Ttype}(b)=-5/2+b$.}
\end{figure}

\section{Algorithmic aspects}
\label{sec-algo}
In this section, we explain how the tropical linear
regression problem can be effectively solved
by using mean-payoff games algorithms. Throughout the
section, we assume that the set of points $\cV$ is given
by as the set of columns the matrix $V$.
By \Cref{cor-center}, in theory, any
algorithm solving mean payoff games in the weakest
sense (deciding the inequality $\chi_i(T)\geq 0$)
can be used. However, some game algorithms lead to
more direct approaches, we next discuss some of these.

Considering the strong duality
result, \Cref{th-strong-duality},
and the result on the existence
of witness points~\Cref{th-witness}, the key algorithmic
issues are:
\begin{itemize}
  \item[(i)] to compute the upper mean payoff,
$\rho(T)$ (which is the opposite of the value
    of the tropical linear regression problem);
    \item[(ii)] to decide whether
there is a finite eigenvector $u\in\R^n$
such that $T(u) = \rho(T)+u$, and to compute
such an eigenvector (when this is so, $-u$ is the center
an an optimal ball included in $\Sp(\cV)$ and the apex
of an optimal regression hyperplane, see~\Cref{rk-rk});
\item[(iii)] to find a sub-eigenvector
$b\in (\rmax)^n\setminus\{\bot\}$, satisfying $T(b) \geq \rho(T)+b$
  (then, $\cH_b$ is an optimal regression hyperplane);
  \item[(iv)] to find a super-eigenvector $c\in \R^n$
    satisfying $T(c)\leq \rho(T)+ c$ (then, $-c$ is the center
    of an optimal ball included in $\Sp(\cV)$.
\end{itemize}
For simplicity of the discussion, we assume that $T$ sends $\R^n\to \R^n$.
The case in which $T$ sends $\R^n$ to $(\rmax)^n$
reduces to this one by considering the action of $T$ on
the parts of $(\rmax)^n$ and looking for invariant
parts. 

Then, problems (i)--(iv) are solved, simultaneously, as
soon as we know an invariant half-line of $T$.
Indeed, we observed after stating~\Cref{th-duality} that if $(u,\eta)$
is an invariant half-line, then $\chi(T) =\eta$.
In this way,
$\rho(T) = \max_{i\in[n]}\chi_i(T)$ is determined,
and this solves issue (i). 
Moreover, by~\Cref{prop-finiteeig}
$T$ admits a finite eigenvector if and only if $\eta$ is a constant
vector, i.e., $\eta=(\lambda,\dots,\lambda)$ for some $\lambda\in\R$,
and $u$ is an eigenvector. This solves issue~(ii).
We observed in the proof of~\Cref{th-duality} that $u$ satisfies
$T(u) \leq \rho(T)+u$, and so, this solves issue (iii).
Finally, setting $I:=\{i\in[n]\mid \chi_i(T)=\rho(T)\}$,
and defining the vector $\bar{u}$ such that $\bar u_i=u_i$
for $i\in I$ and $\bar{u}_i=-\infty$ otherwise, it can be checked
that $T(\bar u)\geq \bar u + \rho(T)$, which solves issue (iv).

More generally, the reduction in the second part
of the proof of~\Cref{cor-center}
shows that  algorithm which returns an optimal policy $\sigma$ of Player Min,
i.e., a policy such that $\chi(T)=\chi(T^\sigma)$, 
can be used to produce a finite vector $c\in \R^n$ such that
$T(c)\leq \overline{\chi}(T)+c$, by reduction
to a tropical eigenvalue problem. Moreover,
any algorithm which returns an optimal policy $\tau$ of Player Max,
i.e., a policy such that $\chi(T)=\chi(^\tau T)$, 
can be used to produce a  vector
$b\in (\rmax)^n\setminus\{\bot\}$, satisfying $T(b) \geq \rho(T)+b$,
see the second part of the reduction in~\Cref{cor-dualcenter}.

We refer the reader to~\cite{chaloupka}
for a comparative discussion of mean payoff
game algorithms. 
The main known algorithms include
the pumping algorithm of~\cite{gurvich}, 
value iteration~\cite{zwick}, 
and different algorithms based on the idea of policy
iteration~\cite{bjorklund,schewe08,DG-06}.
In particular, the algorithm of~\cite{DG-06} returns
an invariant half-line.
The policy iterations algorithms~\cite{bjorklund,DG-06}
were reported in~\cite{chaloupka} to have the best
experimental behavior, although
policy iteration is are generally exponential~\cite{friedman}.

For the present application to tropical linear regression,
we often know in advance that the operator $T$ has a finite
eigenvector; this occurs in particular if the
entries of the matrix $V$ are finite, and more generally,
under the dominion condition of~\Cref{thm-dominions}.
Then, one can use another algorithm,
projective Krasnoselkii-Mann value iteration~\cite{stott2020},
which is straightforward to implement and still effective.
Starting from a vector $v^0=(0,\cdots, 0)^\top$, this algorithm
computes the following sequence:
\begin{eqnarray}
& \tilde{v}^{k+1}=T(v^{k}) - (\max_{i\in [n]}T(v^{k})_i) e, \\
& v^{k+1}= (1-\gamma) v^{k} +\gamma  \tilde{v}^{k+1} .
\end{eqnarray}
where $e=(1,\cdots,1)^\top \in \R^n$, and $\gamma\in (0,1)$ is fixed,
$1-\gamma$ being interpreted as a {\em damping parameter}. In the original
Krasnoselskii-Mann algorithm, one writes simply $v^{k+1}=(1-\gamma)v^k + \gamma T(v^k)$. It follows from~\cite[Coro.~13]{stott2020}, based on a general result of
Baillon and Bruck~\cite{baillonbruck} on the convergence
of the original Krasnoselskii-Mann algorithm in normed
spaces, see also~\cite{cominetti}, that $v^k$ does converge
to an eigenvector of $T$ as soon as such a (finite) eigenvector $u$ exists.
Moreover, $\|T(v^k)-v^k\|_H \leq 2\|u\|_H /\sqrt{\pi \gamma(1-\gamma)k}$.
In practice, we fix a desired precision $\epsilon>0$,
and stop the computation of the sequence $v^k$
when $\|T(v^k)-v^k\|_H\leq \epsilon$.

We now analyze the complexity of the projective Krasnoselskii-Mann
algorithm in our special setting. The following observation,
shows that, notwithstanding the quadratic size
of $\maxset$ in the game associated
with $T$ (see the discussion after~\eqref{e-def-T}),
the operator $T$ can be evaluated in linear time. 
\begin{proposition}\label{prop-eval}
  The operator $T$ can be evaluated in $O(|E|)$ arithmetic operations.
\end{proposition}
\begin{proof}
  We write $T_i(x)= \min_{k\in [p],(i,k)\in E } (-V_{ik}+y_{ik})$ where $y_{ik}=\max_{j\in [n],j\neq i,(j,k)\in E} (V_{jk}+x_k)$.
  First, for each column $k$ of the matrix $V$, we compute the column
  maximum $M_k:=\max_{j\in [n],(j,k)\in E}(V_{jk}+x_k)$ together with an arbitrary index $j_k$
  that achieves this maximum, and also the second column maximum,
  $m_k:=\max_{j\in [n],j\neq j_k, (j,k)\in E}(V_{jk}+x_k)$. This preprocessing
  requires $O(|E|)$ arithmetic
  operations. We observe
  that $y_{ik}=m_k$ if $i=j_k$ and $y_{i_k}=M_k$ otherwise.
  Hence, all the $y_{ik}$ with $(i,k)\in E$ can be computed in $O(|E|)$ arithmetic operations.  Finally,
  the $T_i(x)$ are obtained from the $y_{ik}$ in $O(|E|)$ arithmetic operations.
\end{proof}
We set:
  \[
  W:= \max_{v\in \cV}\|v \|_H \enspace .
  \]
\begin{lemma}\label{lem-bound}
  Suppose that $\cV$ is finite, then any finite eigenvector $u$ of $T$
  satisfies $\|u\|_H\leq W$.
\end{lemma}
\begin{proof}
  By definition of $W$, we have $v\in B_H(0,W)$
  for all $v\in \cV$, and since $B_H(0,W)$ is stable by tropical
  linear combinations, we get $\Sp(\cV)\subset B_H(0,W)$. Moreover,
  by~\Cref{lem-ball}, $B_H(u,-\rho(T))\subset \Sp(\cV)$.
  Hence $u\in B_H(0,W)$, meaning that $\|u\|_H\leq W$.
\end{proof}
\begin{remark}\label{rk-notfinite}
  There are situations (\Cref{sec-appendix})
  in which although some vectors of $\cV$ have infinite entries,
it is still the case that $T$ has
a finite eigenvector. Then, we may still show
that there exists a finite eigenvector
with not too large entries.
To see this, we need to replace
$W$ by $W':=\max_{k\in[p]} \delta(V_{\cdot,k})$,
  where $\delta(V_{\cdot,k})= \max_{i\in [n],(i,k)\in E}V_{ik}-\min_{j\in[n](j,k)\in E}V_{jk}$. We can always choose
such an eigenvector $u$ in such a way that $\|u\|_H = O(n W')$,
by appealing to a Blackwell optimality argument, using
the proof method of~\cite[Lemma~8.51]{skomra:tel-01958741}
(details are left to the reader).
Note that in the special case in which
$V$ has finite entries, the bound on $\|u\|_H$ is improved
by a factor $n$.
\end{remark}
\begin{corollary}[Approximate optimality certificate]\label{cor-km}
  Suppose that $\cV\subset \R^n$ is of cardinality $p$. Then, the projective
  Krasnoselskii-Mann iteration returns in a number of arithmetic
  operations $O(npW/\epsilon^2)$
  a vector $u\in \R^n$ such that $-u$ is both the center of a ball
  of radius $-\rho(T)-\epsilon$ included in $\Sp(\cV)$ and the
  apex of a regression hyperplane, $\cH_{u}$, such that $\dist_H(\Sp(\cV),\cH_u)
  \leq -\rho(T) + \epsilon$.
\end{corollary}
\begin{proof}
  By~\cite[Coro.~13]{stott2020} and~\Cref{lem-bound},
  after $k=O(\lceil W/\epsilon^2\rceil)$ iterations,
  we end up with a vector $u:=v^k$ which satisfies
  $\|T(u)-u\|_H\leq \epsilon$. Moreover, by~\Cref{prop-eval},
  each iteration
  requires $O(np)$ arithmetic operations.
  Setting
  $\underline{\lambda}:= \operatorname{bot}(T(u)-u)$,
where $\operatorname{bot}(x):=\min_i x_i$,
we deduce that $\underline{\lambda} +u\leq
T(u)\leq \underline{\lambda}(T) + \epsilon+u $, which, by~\Cref{th-duality}, entails that $\rho(T)\leq \underline{\lambda}(T)+\epsilon$.
Then, by~\Cref{th-strong-duality},
 $B(-u,-\rho(T)-\epsilon)
  \subset \Col(V)$. The proof that 
  $\dist_H(\Col(V),\cH_u)
  \leq -\rho(T) + \epsilon$ is dual. 
  \end{proof}
The following result shows that the factor in $1/\epsilon^2$ can be replaced
by $1/\epsilon$ if we look separately for the center
of a Hilbert's ball included in $\Sp(\cV)$ and
for the apex of an approximate tropical linear regression
hyperplane (in~\Cref{cor-km}, the apex and the center coincide).

\begin{corollary}
  \label{th-vi}
  Suppose that $\cV\subset \R^n$ is of cardinality $p$.
  Then, an $\epsilon$-approximation of the inner radius of $\Col(V)$,
  as well as vectors $v,z\in\R^n$ satisfying
  $B_H(v,\inradius(\cV)-\epsilon) \subset \cV$
  and $\dist_H(\cV,\cH_z) \leq \inradius(\cV)+\epsilon$
  can be obtained in $O(npW/\epsilon)$ arithmetic operations.
\end{corollary}

\begin{proof}[Proof of~\Cref{th-vi}]
We now rely on the value
iteration approach of~\cite{mtns,skomra:tel-01958741}. The latter
computes the sequence given by $v^0=0$, $v^k:=T(v^{k-1})$,
together with the numbers $\bar{\lambda}^k:= \max_{i\in [n]}v^k_i$, 
$\underline{\lambda}^k:= \min_{i\in [n]}v^k_i$. The sequence
$v^k$ generally does not converge, even up to an additive
constant. So, we rely on the following
``regularized'' sequence~\cite{sgjg04}, 
\begin{align}
  \label{e-trick}
w^k:= \inf(v^0,v^1-\bar{\lambda}^k/k,\dots,v^{k-1}-\bar{\lambda}^k (k-1)/k) \enspace .
\end{align}
Lemma~8.18 of~\cite{skomra:tel-01958741} entails
that $\rho(T)$ satisfies $\underline{\lambda}^k/k\leq \rho(T)\leq\bar{\lambda}^k/k$ with $\bar{\lambda}^k/k-\rho(T)\leq \|u\|_H/k$
and $\rho(T)-\underline{\lambda}^k/k\leq \|u\|_H/k$,
where
$u\in\R^n$ is an arbitrary finite eigenvector of $T$. Hence,
it suffices to execute the algorithm up to the iteration $k:= \lceil W/\epsilon \rceil$ to make sure that $\bar{\lambda}^k \leq \rho(T) + \epsilon$ and $\underline{\lambda}^k\geq \rho(T)-\epsilon$. 
Moreover, Lemma~2 of~\cite{sgjg04} entails that $T(w^k) \leq \bar{\lambda}^k + w^k$. Hence, by~\Cref{lem-ball},
$-w^k$ is the center of a Hilbert's ball of radius $-\bar{\lambda}^k$ included in $\Sp(\cV)$. The construction of the apex of an approximate
optimal regression hyperplane uses a dual argument, replacing
inf by sup in~\eqref{e-trick}.
\end{proof}
\begin{remark}
  The conclusions of \Cref{cor-km} and \Cref{th-vi}
  can be extended to the situation
  in which some vectors of $\cV$
  have infinite entries, provided $T$
    has a finite eigenvector. Using~\Cref{rk-notfinite}, we
    need to replace $W$ by $W'n$ in the bounds of~\Cref{cor-km} and \Cref{th-vi}.
\end{remark}

\section{Illustration: inferring hidden information from equilibria in repeated invitations to tenders}\label{sec-example}

We now illustrate our results
on an example from auction theory, in
which tropical linear regression
allows one to identify secret information
from the observation of prices offered in repeated invitations
to tenders (ITT).

\subsection{Auction model with hidden preference factors}
We suppose a public
decision maker chooses the best offer made by the firms
responding to ITT. In accordance
with market regulations, see e.g.~\cite[Art.~R.2152-7]{codedelacommandepublique}, the best offer is not nessarily the one with the lowest price:
other factors, like technical quality, 
respect of environment, of social impact,
can also be taken into account.  In the presence
of corruption, decisions may be also influenced
by bribes.

We assume that this ITT is done repeatedly for a similar service or product each time and in front of the same local firms. We label the firms by $1,2, \cdots,n$, and we suppose that we have a history of $q$ ITTs with the prices offered by each firm, that are revealed by the decision maker, after
having made her choice.

More precisely, we denote the price offered by firm $i\in [n]$ for the ITT number $j\in [q]$ by $p_{ij}$.
We assume that the decision maker
has a non public preference factor $f_i>0$ for each firm $i$, and that
she selects the firm of index $i$ minimizing the expression:
\begin{equation}\label{e-decision-maker}
\min_{i\in[n]} p_{ij} f_i^{-1} \enspace .
\end{equation}
In this way, the decision maker considers that
for a requested price of $p_{ij}$, the final cost
to be taken into account is $p_{ij}f_i^{-1}$, where $f_i^{-1}\geq 1$ is a proportional penalty depending on her estimate $f_i$
of the technical, environmental, or social
quality of the firm (the larger $f_i$, the better its quality).

The same model
applies to the
situation in which $f_i^{-1} = 1 -\alpha_i \beta$ for some
$0\leq \alpha_i\leq 1$ and $0\leq \beta< 1$. Now, $\alpha_i$
may be interpreted as a proportional bribe: the firm
promises to secretly give back $\alpha_i p_{ij}$ to the decision
maker if its offer is accepted, and the parameter
$\beta$ measures how sensitive is the decision
maker to bribery ($\beta=0$ corresponds to a totally
honest decision maker, and $\beta=1^-$ to a totally
dishonest one). This is a variant of the
classical {\em first-price sealed-bid auction}~\cite{krishna},
incorporating the secret preference.

We suppose that the same firms answer in a recurrent
manner to invitations from the same decision maker, and
that the factors $f_i$ secretly attached to each firm are kept constant.
Then we expect
that the prices to be offered to constitute an {\em equilibrium},
meaning that for each invitation $j\in[q]$, the minimum
$\min_{i\in[q]}p_{ij}f_i^{-1}$ is achieved twice at least. Indeed, if
the firm $i$ that wins the invitation offers a price $p_{ij}$ such that
$p_{ij}f_i^{-1}$ is strictly smaller than $p_{kj}f_k^{-1}$ for all $k\in [n]\setminus\{i\}$,
it may offer a higher price and still win the offer, so,
in the long run, if an invitation of the same type
is made recurrently, the firm will adapt its offer.

This can be modeled in terms of membership to a tropical hyperplane.
We put $V_{ij}=-\log (p_{ij})$ and $a_i = \log f_i$, so that
the decision maker selects the firm of index
achieving the maximum in 
\begin{equation}\label{e-bribery}
\max_{i \in [n]} (V_{ij} + a_i) \enspace.
\end{equation}
Assuming the prices $p_{ij}$ are observed, our goal is to infer
the secret information $f_i$, i.e.\ the preference
factor for firm $i$, or the bribe
offered by this firm.

We first suppose %
that for each invitation, the identity of the firm that wins the contract is not known to us.
We want to infer the hidden information $f=(f_i)_{i\in [n]}$. %
So, we look for a %
tropical hyperplane $\cH_b$ that is the best regression of the set $\cV$ formed by the points $(V_{\cdot j})_{j\in[q]}$ following the analysis of~\Cref{sec-strong},
i..e, we solve a problem
of the form~\eqref{e-regression}.
Following~\Cref{th-strong-duality}, we solve this problem by computing
a super-eigenvector $b\in \R^n$ of $T$, i.e. such that $T(b)\geq \rho(T) + b$, where the operator $T$ is given by~\eqref{e-def-T}.

We note that the decision maker cares only about the {\em relative} preference factors between the firms, in the sense that if all the preference factors $f_i$, $i\in [n]$ are multiplied by the same positive constant, the choices of the decision maker will not change. Therefore, we can suppose without loss of generality that $\max_{j\in[n]}f_j=1$, or equivalently, $\max_{j\in[n]}a_j =0$.

\subsection{Numerical instance and experiments}
In the following toy example, we take $n=3$ firms, and a history of $q=6$ ITTs. 
We suppose that the decision maker attributes to the firms
the preference factors $f=(1, 0.8 , 0.6)$,
and we take $\forall i\in [3], a_i=\log(f_i)$.

We generated the matrix $V_{ij}$ and the prices $p_{ij}=\exp(-V_{ij})$ by the
following structured probabilistic model. We consider six types of products with prices of different order of magnitude. In~\Cref{tab-example} the reference prices of these products are $P=(1,3,9,25,70,130)$.
For each $j\in[6]$, we draw entries $A_{ij}, i\in [3]$ randomly in the interval $K_{ij}=[0.9\times P_j f_i, 1.1\times P_j f_i]$ following a log-uniform law, i.e.\ equal to the exponential of a variable generated uniformly on the logarithm of the interval $K_{ij}$. We choose the log-uniform law because it's in adequacy with Benford's law that is observed in real-life price
instances.  Then, we take $B_{ij}=-\log(A_{ij})$, and we project each column $B_{\cdot j}$ into the tropical hyperplane $\cH_{a}$, to get a vector $C_{\cdot j}$, such that for a given $i\in \argmax_{k\in [n]}(B_{kj}+a_k)$, we take $C_{ij}=\max_{k\neq i}(B_{kj}+a_k) - a_i$ and we take for all $k\neq i$, $C_{kj}=B_{kj}$. 
Now the columns $C_{\cdot j}$ belong to the tropical hyperplane $\cH_a$. To model the inefficiency of the market, we perturb these columns by taking $V_{ij}=C_{ij}+\delta_{ij}$, with $\delta_{ij}$ generated randomly uniformly in $[-\delta,\delta]$, with 
$\delta=0.05$.
Then the prices are given by $p_{ij}=\exp(-V_{ij})$.

To solve our example, we used the projective Krasnoselskii-Mann iteration
described in~\Cref{sec-algo},
with a damping parameter $\gamma = 1/2$.
We take $b=v^N$ that gives the approximation of the preference factors by tropical linear regression: $f^{\textrm{reg}}_i=\exp(b_i)$, $i \in[n]$.

We define the error of the approximation $e$ as the ratio between the Hilbert's distance of the set $\cV$ to the hyperplane $\cH_b$, which measures the ``distance to equilibrium" in this market, and the maximal absolute value of the logarithm of the Hilbert's seminorms of the price vectors $(p_{\cdot j})_{j\in [q]}$:
\[
e := \frac{\dist_H(\cV, \cH_b)}{\max_{j\in[q]} |\log(\|p_{\cdot j}\|_H) |} \enspace.
\]
The following~\Cref{tab-example} shows the preference factors $f_i$, the prices $p_{ij}$ generated with this model and for each invitation we underlined the price of the firm wining that invitation in the sense of achieving the minimum in~\eqref{e-decision-maker}. 
\Cref{tab-example} shows also the prediction $f^{\textrm{reg}}$ of the preference factors that we find by tropical linear regression.

\begin{figure}[htbp]
\center
\includegraphics[scale=0.6]{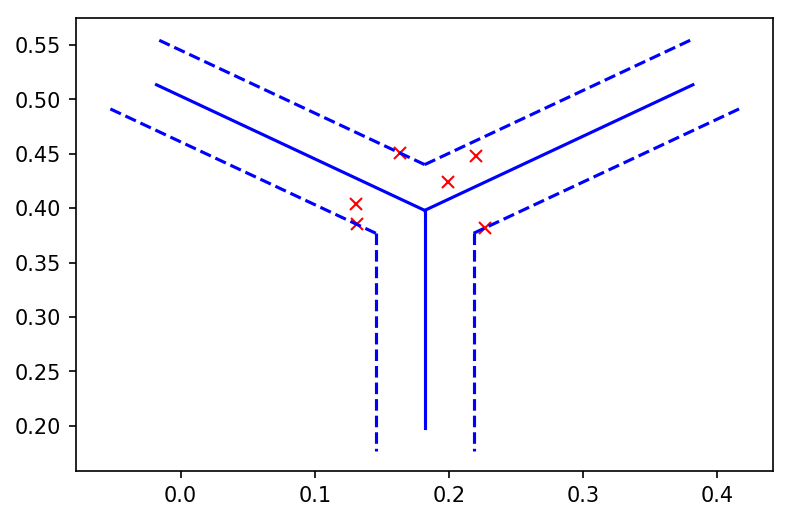}
\caption{The points $(V_{\cdot j})_{j\in [6]}$ in the projective space $\proj(\rmax)^3$, with the tropical hyperplane $\cH_b$ (in blue solid lines) and the points of the space that are at distance equal to $\dist_H(\cV, \cH_b)$ from $\cH_b$ (in blue dashed lines).}
\label{fig-bribery6}
\end{figure}

In this example, we set a target accuracy of $\epsilon=10^{-8}$, and we get that the number of iterations $N$ needed to get $\|T(v^N)-v^N\|_H\leq \epsilon$ is $N=25$. 
By setting $b=v^N$, we have $\dist_H(\cV, \cH_b)=4.21 \times 10^{-2}$ and $\max_{j\in[q]} |\log(\|p_{\cdot j}\|_H) | = 3.84$, and this leads to an error equal to  $e=1.09 \times 10^{-2}$. \Cref{fig-bribery6} shows the points $(V_{\cdot j})_{j\in [6]}$ in the projective space $\proj(\rmax)^3$, with the tropical hyperplane $\cH_b$ (in blue solid lines) and the points of the space that are at distance equal to $\dist_H(\cV, \cH_b)$ from $\cH_b$ (in blue dashed lines). \Cref{fig-bribery6} shows in particular the existence of a witness point in each on the three sectors associated to the tropical hyperplane $\cH_b$.

\begin{table}[H]
	\begin{tabular}{ c | c | c | c | c | c | c | c | c }
		& individual houses & social housing & school & road & stadium & bridge & $f$  & $f^{\textrm{reg}}$ \\
		\hline
		Firm $1$
		& 1.02 &  3.21 & \underline{8.72} & 26.2 & 69.8 & \underline{123} & 1 &  1\\
		\hline
		Firm $2$
		&  0.81 & 2.65 & 7.49 & 20.3 & \underline{53.8} & 106 & 0.8 & 0.81\\
		\hline
		Firm $3$
		&  \underline{0.6} &  \underline{1.86} &  5.5 & \underline{14.7} & 41.8 & 76 & 0.6 &  0.605\\
	\end{tabular}
	\caption{Prices proposed by firms in million euros, the vector of preference factors $f$ and its estimation by tropical linear regression $f^{\textrm{reg}}$ based on the observation of the prices.}
\label{tab-example}
\end{table}

\begin{figure}[htbp]
\center
\includegraphics[scale=0.6]{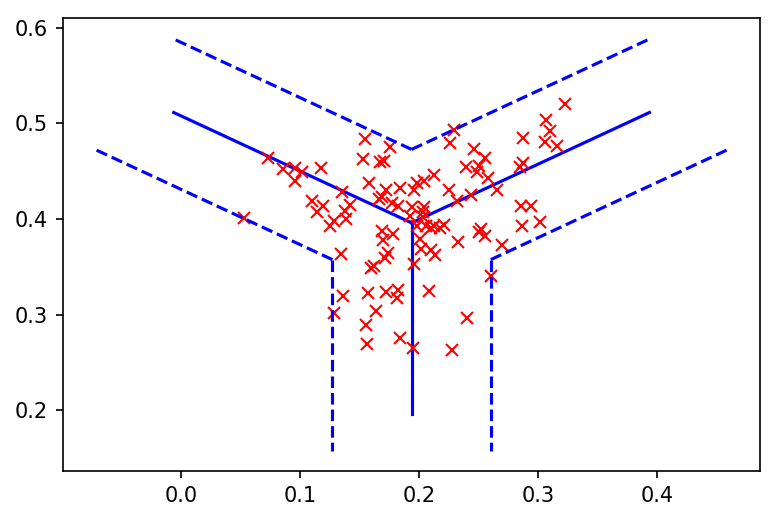}
\caption{The points $(V_{\cdot j})_{j\in [100]}$ in the projective space $\proj(\rmax)^3$, with the tropical hyperplane $\cH_b$ (in blue solid lines) and the points of the space that are at distance equal to $\dist_H(\cV, \cH_b)$ from $\cH_b$ (in blue dashed lines).}
\label{fig-bribery100}
\end{figure}

Now we consider a similar example still with $n=3$ firms, but with $q=100$ invitations to tenders. We use the same generation model, the reference prices $P_j, j\in [100]$, being generated randomly following a log-uniform law on the interval $[1,100]$. 
We set a target accuracy of $\epsilon=10^{-8}$, and we get that the number of iterations $N$ needed to get $\|T(v^N)-v^N\|_H\leq \epsilon$ is $N=24$. 
By setting $b=v^N$, we have $\dist_H(\cV, \cH_b)=7.69 \times 10^{-2}$ and $\max_{j\in[q]} |\log(\|p_{\cdot j}\|_H) | = 3.72$, and this leads to an error equal to  $e=2.06 \times 10^{-2}$, and the approximation of the preference factors that we obtain is $f^{\textrm{reg}}=(1, 0.7994, 0.6018)$.
\Cref{fig-bribery100} shows the points $(V_{\cdot j})_{j\in [100]}$ and the approximation hyperplane $\cH_b$ obtained in this case with a history of $q=100$ invitations. We observe also that we have at least a witness point in each sector defined by the tropical hyperplane $\cH_b$. 

\subsection{Example of regression with types -- in which the identities of the winners of the invitations are known}
We now suppose the decision maker makes public not only the bid
prices that were offered to her, but also the identities
of the firms that won the different invitations $j\in [q]$.
Then, we can write the set of points $\cV$ as a disjoint union
$\cV=\cup_{\ell\in [n]}\cV_\ell$, where $\cV_\ell$ is the set of invitations won by firm $\ell$.
  This information can be exploited through
  the {\em typed} tropical linear regression of~\Cref{subsec-typed}.
  Indeed, if $v=V_{\cdot j}\in \cV_\ell$, and if the market
  is ``at equilibria'', we know not only that the maximum
  $\max_{i \in [n]} (V_{ij} + a_i)$ is achieved twice, but
  that it must be achieved by the firm that won the invitation,
  i.e., $i=\ell$. Thus, the vector $v\in \cV_\ell$
  should be close to the {\em signed} tropical hyperplane
  $\cH_a^\ell$, a finer condition than
  being close to $\cH_a\supset \cH_a^\ell$.
  So, to infer the vector $a$, we now solve
  the typed regression problem~\eqref{eq-reg-types},
  instead of the untyped problem~\eqref{e-regression}.
Following~\Cref{th-typed}, we are looking for a super-eigenvector $b$ such that $\Ttype (b) \geq \rho(\Ttype)+b$, where the operator $\Ttype$ is given by~\eqref{e-T-types}.

We use the same two examples above, and we generate the information of the firm winning each contract $j\in [q]$ by using the information $f$ known by the decision maker. We construct the sets $\cV_i$ and the operator $\Ttype$, and we find a super-eigenvector of $\Ttype$ by using the projective Krasnoselkii-Mann value iteration algorithm described in~\Cref{sec-algo}.

After doing the numerical experiments, we find
that, the apex $b$ found by typed tropical linear regression, taking advantage of the knowledge of which firm won each invitation, is the same as the one found above by tropical linear regression, for both examples with $q=6$ and $q=100$. Hence, here, the additional information provided
by the identity of the winners did not help to improve
the inference of hidden preferences, by comparison
with the basic model in which only the history of the bid
prices is used.

\section{Concluding remarks}

We solved the tropical linear regression problem,
when the metric is of {\em sup-norm} type, and for tropical
linear spaces of codimension $1$ (tropical hyperplanes),
but for a configuration of points of arbitrary cardinality.

Several open problems related to the present work arise when
changing either the class of metrics or of tropical spaces.

For instance, we may replace Hilbert's metric by
the $L_p$-projective metric, i.e., the metric obtained
by modding out the $L_p$ normed space $\R^n$ by the action
of additive constants, or by replacing the Hausdorff distance
in~\eqref{def-hausdorff} by a $L_p$ type distance,
for $p\in [1,\infty)$.  Approaches based on mixed
linear programming, or on local descent, have been proposed
in~\cite{yoshidatropPCAphylo,pageyoshida,hook}
in some specific cases.

Another generalization consists in replacing
hyperplanes by tropical linear spaces
of a codimension not necessarily $1$.
Recall that the {\em tropical Grassmannian}
$\operatorname{Gr}^{\operatorname{trop}}_{k,n}$
can be defined
as the image by a non-archimedean valuation of the
Grassmannian $\operatorname{Gr}_{k,n}(\mathbb{K})$
over an (algebraically closed) non-archimedean
field, under the Pl\"ucker embedding,
see~\cite{SpeyerSturmfels04,rinconfink}.
In this way, an element of
$\operatorname{Gr}^{\operatorname{trop}}_{k,n}$
is represented by its {\em tropical Pl\"ucker coordinates}
$p=(p_I)\in (\R\cup\{-\infty\})^{n\choose k}$.
This vector yields a {\em tropical linear space $L(p)$}, defined
by
\[
L(p)  = \bigcap_{I} \{ x\in (\rmax)^n\mid
\max_{i\in I} (p_{I\setminus \{i\}}) +x_j)\text{ is achieved at least twice}\}\enspace ,
  \]
where the minimum is taken over all subsets of $[n]$ of cardinality $k+1$.
When $k=n-1$, $V(p)$ is a tropical hyperplane. Hence,
a general version of tropical linear regression problem
can be written as
\begin{align}
\min_{p\in \operatorname{Gr}^{\operatorname{trop}}_{k,n}} \max_{v\in \mathcal{V}}
  \min_{x\in L(v)} \|v-x\|_H \enspace .\label{e-Tropgen}
\end{align}
We solved here this problem when $k=n-1$. When
$k=1$, $L(p)$ is reduced to a single point, and it
is not difficult to see that~\eqref{e-Tropgen} reduces to a
linear program. We leave it as an open question to solve
this problem when $1<k<n-1$. The same problem may
be considered when $p$ is a valuated
matroid, or when it is inside the image
of the Stiefel map~\cite{rinconfink},
meaning that $p$ is given by the maximal tropical
minors of a matrix. A version
of the latter problem
(with a $L_1$-type error) is considered in~\cite{yoshidatropPCAphylo}.
One may also replace the linear space $L(p)$ by the {\em column
  space} of a tropical matrix $A$, which boils
down to finding a best approximation by a tropical polyhedral cone
with a fixed number of vertices, see~\cite{hook,pageyoshida}.

   \section*{Acknowledgment}
   We thank Quentin Jacquet and Marie Laclau for helpful comments
   and references, and also Nicolas Vieille for helpful references.

\appendix

\section{Dominions of the two players and existence of a finite eigenvector}\label{sec-appendix}

The strongest form of strong duality (\Cref{th-witness}), with
the existence of {\em witness points}, is valid whenever the
Shapley operator $T$ in~\Cref{e-def-T} has a finite eigenvector.
In this appendix, we provide a sufficient condition
for the existence of this eigenvector,
which is less demanding than the condition of~\Cref{prop-finite}
(requiring $V$ to have only finite entries).

We recall that the operator $T$ represents a game $\Gamma$ with two players Min and Max, such that when we are at state $i$, player Min plays first by choosing a column $k\in [p]$ such that $(i,k)\in E$, then player Max chooses a state $j\in [n]$ such that $j\neq i$ and  $(j,k)\in E$. Moreover, policies can be defined using \eqref{setAforT}. The game $\Gamma$ is played repeatedly starting from a given initial position.

We call {\em dominion} of one player a nonempty subset of states $I\subset [n]$ such that from any initial position in $I$, that player can force the state to remain in $I$ at each stage of the repeated game, whatever actions the other player chooses. This means that there exist a policy of that player such that for any strategy of the other player, a trajectory of the game starting in $I$ is such that the states visited by Min are all contained in $I$.
The next result, which follows from a more general result (which applies to
arbitrary Shapley operators) relates the lack of disjoint dominions of the two players
with the existence of a finite eigenvector of a polyhedral Shapley operator.

\begin{theorem}[Corollary of Thm.\ 1.2 of~\cite{akian2018game}]\label{thm-dominions}
The following assertions are equivalent:
\begin{enumerate}
\item The two players do not have disjoint dominions in the game $\Gamma$;
\item For all $r\in \R^n$, the operator $r+T$ has a finite eigenvector.
\end{enumerate}
\end{theorem}
Deciding the existence of disjoint dominions for (deterministic) mean payoff
games is equivalent to deciding the existence of a non-trivial fixed point
of a monotone Boolean function, which is a NP-complete problem, see the discussion in~\cite{hochart2015}. 
However, we next show that for the restricted
class of games associated to the  Shapley operator $T_V$,
this problem can be solved in polynomial time.

We make the following assumption, which
is required for the operator $T$ to send $\R^n$ to $\R^n$,
and a fortiori, to have a finite eigenvector.
\begin{assumption}\label{assump-two-finite-entries}
Each column of the matrix $V$ contains at least two finite entries.
\end{assumption}

\begin{proposition}\label{prop-dominions}
Suppose that \Cref{assum-finrow} and \Cref{assump-two-finite-entries} hold.
Then, the following assertions are equivalent:
\begin{enumerate}
\item\label{prop-dominions1} There are disjoint dominions for the two players in the game $\Gamma$; 
\item\label{prop-dominions2}  There exist nonempty subsets $I, J$ of $[n]$,  such that $I\cup J = [n]$, $I \cap J = \emptyset$, some columns of $V$ have support included in $I$, and the other columns of $V$ have at least two finite entries in $J$.
\item\label{prop-dominions3}  There exists a subset $K$ of $[p]$, such that $K\neq\emptyset$ and $K\neq [p]$, and such that if we denote by $I_K$ the union of supports of the columns of $V$ in $K$, then all the columns not in $K$ have at least two finite entries that are outside $I_K$. In this case, $I_K$ together with its complement $[n]\setminus I_K$ constitute disjoint dominions of players  Min and Max, respectively.
\end{enumerate}
\end{proposition}

\begin{proof}%
We verify first that the assertion \eqref{prop-dominions2} and the first part of assertion \eqref{prop-dominions3} are equivalent. Indeed, it is straightforward that \eqref{prop-dominions3} implies \eqref{prop-dominions2} by taking $I=I_K$ and $J=[n]\setminus I_K$. Now, if \eqref{prop-dominions2} is true, we take $K=\{k\in [p] \mid \supp V_{\cdot k} \subset I\}$, so $I_K=\cup_{k\in K}\supp V_{\cdot k}\subset I$, and each column $k\not\in K$ has at least two finite entries in $J$, i.e. outside $I_K$.

Now we suppose that assertion \eqref{prop-dominions1} is satisfied, i.e.\ there are disjoint dominions $I$ and $J$ respectively for Player Min and Player Max. 
Let us show that this implies \eqref{prop-dominions2}.
The set $J$ is a dominion for Player Max, then there exists a policy $\tau$ for Max, such that, for all $i\in J$, for any possible action $(i,k)\in E$ of Player Min, the policy $\tau$ sends the state in $J$, that is $\tau((i,k))\in J$.
Since a policy for Max is a map from $E$ to $[n]$ such that $j=\tau((i,k))$ satisfies $j\neq i$ and $(j,k)\in E$, this implies that, for all $(i,k)\in E$ with $i\in J$, there exists $(j,k)\in E$ with $j\in J$ such that $j\neq i$. Therefore, for all $k\in [p]$, 	$\supp V_{\cdot k} \cap J$ is either empty or it contains at least two elements. We take $I'=[n]\setminus J\supset I\neq \emptyset$, then the sets $I',J$ satisfy the assertion \eqref{prop-dominions2}.

Now, we suppose that the first part of assertion~\eqref{prop-dominions3} is true, and show that $I_K$ and $J=[n]\setminus I_K$ are disjoint dominions of players Min and Max respectively, which will imply \eqref{prop-dominions1}. Indeed, if $i \in I_K$, then there exists $k\in K$, such that $i\in \supp V_{\cdot k} \subset I_K$. Let us consider a policy $\sigma$ of Min such that  if $i\in I_K$ then $\sigma(i)=(i,k)$ with $k\in K$. Then, if $i\in I_K$, and if Min  plays the action $(i,k)=\sigma(i)$,  for any possible action of player Max (which exists by~\Cref{assump-two-finite-entries}), that is a choice of $j\in \supp V_{\cdot k} $ such that $j\neq i$, we have 
$j\in I_K$. This shows that $I_K$ is a dominion of Player Min. Now, let $i\in J$, for any action $(i,k)\in E$ of Min (which exists by~\Cref{assum-finrow}), we have $k\in K$, since $i\in \supp V_{\cdot k}$ and $i \not\in I_K$, so by \eqref{prop-dominions3}, there exits $j \in \supp V_{\cdot k}\setminus I_K$, with $j\neq i$. So $j\in J$ and $j$ is a possible action of Max when the game is in state $(i,k)$.
Considering the policy $\tau$ for Max, such that $\tau(i,k)=j$ for $i,k,j$    as before, we get that the set $J$ is a dominion of Player Max.
\end{proof}

From the proof of~\Cref{prop-dominions}, we deduce in a straightforward manner the following observation, which will be used in \Cref{algo-dominions}.
Note that in the present setting (deterministic mean payoff games),
if $I,J$ are disjoint dominions of the two players, then
$[n]\setminus J$ and $J$ are
also dominions of the two players, hence we shall restrict
our search to disjoint dominions that constitute partitions of
$[n]$.

\begin{lemma}\label{lem-dom-min}

If $D^{\operatorname{Min}}, D^{\operatorname{Max}} \subset [n]$ are disjoint dominions of players Min and Max respectively, that constitute a partition of $[n]$, and $K$ is a subset of columns of $V$ such that the set $S=\cup_{k\in K}\supp V_{\cdot k}$ 
satisfies $S\subset D^{\operatorname{Min}}$, then for each column $k\not\in K$ that has only one finite entry $i$ outside of $S$, we have $S\cup \{i\} \subset D^{\operatorname{Min}}$.\hfill\qed
\end{lemma}

\begin{algorithm}[htbp]
	\begin{algorithmic}[1]
	  \For{$k \in [p]$}
	\State $K\leftarrow\{k\}$ 
        \State $S\leftarrow$ the support of column $k$ of $V$
        \State Declare $S$ to be {\sc augmented} (Boolean flag)
        \While{$S$ is declared as {\sc augmented}}
        \State Declare $S$ not to be {\sc augmented}
        \State Declare all the elements of $[p]\setminus K$ to be {\sc unscanned} (Boolean flags)
        \While{$([p]\setminus K)$ contains an {\sc unscanned} element}
        \State $\ell\leftarrow$ smallest {\sc unscanned} element of $[p]\setminus K$, declare $\ell$ to be {\sc scanned}
       \State $S_\ell\leftarrow \{i\in [n]\setminus S \mid V_{i\ell} \text{ is finite}\}$\label{linedefSl} %
       \If{$|S_\ell|= 1$} $\triangleright$ {\em column $\ell$ of $V$ has precisely one finite entry outside $S$}
        \State $K\leftarrow K\cup\{\ell\}$,
	$S \leftarrow S \cup S_\ell$\label{lineupS}        
        \State Declare $S$ to be {\sc augmented}
	\EndIf
	\If{$|S_\ell|= 0$} $\triangleright$ {\em column $\ell$ of $V$ has no finite entries outside $S$}
        \State $K\leftarrow K\cup\{\ell\}$\label{lineupK2}
	\EndIf
        \State $\triangleright$ $S$ is the union of supports of the columns of $K$ \label{line-invariant}
        \EndWhile
	\EndWhile
	\If{$S\neq [n]$}
	 \Return\label{line-return} $S$ and $[n]\setminus S$ which are disjoint dominions of Players Min and Max respectively
	\EndIf
	\EndFor
	\State There are no disjoint dominions
	\end{algorithmic}
	\caption{Detecting dominions in the game arising from the tropical linear regression problem, for an input matrix $V\in (\rmax)^{n\times p}$.}
\label{algo-dominions}
\end{algorithm}

\begin{theorem}\label{th-correct}
  \Cref{algo-dominions}, which decides the existence
  of disjoint dominions in the game $\Gamma$ associated
  to a matrix $V\in (\rmax)^{n\times p}$,
  is correct, and it makes
  $O(n^2 p^2)$ arithmetic operations.
\end{theorem}

\begin{proof}
  The algorithm looks for a set of columns $\bar{K}$ satisfying
  the last statement of~\Cref{prop-dominions}. Since
  the set $\bar{K}$ is required to be nonempty, it suffices
  for each $k\in [p]$, to verify whether
  there is such a set $\bar{K}\ni k$ (for loop of the algorithm).

  We next show that the algorithm admits the following
  invariants.
  \begin{enumerate}
  \item At Line~\ref{line-invariant},
     $S$ is the union of supports of the columns of $K$.
  \item If there is a subset $\bar{K}\ni k$ satisfying
    the last statement of~\Cref{prop-dominions}, with associated
 then at line~\ref{line-invariant} of the algorithm,
    the set $K$ satisfies $K\subset \bar{K}$ and the set
    $S$ satisfies $S\subset D^{\operatorname{Min}}$.
    \end{enumerate}
  The first invariant is enforced by lines~\ref{linedefSl},~\ref{lineupS} and~\ref{lineupK2}. We prove that the loop invariant at line~\ref{line-invariant}
  holds by induction
  on the cardinality of $S$. Let us assume that the condition
  of the first ``if'', i.e., $|S_\ell|=1$ is satisfied.
  Then, by~\Cref{lem-dom-min}, and by the induction assumption,
  we must have $S\cup \{i\}\subset D^{\operatorname{Min}}$.
Moreover, the last statement of~\Cref{prop-dominions}
entails that $K\cup \{\ell\}\subset \bar{K}$, and so,
the loop invariant is valid in this case.
Moreover, if the condition
of the second ``if'', i.e., $|S_\ell|=0$ is satisfied,
then, the second invariant is still valid. This
shows that the loop invariant
is always valid.

At the exit of the outer while loop, at line~\ref{line-return}, we have by construction
that every column of $V$ with index outside $K$
has at least two finite entries outside $S$.
Then, by the last statement of~\Cref{prop-dominions},
if $S\neq [n]$, $S$ and $[n]\setminus S$ provide
disjoint dominions of Players Min and Max, whereas
if $S=[n]$, there are no dominions arising
from a set $\bar{K}\ni k$. This shows the correctness
of the algorithm.

Each iteration of the inner ``while'' loop makes $O(n)$ arithmetic operations, and every outer ``while loop'' executes the inner while loop $O(p)$ times. Moreover, the number of outer ``while loop'' iterations is at most $n-1$. Finally, we have at most $p$ iterations in the ``for'' loop, which leads to a complexity bound of $O(n^2 p^2)$ arithmetic operations for the algorithm.
\end{proof}
We call {\em Boolean pattern} of the matrix $V\in (\rmax)^{n \times p}$ the matrix with entries in $\{0,-\infty\}$, obtained
by replacing each finite entry of $V$ by $0$.
\Cref{thm-dominions} provides a
sufficient condition involving the Boolean pattern on $V$,
which guarantees that for all matrices $V$ with this
pattern, the operator $T$ admits a finite eigenvector.
This condition is not necessary.
Consider the following Boolean pattern:
\begin{equation}\label{exp-pattern}
  \small
  \begin{pmatrix}
    0_{3,2} & 0_{3,2}\\
    (-\infty)_{3,2}& 0_{3,2}
    \end{pmatrix}
    \enspace ,                  %
\end{equation}
where for $\alpha\in \{0,-\infty\}$,
$\alpha_{p,q}$ denotes the $p\times q$ matrix
with entries identically equal to $\alpha$.
\begin{proposition}
  If $V$ is a matrix with Boolean pattern~\eqref{exp-pattern},
  then, the operator $T$ has a finite eigenvector,
  but the associated game admits disjoint dominions.
\end{proposition}

\begin{proof}
First we have that the set $K=\{1,2\}$ satisfies the condition \eqref{prop-dominions3} of \Cref{prop-dominions}, and from the proof of~\Cref{prop-dominions}, we have that the sets $I=\{1,2,3\}$ and $J=\{4,5,6\}$ are disjoint dominions of Player Min and Player Max respectively.%

To show that $T$ has a finite eigenvector, by~\Cref{prop-finiteeig},
it suffices to check that $\chi(T)=0$. The inequality
$\chi(T)\leq 0$ follows from~\Cref{rem-nonnegative}. 
We next show that $\chi(T)\geq 0$.

If the game starts from a state $i\in \{4,5,6\}$, Player Min
must choose the next state to be a pair $(i,k)$ with
$k\in \{3,4\}$, and Player Max can respond by choosing
the next state $j$ to belong to $\{4,5,6\}$.
So, Player Max can force Min
to play the same game as the one defined
by the submatrix $X:=(V_{ij})_{i\in \{4,5,6\},j\in \{3,4\}}$. Since the
matrix $X$ consists of only $2$ columns of $(\rmax)^3$,
it follows from~\Cref{cor-simplicial} that
the inner radius of $\Col(X)$ is equal
to $0$. Then by~\Cref{th-duality}, $\rho(T_X)=0$, and this
entails that Player Max can ensure a payment equal to $0$ in the original game,
so that
$\forall i\in \{4,5,6\}, \chi_i(T)=0$.

Suppose now that the initial state $i\in \{1,2,3\}$.
Since $\chi(T)=\min_{\sigma} \chi(T^\sigma)$ where the minimum is taken over the stationary policies of Player Min, it suffices to show that for any
such policy, and for $i\in \{1,2,3\}$, $\chi_i(T^\sigma)\geq 0$.
If this policy of Player Min chooses the column $3$ or $4$,
Player Max can again enforce Player Min to play the game associated
to the submatrix $X$, and then Player Max can ensure a payment $0$ as before.
Now, if the policy of Player Min does not choose the columns $3$ and $4$,
Player Max is forced to play a subgame correspinding to the the submatrix $Y:=(V_{ij})_{i\in \{1,2,3\},j\in \{1,2\}}$, and by the same reasoning as before,
we know that the value of this game is equal to $0$.
Then $\forall i\in \{1,2,3\}, \chi_i(T)=0$. 
\end{proof}
We leave it as an open question to characterize
the Boolean patterns of $V$ which guarantee that
the operator $T$ has a finite eigenvector.

\bibliography{tropicallinearregressionforarxivv2}


\bibliographystyle{alpha}
\end{document}